\documentclass[a4paper, oneside, english,reqno]{amsart}
\usepackage[utf8]{inputenc}

\usepackage[english]{babel}

\usepackage{amsmath,amsfonts,amssymb,mathrsfs,amscd,comment,latexsym,bbold}
\usepackage{amsthm, mathtools, stmaryrd}
\usepackage{textcomp, url, enumerate, latexsym, graphicx, varioref, hyperref, multirow, layout, siunitx, booktabs}
\usepackage{pdflscape}
\usepackage{verbatim}

\usepackage[matrix,arrow,curve]{xy}
\usepackage{hyperref}

\usepackage{enumitem}

\usepackage{csquotes}

\usepackage[top = 3cm, bottom = 4cm, left = 3cm, right = 3cm]{geometry}

\usepackage{microtype}

\usepackage{cleveref}
\RequirePackage[backend    = biber,
                sortcites  = true,
                giveninits = true,
                doi        = false,
                isbn       = false,
                url        = false,
                maxnames   = 6,
                style      = numeric]{biblatex}
\DeclareFieldFormat*{title}{#1}
\usepackage{biblatex}
\renewbibmacro{in:}{}

\usepackage{amsthm}

\theoremstyle{plain}
\newtheorem{theorem}{Theorem}
\newtheorem*{nonum-theorem}{Theorem}
\newtheorem{proposition}{Proposition}
\newtheorem{lemma}{Lemma}
\newtheorem{lemma-remark}{Lemma-remark}

\newtheorem{corollary}{Corollary}
\newtheorem{nonum-corollary}{Corollary}

\newtheorem*{notation*}{Notation}

\newtheoremstyle{handleNumber}{}{}{\itshape}{}{}{}{\newline}{{\bf #1} \thmnote{#3}}
\theoremstyle{handleNumber}
\newtheorem*{handnum-theorem}{Theorem}

\theoremstyle{definition}
\newtheorem{definition}{Definition}

\newtheorem*{nonum-definition}{Definition}

\newtheorem{example}{Example}

\newtheorem{conjecture}{Conjecture}[section]

\theoremstyle{remark}
\newtheorem{remark}{Remark}

\newcommand{\nocontentsline}[3]{}
\newcommand{\tocless}[2]{\bgroup\let\addcontentsline=\nocontentsline#1{#2}\egroup}

\renewcommand{\leq}{\leqslant}
\renewcommand{\geq}{\geqslant}

\newcommand{\PP}{\mathbb P}

\newcommand{\SH}{\mathbf{SH}}
\newcommand{\rmSH}{\mathrm{SH}}

\newcommand{\HH}{\mathbf{H}}
\newcommand{\rmHH}{\mathrm{HH}}

\newcommand{\Spec}{\operatorname{Spec}}

\newcommand{\Fr}{\mathrm{Fr}}
\newcommand{\fr}{\mathrm{fr}}
\newcommand{\fraff}{\mathrm{fraff}}
\newcommand{\tfr}{\mathrm{tfr}}
\newcommand{\nfr}{\mathrm{nfr}}
\newcommand{\pfr}{\mathrm{pfr}}
\newcommand{\tgfraff}{\mathrm{tgfr}}
\newcommand{\tghfraff}{\mathrm{tghfr}}

\newcommand{\lownrfraff}{\mathrm{lownr}}

\newcommand{\Aff}{\mathrm{Aff}}
\newcommand{\Sch}{\mathrm{Sch}}
\newcommand{\SchPair}{\mathrm{Sch}^\mathrm{pair}}
\newcommand{\Sm}{\mathrm{Sm}}
\newcommand{\Smat}{\mathrm{Sm}^\mathrm{cci}}
\newcommand{\SmAff}{\mathrm{SmAff}}
\newcommand{\AffSm}{\mathrm{SmAff}}
\newcommand{\Sh}{\mathrm{Shv}}
\newcommand{\Pre}{\mathrm{Pre}}
\newcommand{\Prefr}{\Pre^\fr}
\newcommand{\PrefrA}{\Pre^\fr_{\A^1}}

\newcommand{\Corr}{\mathrm{Corr}}

\newcommand{\syncor}{\mathrm{syncor}}
\newcommand{\Nis}{\mathrm{Nis}}

\newcommand{\stau}{{\mathrm{s}\tau}}
\newcommand{\wtau}{{\mathrm{w}\tau}}
\newcommand{\mot}{\mathrm{mot}}
\newcommand{\smot}{\mathrm{smot}}
\newcommand{\tfsmot}{{\tf\text{-}\smot}}
\newcommand{\nissmot}{{\Nis\text{-}\smot}}

\newcommand{\K}{\mathrm{K}}

\newcommand{\CMon}{\mathrm{CMon}}

\newcommand{\Id}{\mathrm{Id}}

\newcommand{\Shv}{\mathrm{Shv}}
\newcommand{\sShv}{\mathrm{sShv}}
\newcommand{\bfShv}{\mathbf{Shv}}

\newcommand{\calO}{\mathcal O}

\newcommand{\defeq}{:=}

\newcommand{\cofib}{\operatorname{cofib}}
\newcommand{\fib}{\operatorname{fib}}

\newcommand{\tf}{\mathrm{tf}}
\newcommand{\nis}{\mathrm{nis}}
\newcommand{\A}{\mathbb{A}}
\newcommand{\Gm}{\mathbb{G}_m}

\newcommand{\overarrow}[1]{\overline{#1}}
\newcommand{\Lrep}{\mathcal{L}}
\newcommand{\gp}{\mathrm{gp}}

\newcommand{\Spt}{\mathrm{Spt}}
\newcommand{\PSpt}{\mathrm{PSpt}}
\newcommand{\Shift}{\mathrm{Shift}}

\newcommand{\id}{\mathrm{id}}

\newcommand{\calC}{\mathcal{C}}
\newcommand{\Map}{\mathrm{Map}}
\newcommand{\Cor}{\mathrm{Cor}}
\newcommand{\DM}{\mathbf{DM}}
\newcommand{\Cone}{\operatorname{Cone}}

\newcommand{\tr}{\mathrm{tr}}

\newcommand{\calF}{\mathcal{F}}

\newcommand{\fth}{{1\mathrm{th}}}

\newcommand{\ZF}{\mathbb Z\mathrm{F}}

\DeclareMathOperator{\Coker}{Coker}

\newcommand{\Func}{\mathrm{Func}}

\newcommand{\Ho}{\mathrm{Ho}}

\newcommand{\tgfr}{\mathrm{tgfr}}
\newcommand{\tghfr}{\mathrm{tghfr}}
\newcommand{\lownrfr}{\mathrm{lownhfr}}

\newcommand{\uppnrfr}{\mathrm{uppnfr}}

\newcommand{\Ab}{\mathrm{Ab}}

\DeclareMathOperator{\coker}{coker}

\title{
Stable motivic trivial fibre decomposition and framed motives over integers.
}
\author{Andrei Druzhinin}
\date{December 2021}
\address{Andrei Druzhinin, Chebyshev Laboratory, 
St. Petersburg State University, and
St. Petersburg Department of Steklov Mathematical Institute 
of Russian Academy of Sciences, Russia
}
\thanks{
The author was supported by a Young Russian Mathematics award.
}

\subjclass[MSC2020]{14F42,14G45}
\keywords{stable motivic homotopy group sheaves, base scheme case, localisation theorems, framed correspondences}

\begin{document}

\maketitle

\begin{abstract}
Using the trivial fiber topology we describe 
motivic $\infty$-loop spaces and fibrant replacements in the motivic stable homotopy category $\SH_{\A^1,\nis}(B)$ defined over one dimensional base schemes $B$. The part of the proof is the decomposition the stable motivic localisation into the one with respect to the $\tf$-topology and the Nisnevich localisation.
\end{abstract}

\tableofcontents

\section{Introduction}
The motivic categories listed below, namely, the Voevodsky motives and the Morel-Voevodsky motivic homotopy category in the stable and the liner version, are defined and studied over various base schemes in \cite{Voe-motives,Cisinski-Deglise-Triangmixedmotives,Morel-Voevodsky,morel-trieste,Jardine-spt,mot-functors} and numerous further works.
Whenever the functor \[Y\mapsto M(Y); \Sm_S\to \DM_{\A^1,\nis}(S),\text{ or }\mathbf{D}_{\A^1,\nis}(S)\text{ or }\SH_{\A^1,\nis}(S)\]
is constructed, there is the question on computations or presentations for (1) the motive $M(Y)$, or (2) the motivic homotopy groups of $Y$.
This means, in other words, a systematization and parametrization of motivic homotopies and loops of various types and a reconstruction of the motive and its homotopy groups via some procedures more elementary in comparing with the definition. 

For the case of base fields that satisfy the strict homotopy invariance theorem \cite{Voe-hty-inv,hty-inv} Vovodsky's theory \cite{Voe-hty-inv} and Garkusha{\&}Panin's framed motives theory \cite{Framed} give the formulas for $M(Y)$
\begin{equation}\label{eq:CorSigmaYFrSigmaY}
L_\nis\Cor(-\times\Delta^\bullet_k,\Sigma_{\Gm}^\infty Y)\in \DM_{\A^1,\nis}(k),\quad
L_\nis \Fr(-\times\Delta^\bullet_k,\Sigma^\infty_T Y)^\gp\in \SH_{\A^1,\nis}(k),\quad
\end{equation}
and imply that the stalks of motivic homologies and stable homotopies in $\DM_{\A^1,\nis}(k)$ and $\mathbf{D}_{\A^1,\nis}(k)$, $\SH_{\A^1,\nis}(k)$ on essentially smooth local scheme $U$ equal homologies and homotopies of the following complexes and the simplicial space respectively
\begin{equation}\label{eq:CorZFFrDeltaY}
\Cor(U\times\Delta^\bullet_k,Y),\; \ZF(U\times\Delta^\bullet_k,Y),\; \Fr(U\times\Delta^\bullet_k,Y)^\gp,
\end{equation}
where the bi-functors $\Cor$, and $\ZF$, $\Fr$ have precise geometrical description.
As shown in \cite{MVW,Nesh-FrKMW} the above formulas allow to prove the isomorphisms 
\[H^{n}_\mathrm{mot}(U,\mathbb Z(n))=\mathrm{K}^\mathrm{M}_n(U), \quad
\pi^{n,n}_\mathrm{smot}(\mathbb S)(U)=\mathrm{K}^\mathrm{MW}_n(U),\]
note that this is the original proof for the left one, and the right one is originally proven by Morel using different technique, see \cite{Morel-sphere-spt}.

The framed transfers formula \eqref{eq:CorSigmaYFrSigmaY} holds for any smooth scheme $Y\in \Sm_k$ providing reconstructions of the categories $\SH_{\A^1,\nis}(k)$, $\SH^\mathrm{eff}_{\A^1,\nis}(k)$, $\SH^\mathrm{veff}_{\A^1,\nis}(k)$ \cite{Framed,five-authors,BigFrmotives,FramedGamma}, and implies the $\PP^1$-stable connectivity theorem by Morel \cite{morel-trieste,Morel-connectivity}
that do not hold in general over base schemes because of Ayoub's contrexample \cite{AyoubcontrexempleA1connexite}.
This means that over base schemes, there exist some other type of motivic loops that is not reflected in the formulas \eqref{eq:CorSigmaYFrSigmaY}, \eqref{eq:CorZFFrDeltaY}, and make a non-trivial effort into the negative motivic homotopy groups of the relative smooth schemes $Y$ over the base scheme $S$ in general.

\subsection{The aim of the project and the present article.}\label{section:Plan-W-P} 
The project aims to generalise in a correct way the formulas \eqref{eq:CorSigmaYFrSigmaY} to the case of base schemes and get, in particular, the description of the lowest non-trivial negative homotopy group stalks in terms of generators and relations.
This means extracting the additional 'negative' obstructions skipped by the original formulas and fixing the formula using this extra type of loops.

The plan:\\ \vspace{1pt}
\begin{tabular}{|l|l|l|}
\hline 
& Formal Statement & Essential action \\
\hline
1)
& 
to generalise  
& 
to move or contract the part of
\\
&
the framed motives theorems \cite{Framed,hty-inv,framed-cancel,ConeTheGNP}
& 
motivic homotopy loops
\\
& 
in the appropriate way over base schemes
& 
over base schemes that behave
\\
&
& 
like in the case of a base field;
\\ \hline
2)
& 
to decompose 
&
to extract the rest
\\
& 
the stable motivic localisation into 
&
part of the obstructions
\\
&
the part similar to the base field case and
&
and organise it
\\
&
the part that concentrates the difference;
&
for further computation;
\\ \hline
3)
& 
to improve the formulas \eqref{eq:CorSigmaYFrSigmaY},
&
to parametrise generators and relations
\\
& 
to reconstruct 
&
in the negative
\\
& 
the categories $\SH_{\A^1,\nis}(S)$ and $\SH^\mathrm{veff}_{\A^1,\nis}(S)$.
&
stable motivic homotopy groups;
\\ \hline
\end{tabular}
\vspace{5pt}
\noindent

The project is started in \cite{DKO:SHISpecZ}, and we also refer to the talks and notes 
\cite{notesFramedmotivesoverSpecmathbbZ,StrictHomotopyInvarianceoverZandMotivicCategoriesoverBasQuestions}, where some principles and results are already explained and presented. In the present article, we document results that belong to points (1) and (2) and were not covered yet.
The basic ingredient that we use in the article to extract the additional type of loops is the trivial fibre subtopology of the Nisnevich topology defined in \cite[Definition 3.1]{DKO:SHISpecZ}; we use the notation $\mathcal L_{\tf}$ for the $\tf$-localisation endofunctor on the homotopy category of pointed simplicial presheaves $\Pre_\bullet(S)$, i.e. the localisation with respect to the trivial fibre topology,
and the notation $\mathcal L_{\A^1,\tf}$ for the unstable $(\A^1,\tf)$-motivic localisation endofunctor.

\subsection{Reduction to $\tf$-motivic localisation}
\label{subsect:stabletfmotivicdecomposition}
Recall that 
there are natural stable motivic weak equivalences of simplicial presheaves bi-spectra and $T$-spectra respectively
\[\Sigma^\infty_{\Gm,S^1}Y\simeq \Fr(-,\Sigma^\infty_{\Gm,S^1} Y),\quad \Sigma^\infty_{T}Y\simeq \Fr(-,\Sigma^\infty_{T} Y),\]
see \Cref{citeth:GP14MotivcEq} for details,
and moreover, any scheme-wise connective bi-spectrum or $T$-spectrum of simplicial presheaves, where $T=\A^1_S/(\A^1_S-0)$, 
has a stably motivically equivalent natural replacement by a scheme-wise connective
radditive bi-spectrum or $T$-spectrum of quasi-stable framed presheaves $\calF$ in sense of \cite{Framed},
see \Cref{rem:framedspectrarepalemcet}.
\begin{theorem}[\Cref{th:LsmotLnisLtfsmotFPretfr}]
Let $S$ be a noetherian separated scheme of Krull dimension $d$ 
that satisfies assumptions of \Cref{sect:AssumptionsBase}.
Let $\Lrep_{\nissmot}$ denote the the stable motivic localisation endofunctor on presheaves bi-spectra or $T$-spectra over $S$.

For any radditive bi-spectrum or $T$-spectrum of quasi-stable framed presheaves $\calF$, 
there is the canonical scheme-wise level-wise equivalence
\begin{equation}\label{eqintro:LnismotdecompostionLnisLtfsmot}\Lrep_{\nissmot}(\calF)\simeq \Lrep_\nis\Lrep_{\tfsmot}(\calF),\end{equation}
where $\Lrep_{\tfsmot}$ is the stable $\tf$-motivic localisation with respect to $\tf$-topology instead of the Nisnevich topology,
$\Lrep_{\nis}$ is equivalent to the level-wise Nisnevich localisation.
\end{theorem}
Note that 
\begin{itemize}
\item $\Lrep_{\tfsmot}$ takes scheme-wise level-wise connective objects to scheme-wise level-wise $(-d)$-connective ones,
\item $\Lrep_{\nis}$ preserves Nisnevich locally connective objects.
\end{itemize}
Consequently, the above theorem allows to recover the connectivity result of \cite{ConnDodekindDomains,ConnGabPresLemNoethDominffield,ConnBase} for motivic $\PP^1$-spectra over the base schemes that satisfy assumptions from \Cref{sect:AssumptionsBase}.
\begin{corollary}
\label{th:connectivitysmotsheaves}
Let $S$ be as in \Cref{sect:AssumptionsBase}.
For any $Y\in \Sm_S$ the assocated Nisnevich sheaves of stable motivic homotopy groups $\pi_{i,j}(Y)$ are trivial for $i<j-\dim S$.
\end{corollary}
Moreover, 
$\Lrep_{\tfsmot}$ concentrates the negative direction obstructions and has essentially simpler structure properties in comparing with $\Lrep_{\nissmot}$ because for any scheme-wise connective $\calF$ there is the vanishing of presheaves \[[-,\calF\wedge S^l]_{\SH^\fr_{\A^1,\tf}(S)}=\pi_{-l}(\Lrep_{\tfsmot}(Y))=0\] for all $l>d$, while for the Nisnevich case this holds only Nisnevich locally. 
Note that the latter property of $\Lrep_{\tfsmot}$ is proven in \cite[\S 14]{DKO:SHISpecZ} already for any base scheme,
and holds in short because the $\tf$-cohomological dimension of an arbitrary scheme $U\in \Sch_S$ is bounded by $\dim S$.
The proof of the equivalence \eqref{eqintro:LnismotdecompostionLnisLtfsmot} uses the assumption from \Cref{sect:AssumptionsBase} and
requires mentioned in point (1) of the plan appropriate lifting to the base scheme case for same set of essential results of Grakusha{\&}Panin's theory.

\subsection{$\infty$-categorical form}
\label{subsect:categoricalform}
To compute the unit of the adjunction \[\Sigma^\infty_{\PP^1}\colon \Pre(\Sm_S) \rightleftarrows  \SH_{\A^1,\nis}(S)\colon \Omega^\infty_{\PP^1}\] we decompose it into the sequence
\begin{equation}\label{eq:PreSHstAtffrSHstAnisfrSHstSHP1}\Pre(\Sm_S) \rightleftarrows \SH^{s,t}_{\A^1,\tf}(\Corr^\tfr_S) \stackrel{L_\nis}{\rightleftarrows} \SH^{s,t}_{\A^1,\nis}(\Corr^\tfr_S)\simeq \SH^{s,t}_{\A^1,\nis}(S)\simeq \SH_{\A^1,\nis}(S),\end{equation} where the latter adjunction is the equivalence of stable motivic homotopy categories of $\PP^1$-spectra and (s,t)-bi-spectra, 
the second right one is proven by \cite[Thorem 18]{Hoyois-framed-loc}, where $\Corr^\tfr(S)$ is the $\infty$-category of framed correspondences \cite{five-authors},
and by the result of the present article the unit of the adjunction denoted $L_\nis$ agrees with the usual Nisnevich local replacement on $\Pre(\Sm_S)$. 

\begin{theorem}[\Cref{th:SHnisA1S1GmSHfrtfA1S1Gm}]
Let $S$ be a base scheme satisfying assumptions from \Cref{sect:AssumptionsBase}, for example, $\Spec \mathbb Z$.

There are reflective embeddings
\[\SH_{\A^1,\nis}(S)\to \SH_{\A^1,\tf}(\Corr^\tfr_S), \quad\SH^\mathrm{veff}_{\A^1,\nis}(S) \to \SH^\mathrm{veff}_{\A^1,\tf}(\Corr^\tfr_S)\]
with images spanned by Nisnevich sheaves in the right-side categories; the left adjoint functors preserve the motives of smooth schemes,
and for a local henselian essentially smooth $U$, and any $Y\in \Sm_S$ \[\pi_{*,*}^{\SH_{\A^1,\nis}(S)}(Y)(U)\simeq\pi_{*,*}^{\SH^\fr_{\A^1,\tf}(S)}(Y)(U).\]
where $\SH^{\fr}_{\A^1,\tf}(S)=\SH_{\A^1,\tf}(\Corr^\tfr_S)$
\end{theorem}

\subsection{The result on $\infty$-loop spaces.}\label{sect:ThoremOmegaLnisLGmLA1tfFr}
The ${\PP^1}$-motivic $\infty$-loop space of $Y\in \Sm_S$ is the image along with the composition 
\begin{equation}\label{eq:OmegaSigmaPw1}\Sm_S\xrightarrow{\Sigma_{\PP^{\wedge 1}}^\infty} \SH_{\A^1,\nis}(S)\xrightarrow{\Omega^\infty_{\PP^{\wedge 1}}}\Shv_\bullet(S),\end{equation}
where $\mathrm{Shv}_\bullet(S)$ is the homotopy category of pointed Nisnevich simplicial sheaves. 
The principle of computations of the motivic $\infty$-loop space that we discuss is based on the presentation of the latter space in terms of simpler and more elementary $\infty$-loop spaces.
The zero term of spectra \eqref{eq:CorSigmaYFrSigmaY} computes the stalks of ${\PP^1}$-motivic $\infty$-loop space of $Y$ as
\[\Omega^\infty_{\PP^{\wedge 1}}\Sigma^\infty_{\PP^{\wedge 1}}(Y)(U)\simeq\Fr(U\times\Delta^\bullet,Y)^\gp.\]
This result can be naturally explained by the fact that the framed correspondences sheaf according to Voevodsky's lemma, see \cite{Voe-notes,Framed} and \cite{five-authors} for the case of base schemes, is given itself by a certain $\infty$-loop space construction namely,
\begin{equation}\label{eq:VoevodskyLemma}\Fr(-,Y\wedge T^r)= \Omega_{\PP^{\wedge 1}}^l \Sigma^{l+r}_{T}(Y)=\varinjlim_l \mathrm{Shv}_\bullet(-\wedge \PP^{\wedge l}, Y\wedge T^{\wedge l+r}).\end{equation}

As was mentioned above, over a base scheme $S$, there should exist some other type of motivic homotopy loops not parametrised by the framed correspondences.
In the following theorem, we reconstruct the motivic $\infty$-loop space \eqref{eq:OmegaSigmaPw1} of $Y\in \Sm_S$ via the cascade combination of framed correspondences presheaf \eqref{eq:VoevodskyLemma} with one of the following $\infty$-loop space endofunctors on the category of presheaves with framed transfers
\begin{equation}\label{eq:zerothtermLGmA1tf}  
\varinjlim_l \Omega_{\Gm}^l \Lrep_{\A^1,\tf}\Sigma^l_{\Gm},\;\;
\varinjlim_l \Omega^l_{\Gm^{\wedge 1}\wedge S^1_S}\Lrep_{\tf}\Sigma^{l}_{\A^1_S/(\A^1_S-0)},\end{equation}
where $S^1_S=\Delta^1_S/\delta\Delta^1_S$.
All the formulas \eqref{eq:VoevodskyLemma}, \eqref{eq:zerothtermLGmA1tf} are kinds of motivic $\infty$-loop spaces in the sense that they are compositions of the infinite suspension and the loop functors for some models of motivic sphere
and a certain motivic equivalence in between, 
that is hidden in the changing of the motivic sphere $\PP^{\wedge 1}\to T$ for \eqref{eq:VoevodskyLemma} and $\Gm^{\wedge 1}\wedge S^1_S\simeq \A^1_S/(\A^1_S-0)$ for \eqref{eq:zerothtermLGmA1tf}, and the functors $\mathcal L_{\A^1,\tf}$, or $\mathcal L_{\tf}$.

\begin{theorem}\label{th:Lacalisationfunctordecomposition}
Let $S$ be as in \Cref{sect:AssumptionsBase}, and $Y\in \Sm_S$.
The $\infty$-loop space in $\Shv_\bullet(\Sm_S)$ of the motivic suspension spectrum of $Y$ in $\SH_{\A^1,\nis}(S)$
is Nisnevich locally equivalent to 
\begin{equation*}\label{eq:OmegastTGSAOmegastTTinftyLmot(Y)}\begin{array}{lllll}
\varinjlim_{l}&\Omega^l_{\Gm^{\wedge 1}}&(\Lrep_{\A^1,\tf}&\Fr(Y\wedge \Gm^{\wedge l}))^\gp&\simeq\\
\varinjlim_{l}&\Omega^l_{\Gm^{\wedge 1}\wedge S^1}&\Lrep_{\A^1,\tf}&\Fr(Y\wedge T^{\wedge l})&\simeq\\ 
\varinjlim_{l}&\Omega^l_{\Gm^{\wedge 1}\wedge \Delta^1_S/\delta\Delta^1_S}&\Lrep_{\tf}&\Fr(Y\wedge T^{\wedge l})&\in \Shv_\bullet(\Sm_S),
\end{array}\end{equation*}
where $\Fr(Y)=\Fr(-,Y)$ is the presheaf of framed correspondences in sense of \cite{Framed}.
\end{theorem}
\begin{remark}
The above theorem implies that the functor $\Omega^\infty_{\Gm}\colon \DM(S)\to \mathbf{D}(\Shv^\Ab_\nis(\mathrm{Cor}_S))$ takes the motive of $Y$
to the complex Nisnevich locally quasi-isomorphic to 
\begin{equation*}
\varinjlim_{l}\Omega^l_{\Gm^{\wedge 1}}\Lrep_{\A^1,\tf}(\mathbb Z_\tr(Y\wedge \Gm^{\wedge l}))\in \mathbf{D}(\Shv^\Ab_\nis(\mathrm{Cor}_S)).
\end{equation*}
\end{remark}

\subsection{$\tf$-motivic localisation, description of stalks, and
further results}
For a base field $k$ the $\tf$-topology is trivial on $\Sm_k$, so $\Shv_\tf(k)\simeq \Pre(k)$, $\Lrep_{\A^1,\tf}F=F(\Delta^\bullet\times-)$. 
The stable and unstable $\A^1$-$\tf$-motivic localisation is the subject of the further study according to the plan from \Cref{section:Plan-W-P} and its behaviour is simpler than for the case of the Nisnevich topology because 
the functors $-\times\A^1$, and $-\times\Gm$ do not change the $\tf$-cohomological dimension of smooth schemes of positive relative dimension over $S$, for any base scheme $S$.
We have in a joint work with Kolderup and {\O}stv{\ae}r in process
as was mentioned in \cite{DKO:SHISpecZ}, 
and \cite{notesFramedmotivesoverSpecmathbbZ} outlined some preparation for the proof,
the computation by the formula 
\begin{equation}\label{eq:LA1tf}\Lrep_{\A^1,\tf} = 
(\Lrep_{\A^1}\Lrep_{\tf})^d\Lrep_{\A^1}, d=\dim S.
\end{equation}
Such finite-size formula in combination with \Cref{th:Lacalisationfunctordecomposition} allows us according to the above plan
to get a reconstruction of the category $\SH^\mathrm{veff}_{\A^1,\nis}(S)$ and get generators and relations in the lowest motivic homotopy groups in terms of some algebrogeometric data
that description and analysis are left for future work.

In the present work, nevertheless, already as a consequence of \Cref{th:Lacalisationfunctordecomposition} we can get some weaker description.
\begin{corollary}
Let $S$ be a local scheme that satisfy assumptions of \Cref{sect:AssumptionsBase}, and let $\dim S=1$.
Let $z$ denote the closed point
Then for any $Y\in \Sm_S$ there is the isomorphism of Nisnevich stalks
\[\pi_{-1,0}^{\SH_{\A^1,\nis}(S)}(Y)\simeq \varinjlim_{l} H_\tf^1( (\Gm^{\wedge 1}\wedge S_{S}^1)^{\wedge l}\times -, \ZF(Y\wedge T^{\wedge l})_\tf ), 
\]
and note that for any presheaf $F$ there is the isomorphism
\[ H_\tf^1(U,F_\tf) \cong \coker \left( F( U\times_S(S-z) )\oplus F( U^h_z ) \to F( U^h_{z}\times_S(S-z) ) \right) .
\]
So the formula above is the description of the stalks of the sheaf $\pi^{\SH_{\A^1,\nis}(S)}_{-1,0}(Y)$ in terms of generators and relations parametrised by countable set of algebrogeometric parameters that satisfy countable set of relations and conditions.
\end{corollary}

As mentioned above, 
the additional type of motivic homotopy loops that we discuss throughout the introduction are essential for the structural properties of categories $\SH_{\A^1,\nis}(S)$ and $\DM_{\A^1,\nis}(S)$, they affect the shift of t-structure along with the stable motivic localisation and appear in the computations of negative motivic homotopy groups of smooth schemes in general.
Nevertheless, these loops do not appear in the answers for the $\infty$-loop spaces for a wide class of often studied motivic spectra and cohomology theories in particular for the sphere spectrum.
\begin{conjecture}[Particular case of \protect{\cite[conjecture ]{tnpsCompactificationsacCousin}}]\label{conj:S_S} 
\item[] For a regular base scheme $S$ of a finite Krull dimension
\[\Lrep_\smot(
\mathbb S_S^{j,j})\simeq \Lrep_\nis\Fr(\Delta^\bullet\times -, \Gm^{\wedge j})^\gp\]
\item[] the stable motivic homotopy group sheaves $\underline{\pi}_{i,j}(\mathbb S_S)$ in $\SH_{\A^1,\nis}(S)$ and $\underline{\pi}_{i}(\mathbb S_S)$ in $\SH^{S^1}(S)$, for $i<j$ are trivial, and 
$\pi_{-j,-j}(\mathbb S_S)\simeq \ZF(-, \Gm^{\wedge j})/\ZF(\A^1\times -, \Gm^{\wedge j})$.
\end{conjecture}



\subsection{Conditions on the base scheme.}
\label{sect:AssumptionsBase}

In the article, we assume and use the following two properties of the base scheme $S$
proven in the literature 
at the moment for schemes $S$ of Krull dimension 1. 
\begin{itemize}
\item[(1)] For any $\sigma\in S$ and $\A^1$-invaraint $F\in \Spt_s^\fr(S)$ the object $L_\nis F$ is $\A^1$-invariant.
\item[(2)] For any $\eta\in S$, $X\in \Sm_S$, $z\in X$, and quasi-stable $\A^1$-invariant $F\in \Spt^\fr_s(S)$ there canonical morphism $F( (X^h_z)_\eta )\to L_\nis F( (X^h_z)_\eta ) $ is a stable equivalence.
\end{itemize}
The strict homotopy invariance theorem for $\A^1$-invariant additive abelian presheaves with framed transfers is already proven for all base fields in \cite{SHIfrfield}, and consequently the property of point (1) holds for all schemes $S$.
\begin{remark}
The work \cite{hty-inv} covers the case of 
the above strict homotopy invariance theorem 
for perfect infinite fields $k$, $\operatorname{char}k\neq 2$;
the work \cite{surj-etale-exc} covers the case of $k$ with $\operatorname{char}k=2$ and finite fields;
the works \cite{DrKyllfinFrpi00,five-authors} both cover the case of finite fields.
\end{remark}
\noindent The point (2) is proven in \cite[\S 8]{DKO:SHISpecZ} for $S$ of Krull dimension 1 under the assumption of the strict homotopy invariance theorem over the generic point, so it holds for schemes $S$ of Krull dimension 1.

\subsection{Novelty of the present article.}

Summarising, the results of the present article reduce 
the stable motivic localisation in $\SH(S)$ to the unstable $\tf$-motivic localisation $\Lrep_{\A^1,\tf}$,
and 
additionally reduce the motivic infinite loop space of the form $\Omega_{S_{\A}^1\wedge \Gm}\Sigma^\infty_T$ to $\Lrep_\tf$.
The present article continues the series started by \cite{DKO:SHISpecZ}.
The latter article proves the results for the case of $S^1$-spectra, and contains some preparation for (s,t)-bi-spectra.
The present one completes the study for bi-spectra, covers the case of $\PP^1$-spectra, and extends the generality to the case of quotient spaces $X/U$ for an open subscheme $U$ in a smooth scheme $X$. 
The novelties can be summarised as follows: \begin{itemize}
\item
the scheme of the reasoning to achieve the results for bi-spectra in \Cref{section:LocalisationTheoremtf,sect:stablemotiviclocalisation,section:SplittingdiagramLocalisation} differs from the one used for $S^1$-spectra and suggested for bi-spectra in \cite{DKO:SHISpecZ},
the original way would be longer and technically complicated to achieve the results of the present article;
\item
the cone theorem, see \Cref{section:ConeThoreom}, is
proven for any base scheme and an arbitrary smooth open pair
to prove the results on $T$-specrta and $\PP^1$-spectra and on the suspension spectrum of the motivic space $X/U$;
\item
the counterexample for the cancellation theorem is obtained, see \Cref{section:ContrexampleCancellationDMeff}.
\end{itemize}
We should note:

1) 
The main result of \cite{DKO:SHISpecZ} is the strict homotopy invariance theorem \cite[Theorem 11.2]{DKO:SHISpecZ} implies the result like in \Cref{subsect:stabletfmotivicdecomposition,subsect:categoricalform} for the case of $S^1$-spectra, see \cite[\S 14.1]{DKO:SHISpecZ}.
\cite[Theorem 11.2]{DKO:SHISpecZ} does not imply the result for (s,t)-bi-spectra, 
\cite[\S 14.3]{DKO:SHISpecZ} demonstrates an ability to apply the method and the technique of the work to prove the equivalence $\Omega_{\Gm}\Lrep_\mot\simeq \Lrep_\mot\Omega_{\Gm}$.
Such a way would allow to achieve the first equivalence in \eqref{eq:OmegastOmegaPinftyLmot(Y)};
at the same, this would suggest to repeat all the arguments written in \cite[\S\S 7,9,10,11]{DKO:SHISpecZ} for the case of $\Omega_{S^1,\Gm}$-bi-spectra instead of $\Omega_{S^1}$-spectra.
The involving of $\Gm$-loops throughout the arguments increases the length, 
and would require to add more details and change some inner argument formulations inside the proofs that is not short already for $S^1$-spectra.
We cover the case of bi-spectra directly via alternative and formally independent argument.

The most essential changing that we use to simplify some part of the argument we combine $\Gm$-loops functor and the functor $F\mapsto F(-\times\Delta^\bullet)$ into the one stablisation procedure, and using the reformulations provided by \Cref{sect:reformlemma} prove the result in different equivalent forms for different functors of the diagram \eqref{eq:Smaf(Z)A1Smaf(SZ)nisSmaf(S)(U)}.
The original argument followed the standard viewpoint that suggests to consider the $\A^1$-localisation and the stablilisation
after the Nisnevich localisation.

Additionally, to prove the required properties of the adjunction 
$\Pre^\fr(\Sm_Z)\leftrightarrows \Pre^\fr(\Sm_{S,Z})$ in \Cref{sect:DeformationZSZ}
in distinct to \cite[\S 10]{DKO:SHISpecZ} 
we formally use the fibrancy of the $\infty$-categorical or simplicially enriched categories functor $\Pre^\fr_{\A^1}(\Smat_{S,Z})\to \Pre^\fr_{\A^1}(\Smat_{Z})$ formulated in terms of the $\infty$-category $\Corr^\fr(S)$ from \cite{five-authors}. 
The argument in terms of the classical category $\Fr_+(S)$ \cite[\S 10]{DKO:SHISpecZ} concerned the rigidity property of presheaves on categories $\Smat_{S,Z}$ avoided in the present text.
Nevertheless, we should say the arguments in \cite[\S 6]{DKO:SHISpecZ} that prove the equivalence $\HH_{\mot}(\Fr_+(\Sm_{B,Z}))\simeq \HH_{\mot}(\Fr_+(\Sm_{Z}))$ and provided the basement for \cite[\S 10]{DKO:SHISpecZ} are morally equivalent to  the mentioned $\infty$-categorical fibrancy.

The article \cite{DKO:SHISpecZ} and the present one are written in the formally different but equivalent settings of model structures and $\infty$-categories with the use of the Voevodsky framed correspondences that study is developed by Garkusha and Panin in \cite{Framed} and the tangentially framed correspondences  \cite{five-authors} by Elmanto, Hoyois, Khan, Sosnilo, Yakerson respectively.

2) 
The cone theorem in the form of the $S^1$-stable motivic equivalence
$\Fr(-,\Sigma^\infty_{S^1} X/U)\simeq \Fr(-,\Sigma^\infty_{S^1} X)/\Fr(-,\Sigma^\infty_{S^1} U)$
is proven directly for an arbitrary smooth $Y$ and open subscheme $U$ over
noetherian separated base scheme $S$ independently from earlier known cases and basing on the different strategy.
Instead of the constructions of $\A^1$-homotopies like it was in \cite{ConeTheGNP} and \cite{SmAffOpPairs} we prove some lifting properties and the technique of \cite{five-authors}. 
Note that we consider that the geometrical part of the original proof for the case of $\A^n_k/(\A^n_k-0)$ in \cite{ConeTheGNP} could be applied at least to finite dimensional base schemes too.

Nevertheless, there is an alternative strategy to deduce the result for any $Y/U$ for $Y\in \Sm_S$ at least for finite dimensional base scheme $S$ with the use of the localisation theorem and the results over base fields. 
The latter one is proven in its turn for affine schemes $Y$ over infinite base fields and then deduced for finite fields, see  \cite{SmAffOpPairs} and the first arXiv version for finite fields argument generalised later and moved to \cite{DrKyllfinFrpi00}, and see \cite{SmModelSpectrumTP}, or \cite{framed-notes}, or \cite{Segalthopenpairsofsmoothschemes} for non-affine schemes case.
The situation is close to the one with the equivalences $\SH(S)\simeq\SH^\fr(S)$ proven in \cite{Hoyois-framed-loc} with the use of the localisation theorems for both $\HH^\fr(S)$ and $\HH(S)$, and the results of the theory of framed motives over perfect fields that are proven in \cite{Framed,five-authors} in its turn firstly over infinite fields and consequently deduced for the finite fields case.
At the same time, the equivalence $\SH(S)\simeq\SH^\fr(S)$ is a consequence of \cite[Theorem 4.1]{Framed} that argument. as we consider, holds over any base scheme.

\section{Notation}

We denote by $\Pre(S)$ and $\Pre^\fr(S)$ the categories of presheaves on $\Sm_S$ and  $\Corr^\fr(S)$.
For any category as above we consider the categories of prespeactra and spectra in sense of in sense of \Cref{den:PreGmT(SmS)}. 
We use notation
\par
\begin{tabular}{l|l}
$\Pre^{\mathcal S}(S) = \PSpt^{\mathcal S}$ & 
prespectra in sense of \Cref{den:PreGmT(SmS)} with respect to $\mathcal S$, \\ &where $\mathcal S = S^1,\Gm,T$, and $T= \A^1_S/(\A^1_S-0)$;\\
$\Spt^{\mathcal S}(S) = \Pre^{\mathcal S}_{\mathcal S}(S)$ & the categories of spectra with respect to $\mathcal S$;\\
$\mathcal P_{\A^1}(S)$, $\mathcal P_{\tau}$ & the subcategories of $\A^1$-invariant and $\tau$-local objects;\\
$\mathcal P_{\A^1,\tau}=\mathcal P_{\tau}\cap \mathcal P_{\A^1}(S)$& the subcategory of motivic objects\\
$\Shv_{\tau}(S) = \Pre_{\tau}(S)$ &  the $\infty$-categories of sheaves \\
$\bfShv_{\tau}(S) = \mathrm{Ho}\,\Shv_{\tau}(S)$ &  the homotopy categories of sheaves \\
$\rmHH_{\A^1,\tau}(S) = \Pre_{\A^1,\tau}(S)$ & the $(\A^1,\tau)$ motivic $\infty$-category  \\
$\rmSH^{\mathcal S}_{\A^1,\tau}(S) = \Spt^{\mathcal S}_{\A^1,\tau}(S)$ & the $\mathcal S$ stable $(\A^1,\tau)$ motivic $\infty$-category \\
$\HH_{\A^1,\tau}(S) = \mathrm{Ho}\,\rmHH_{\A^1,\nis}(S)$ & the $(\A^1,\tau)$ motivic homotopy category \\
$\SH_{\A^1,\tau}(S) = \mathrm{Ho}\,\rmSH^{\mathcal S}_{\A^1,\nis}(S)$ & the $\mathcal S$ stable $(\A^1,\tau)$ motivic homotopy category \\
\end{tabular}
\par We use similar notation for categories related to presheaves on $\Corr^\fr(S)$. We use the symbol $\gp$ to indicate the group completions.

We use notation for the localisation functors 
\[L^{\mathcal P}_{\alpha}\colon \mathcal P\to \mathcal P_{\alpha},\]
and the notation 
\[\Lrep^{\mathcal P}_{\alpha}\colon \mathcal P\to \mathcal P_{\alpha}\to \mathcal P\]
for the corresponding endo-functor on $\mathcal P$, where 
$\alpha=\A^1$, or $\alpha=\tau$, or $\{\A^1,\tau\}$.
We use notation 
\[\begin{array}{lcl} L_{\mathcal S}&\colon& \Pre^{\mathcal S}(S)\to \Spt^{\mathcal S}(S)\\
L_{\mathcal S}&\colon&\Pre^{\mathcal S}(S) \to \Pre^{\mathcal S}(S)\end{array}\]
for the localisation functor and endofunctor, see \Cref{section:StableLocalisationDefinitionsNotations}, and use the notation $\simeq_{\mathcal S}$ for the corresponding type of equivalences.

\section{Framed correspondences and category $\SH^\fr(S)$.}

In the unpublished notes \cite{framed-notes}, Voevodsky suggested a computational approach for the stable motivic homotopy category $\SH(k)$ 
based on so called framed correspondences, which would replace $\mathrm{Cor}$-correspondences used for Voevodsky's motives category $\DM(k)$. 

\begin{definition}[Voevodsky's framed correspondences.]\label{def:FrCor}
Let $X,Y\in \Sm_{S}$ for a base scheme $S$, $W$ be an open subscheme in $Y$.
A level $n$ explicit framed correspondence from $X$ to $Y/W$ 
is given by the data $(Z,V,\phi,g)$:\begin{itemize}
\item[(i)] a closed subscheme $Z\subset \A^n\times X$ finite over $X$;
\item[(ii)] an \'etale neighbourhood $V$ of $Z$ in $\A^n\times X$, 
\item[(iii)] regular functions $\varphi_i\in \calO( V )$, $i=1,\dots n$; 
\item[(iv)] 
and a regular map $g\colon V\to Y$
such that
$V \times_{\varphi_1,\dots \varphi_n,g} (0\times (Y\setminus W)) = Z$.
\end{itemize}

Denote by 
$\Fr_n( X , Y )=\Fr_n( Y )( X )$ the pointed set of equivalence classes of framed correspondences of level $n$ from 
$ X $ to $ Y $ with respect to the shrinking of the \'etale neighbourhood $V$ of $Z$, and pointed at the class of framed correspondences with empty $Z$.
Denote $\Fr(X,Y)= \Fr( Y )( X )=\varinjlim_{n}\Fr_n( Y )( X )$, see \cite[Definition 2.8]{Framed} for the transition maps.
\end{definition}
We refer to \cite[Definitions 2.3, 2.7]{Framed} for the composition of framed correspondences and the graded category of framed correspondences $\Fr_*(S)$ with mapping sets $\bigvee_n\Fr_n(X,Y)$. 
Any framed correspondence over a base $S$ defines a $\Cor$-correspondence, see \cite[\S 13]{DKO:SHISpecZ}.
Moreover, the following Voevodsky's lemma implies, in particular, that any $\SH(S)$-representable cohomology theory has framed transfers, i.e. defines a presheaf on $\Fr_*(S)$.
\begin{lemma}[Voevodsky's lemma]
There is a natural isomorphism
\[\Fr_n(X,Y)\simeq \Shv^{\mathrm{Set}}_{\nis,\bullet}(X_+\wedge \PP^{\wedge n},Y_+\wedge T^{\wedge n}),\]
where $\PP^{\wedge n}=(\PP^1/\infty)^{\wedge n}$, $\Shv^{\mathrm{Set}}_{\nis,\bullet}$ denotes the category of Nisnevich sheaves of pointed sets.
\end{lemma}
Voevodsky suspected that the use of the framed correspondences over perfect base fields would 
provide formulas for the motive of smooth schemes in $\SH(k)$ similar to the formula for $\DM(k)$
\begin{equation}\label{eq:C(Ztr(X))DM}M(X)= \Lrep_\nis C^*(\mathbb Z_\tr(X))\in \DM_{\mathrm{eff}}(k)\hookrightarrow D(\mathrm{Pre}^\tr(k)),\end{equation}
where $\Lrep_\nis\colon D(\Pre^\tr(k))\to D(\Shv^\tr(k))\to D(\Pre^\tr(k))$,
and $\Pre^\tr(k)=\Pre^\Ab(\mathrm{Cor}(k))$.
This was realised by Garkusha-Panin's theory of framed motives.
Consider the suspension $T$-spectrum 
\begin{equation}\label{eq:SigmainftyT}
\Sigma^\infty_T X_+=(X_+, X_+\wedge T, X_+\wedge T^{\wedge 2}, \dots ),
\end{equation} and the $T$-spectrum
\begin{equation*}
\Fr(\Sigma^\infty_T X_+)=(\Fr(X_+), \Fr(X_+\wedge T), \Fr(X_+\wedge T^{\wedge 2}),\dots, )
\end{equation*}
with structure morphisms induced by the ones of \eqref{eq:SigmainftyT}.
\begin{theorem}[\protect{Garkusha-Panin \cite{Framed}}]\label{citeth:GP14MotivcEq}
Let $S$ be a noetherian separated scheme of finite Krull dimension. 
The canonical morphism 
\begin{gather}\label{eq:SigmaCorSGXLnisLAFrSigmaTX}
\begin{array}{lcl}
\Sigma_{\PP^1}^\infty X
&\to& 
\Lrep_\nis\Lrep_{\A^1}\Fr(\Sigma^\infty_T X_+),\\
\Sigma_{S^1,\Gm}^\infty X
&\to& 
\Lrep_\nis\Lrep_{\A^1}\Fr(\Sigma^\infty_{S^1,\Gm} X_+).
\end{array}
\end{gather}
is a stable motivic equivalence of $\PP^1$-spectra and bi-spectra respectively,
and each space $\Lrep_{\A^1}\Fr(X_+\wedge T^l )$ has natural structure of special Segal's $\Gamma$-space.
\end{theorem}
\begin{proof}[First proof]
For the case of an infinite prefect base field 
the first claim is formulated in \cite[Theorem 4.1(1)]{Framed},
the second one is given by \cite[\S 11, Theorems 11.4]{Framed},
and the third by \cite[\S 6-7, Theorem 6.5]{Framed}.
The arguments for the first equivalence hold word by word for the base scheme case 
with the use of the Voeovdsky lemma proven in \cite[Appendix A.1]{five-authors}.
The arguments for the second equivalences follows from the first word by word as in \cite[\S 11]{Framed} with the use of
\Cref{th:Fr(Y/V)}.
See also \cite{SmAffOpPairs} for details on the base scheme case.
\end{proof}
\begin{proof}[Second proof]
The claim for bi-spectra follows by \cite[Thorem 18]{Hoyois-framed-loc} in combination with \cite[Corollary 2.3.25]{five-authors}, and the case of $T$-spectra follows for the bi-spectra case by \Cref{th:Fr(Y/V)}.
\end{proof}
\begin{remark}\label{rem:framedspectrarepalemcet}
The mappings \Cref{eq:SigmaCorSGXLnisLAFrSigmaTX} extend to the category of simplicial schemes, and allow to make a 
stably motivically equivalent replacement of an arbitrary scheme-wise connective bi-spectrum or $T$-spectrum by a 
radditive bi-spectrum or $T$-spectrum of quasi-stable framed presheaves $\calF$ in sense of \cite{Framed}
as explained in \cite[\S 11]{Framed}.
\end{remark}
\begin{theorem}[\protect{Garkusha-Panin \cite[Theorem 4.1, Theorem 11.7]{Framed}}]\label{citeth:GP14MotivcReplacement}
Let $k$ be a perfect field.
The right side of 
\eqref{eq:SigmaCorSGXLnisLAFrSigmaTX}
is 
a $\A^1$-invariant Nisnevich local $\Omega_{\PP^1}$-spectrum in positive degrees,
and
a $\A^1$-invariant Nisnevich local $\Omega_{s,t}$-bi-spectrum in positive degrees with respect to $S^1$
respectively.
\end{theorem}
\begin{remark}
In the above formulation Garkusha-Panin's results from \cite{Framed} are combined with the results of \cite{DrKyllfinFrpi00} and \cite{surj-etale-exc} that cover the cases of finite fields, and fields of characteristic 2.
\end{remark}
\begin{remark}
The results in \cite{Framed} are formulated in terms of the Morel-Voevodsky model structure on the category of simplicial Nisnevich sheaves.
\end{remark}
So the mapping
\begin{equation}\label{eq:LnisLA1Fr(SigmaTX)gp}\Sigma_{\PP^1}^\infty X\to 
\Lrep_\nis(\Lrep_{\A^1}\Fr(\Sigma^\infty_T X_+))^\gp,\end{equation}
where $(\Lrep_{\A^1}\Fr(\Sigma^\infty_T X_+))^\gp$ denotes the group-completion of special Segal's $\Gamma$-space, 
take $\Sigma_{\PP^1}^\infty X$ 
to an $\A^1$-invariant Nisnevich local 
motivically equivalent replacement that, and consequently
\[\Omega^\infty_{\PP^1}\Sigma^\infty_{\PP^1} X_+\simeq\Lrep_\nis(\Lrep_{\A^1}\Fr(X_+))^\gp.\]
Later the works \cite{BigFrmotives,FramedGamma} by Garkusha and Panin, and Garkusha and Panin and {\O}stv{\ae}r provided a reconstructions of $\SH^\mathrm{veff}(k)$, and $\SH(k)$ via framed motives approach over the same class of base fields like as in \Cref{citeth:GP14MotivcReplacement}.

In parallel to \cite{Framed,BigFrmotives,FramedGamma} the work \cite{five-authors} by Elmanto, Hoyois, Khan, Sosnilo, Yakerson 
provided a construction of an alternative variant of framed correspondences, and an alternative picture of the theory with the use of the main results of the works of Garkusha-Panin's project  \cite{hty-inv,framed-cancel,ConeTheGNP}.
\begin{definition}[tangentially framed correspondences, \protect{EHKSY \cite{five-authors}}]
A \emph{tangentially framed correspondence} 
from $X$ to $Y$, where $X,Y\in \Sm_{S}$, 
is defined by the data:
\begin{itemize}
\item[(i)] 
A finite syntomic $f\colon Z\to X$,
and a morphism $g\colon Z\to Y$;
\item[(ii)] 
a trivialisation $0\simeq \mathcal L_{f}$ of the cotangent complex $\mathcal L_{f}$ in $K$-theory $\infty$-gruppoid on $Z$.
\end{itemize}
For a given $Y$ 
define a presheaf of topological spaces $h^\tfr(Y)$ that value on $X$ is given by the coproduct of the path-spaces in the $K$-theory $\infty$-gruppoid from $0$ to $\mathcal L_{f}$, i.e.
\[h^\tfr(Y)(X) = \coprod_{(f,Z,g)} 0\times_{\mathrm{K}(Z)}\mathcal L_{f}\]
with the coproduct along the triples $(f,Z,g)$ described in (i) above.
\end{definition}
The article \cite{five-authors} provides a construction of
the $\infty$-category $\Corr^\tfr(S)$ with objects being smooth schemes over $S$  
such that the presheaves $h^\tfr(X)$ defined above 
are the mapping spaces presheaves $\mathrm{Map}_{\Corr^\tfr(S)}(-,X)$ in the category $\Corr^\tfr(S)$.
This leads to an $\infty$-categorical revision of the theory of framed motives, recognition result for $\SH^\mathrm{veff}(k)$,
and makes the framed correspondences and $\Cor$ correspondences approaches being formally parallel.

Denote 
$\Pre^\tfr(S)=\Pre(\Corr^\tfr(S))$, and denote by
$\Pre^\tfr_\nis(S)$, $\Pre^\tfr_{\A^1}(S)$, $\Pre_\nis^\tfr(S)_\mathrm{gp}$
the subcategories of Nisnevich sheaves, $\A^1$-invariant, and group-like objects in $\Pre^\tfr(S)$.
Define $\mathbf H^\fr(S)$ as the homotopy category of the intersection $\Pre_{\nis,\A^1}^\tfr(S)_\mathrm{gp}$ of the above categories
\[
\mathbf H^\fr(S)=\Ho(\Pre_{\A^1,\nis}^\tfr(S)_\mathrm{gp})
.\]
Define $\SH^{\fr}(S)$ as the homotopy category of the $\PP^{\wedge 1}$-stabilisation
\[
\SH^\fr(S)=\SH^{\PP^{\wedge 1},\fr}(S)=\Ho(\Pre_{\A^1,\nis}(\Corr^\tfr(S))_\mathrm{gp}[\PP^{\wedge -1}]),
\]
and $\SH^{s,t,\fr}(S)$ as the bi-stabilisation
\[
\SH^{s,t,\fr}(S)=\Ho(\Pre_{\A^1,\nis}(\Corr^\tfr(S))_\mathrm{gp}[(S^1)^{\wedge -1}][\Gm^{\wedge -1}])
.\]
\begin{theorem}[\protect{EHKSY \cite{five-authors}}]
For a prefect field $k$ there are equivalences
\begin{equation}\label{eq:SHveffPrefrgpANisSHPrefrgpANisP}
\SH^\mathrm{veff}(k)\simeq \mathbf H^\fr_\mathrm{gp}(k),\; \SH(k)\simeq \SH^\fr(k)
\end{equation}
The canonical functors
$\Pre(\Corr^\fr(k))\to \HH^\fr(k)\to \HH^\fr_\mathrm{gp}(k)$
are given by 
$F \mapsto L_{\nis} L_{\A^1} F\mapsto L_{\nis} L_{\A^1} F^\gp$.
\end{theorem}
The advantages of tangentially framed correspondences are 
that the presheaf $h^\tfr(X)$ is an $\infty$-commutative monoid, 
and is representable in a certain $\infty$-category,
while such structures for $\Fr(-,X)$ are defined up to $\A^1$-homotopies.
This allows to apply the group-completion before the $\A^1$-localisation in the formula for the motivic replacement like \eqref{eq:LnisLA1Fr(SigmaTX)gp} above.
At the other side, the definition of $h^\tfr(X)$ is less geometrical and elementary in comparing with $\Fr(-,X)$ because it involves the $K$-theory space and the cotangent complex.

The first equivalence in \eqref{eq:SHveffPrefrgpANisSHPrefrgpANisP} 
does not hold over positive dimensional bases because of conterexamples for the cancellation theorem, see \Cref{section:ContrexampleCancellationDMeff}.
The reconstructions from \cite{BigFrmotives,FramedGamma} fail too.
Nevertheless, the second equivalence in \eqref{eq:SHveffPrefrgpANisSHPrefrgpANisP} holds for any base scheme.
\begin{theorem}[\protect{Hoyois \cite[Thorem 18]{Hoyois-framed-loc}}]
For any base scheme $S$ there is a canonical equivalence
\[
\SH(S)\simeq\SH^\fr(S).
\]
\end{theorem}
\begin{proof}[Alternative proof for noetherian separated base schemes of finite Krull dimension.]
There are equivalences 
\[\SH^{\PP^{\wedge 1}}(S)\simeq \SH^{s,t}(S), \SH^{\PP^{\wedge 1},\fr}(S)\simeq \SH^{s,t,\fr}(S),\]
where the left one is well known, and holds due to the motivic equivalence $S^1\wedge \Gm^{\wedge 1}\simeq \PP^{\wedge 1}$,
and the right one follows by the same reason.
Hence it is enough to prove the claim for bi-spectra.

Consider the adjunction
\[\gamma^*\colon \SH^{s,t}(S)\rightleftarrows \SH^{s,t,\fr}(S)\colon \gamma_*,\]
For $Y\in \Sm_B$ the unit of the adjunction
\begin{equation}\label{eq:IdGammausls}\mathrm{Id}_{\SH^{S^1,\Gm}(S)}(\Sigma^\infty_{s,t}Y) \to \gamma^*\gamma_*(\Sigma^\infty_{S^1,\Gm}Y)\end{equation}
is equivalent to the second row morphism in \eqref{eq:SigmaCorSGXLnisLAFrSigmaTX},
because of the motivic equivalence of presheaves $h^\tfr(Y)\simeq \Fr(-,Y)$ provided by \cite[Corollary 2.3.25]{five-authors}.
Hence by \Cref{citeth:GP14MotivcEq} the morphism \eqref{eq:IdGammausls} is an equivalence.
Then since the functor $\gamma_*$ is conservative,
the claim follows.
\end{proof}

The above theorem in \cite{Hoyois-framed-loc} is proven as an application of the Localisation Theorem for $\SH^\fr(B)$ from the same work and more classical Loccalisation Theorem for $\SH(B)$ that appeared originally in Voevodsky-Morel's article \cite[Theorem 3.2.21]{Morel-Voevodsky}, and was recovered later by Ayoub \cite{Ayo-Loctheorem}. 

\section{Localisation theorem for $\SH^{\fr}_{\tf}(\SmAff_S)$.}\label{section:LocalisationTheoremtf} 

In the section, we 
recall and recover 
some localisation theorems from \cite{Morel-Voevodsky,Ayo-Loctheorem,Hoyois-framed-loc,DKO:SHISpecZ,SHfrzarnis}, and prove also one version.
We prove the localisation theorem for categories 
$\HH^\fr_\tau(\SmAff_S)$ with respect to some class of topologies $\tau$
that is close, but 
simpler defined than the class considered in \cite{SHfrzarnis}.
The proof is based on the principles of \cite{DKO:SHISpecZ} or other works, but is shorter and is formally independent. 
The Loclaisation Theorem is presented here in order of preparations for the next section.

\tocless\subsection{}

Let us recall the trivial fibre topology, $\tf$-topology, defined in \cite[Definition 3.1]{DKO:SHISpecZ}.
\begin{definition}
We say that a morphism $\widetilde X\to X$ in $\Sch_S$ is a $\tf$-covering if it is \'etale affine and for each $z\in B$ there is a lifting $X\times_S z\to \widetilde X$.
We call by $\tf$-square the pullback square of the form 
\[\xymatrix{
U^\prime\ar@{^(->}[r]\ar[d] & X^\prime\ar[d]^{f}\\
U\ar@{^(->}[r]^j & X
}\]
such that 
$f$ is \'etale affine, and there is a closed subscheme $Z$ in $B$ such that
$U=X-Y$, where $Y=Z\times_S X$, $j$ is the open immersion, $Y\times_X\widetilde X\simeq Y$.
\end{definition}
As shown in \cite[\S 3]{DKO:SHISpecZ}
$\tf$-squares define regular bounded cd-structure,
and generate complete decomposible topology on $\Sch_B$. 
Since any $\tf$-suqare is a Nisnevich square, $\tf$-topology is a subtopology of the Nisnevich one.
%
According to \cite[Remakr 3.7]{DKO:SHISpecZ} 
$\tf$-topology on $\Aff_S$ over affine $S$ 
is the strongest subtopology of the Nisnevich topology that is trivial over the residue fields.

\tocless\subsection{}
Given an affine base scheme $S$ and a closed subscheme $Z\subset S$.
For a scheme $X\in \SmAff_B$ denote by
$X_Z=X\times_B Z$ the fibre product with $Z$, by
$X^h_Z=X^h_{X_Z}$ the henselisation 
of $X$ along $X_Z$, and by
$X_{S-Z}=X\times_S (S-Z)$ the open complement.
Define the subcategory $\Aff_{S,Z}$ in $\Aff_S$ spanned by schemes of the form $X^h_Z=X^h_Z$ for $X\in \Aff_S$,
and similarly define $\mathcal S_{S,Z}$ for any subcategory $\mathcal S_S$ in $\Aff_S$.

\begin{definition}
For a topology $\tau$ on $\AffSm_S$ define the topology 
$\wtau$ on $\SmAff_{S,Z}$ 
as the 
weakest topology such that the functor $\SmAff_S\to \SmAff_{S,Z}\colon X\mapsto X^h_Z$
is continuous.

For a topology $\tau$ on $\AffSm_Z$ define $\stau$ on $\SmAff_{S,Z}$ as the strongest topology such that the functor $\SmAff_{S,Z}\to \Sm_Z\colon X\mapsto X_Z$ is continuous.
\end{definition}
Note that the topology $\stau$ above is stronger then $\wtau$,
since $\Sch_S$ the base change functor $\SmAff_B\to \SmAff_Z$ is continuous with respect to $\tau$.



\tocless\subsection{}
We consider the cases of presheaves on $\Sm_S$, and the $\infty$-category of framed correspondences. We use the $\infty$-category
$\Corr^\fr(\Sch_S)$ that objects are separable noetherian $B$-schemes and mapping spaces being tangentially framed correspondences defined in \cite{five-authors}, for a subcategory $\mathcal S_S$ of $\Sch_S$ denote by $\Corr^\fr(\mathcal S_S)$ the subcategory of $\Corr^\fr(\Sch_S)$ spanned by the objects of $\mathcal S_S$.
Denote by $\Corr^\fr(\Aff_{S,Z})$ the subcategory spanned by the objects of $\Aff_{S,Z}$ in $\Corr^\fr(\Aff_S)$, 
and similarly for any subcategory $\mathcal S_*$ of $\Aff_*$.
%
Denote $\Pre^\fr(\mathcal S_*)=\Pre(\Corr^\fr(\mathcal S_*))$.

\begin{definition}
A presheaf $F\in \Pre(\SmAff_{B,Z})$ or $\Pre^\fr(\SmAff_{B,Z})$ 
is called \emph{$\A^1$-invariant}
if $F((\A^1\times X)^h_Z)\simeq F(X^h_Z)$ for all $X\in \SmAff_S$.
\end{definition}
The subcategory $\Pre_{\A^1}(\SmAff_{B,Z})$ of $\A^1$-invariant presheaves is reflective, 
and the localisation functor 
$L_{\A^1}\colon\Pre_{\A^1}(\SmAff_{B,Z})\to\Pre_{\A^1}(\SmAff_{B,Z})$ is defined by \[L_{\A^1}F(X)=F((\Delta^\bullet_S\times_S X)^h_Z),\] 
similarly for $\Pre^\fr_{\A^1}(\SmAff_{B,Z})$.


\tocless\subsection{}
We proceed with the localisation property for the homotopy categories of sheaves on $\Sm_{S,Z}$, $\Sm_{S}$, $\Sm_{S-Z}$.
\begin{lemma}\label{lm:tildei*!j**preserveA1loctauloc}
Let $\tau$ be a topology on $\SmAff_S$. 
Then the functors
\[\Pre^\fr(\SmAff_{S,Z})\rightleftarrows \Pre^\fr(\SmAff_{S}) \rightleftarrows \Pre^\fr(\SmAff_{S-Z})\]
\[\begin{array}{ll}\tilde i^!(F)(X^h_Z)=\fib (F(X^h_Z)\to F(X^h_Z-X_Z)) &
\tilde i_*(F)(X)=F(X^h_Z) \\
j^*(F)(X)=F(X) &
j_*(F)(X)=F(X_{S-Z}).
\end{array}\]
preserve 
$\tau$-sheaves in sense of $\wtau$-topology on $\SmAff_{S,Z}$.
If $\tau\supset \tf$
then the functors 
\[
\Sh_{\wtau}(\SmAff_{S,Z}) \rightleftarrows \Sh_\tau(\SmAff_S)\rightleftarrows \Sh_{\tau}(\SmAff_{S-Z})
\]
induced by the restriction
preserve $\A^1$-invariant objects.
\end{lemma}
\begin{proof}
The claim is well known for $j^*$ and $j_*$.

To prove that $\tilde i^!$ preserve $\tau$-sheaves 
consider $F\in \Sh_{\tau}(\SmAff_S)$. 
By definition any $\wtau$-covering 
in the category $\SmAff_{B,Z}$ is given by the morphism $\widetilde{X}^h_Z\to X^h_Z$
for a $\tau$-covering $\widetilde{X}\to X$ in $\SmAff_S$. 
Since $\widetilde X\to X$ and $\widetilde X_{S-Z}\to X_{S-Z}$ are $\tau$-coverings
it follows that $i^!(F)$ is $\wtau$-sheaf.
If $F$ is $\A^1$-invariant then for any $X\in \SmAff_S$
\begin{multline*}
\fib(F((\A^1\times X)^h_Z)\to F((\A^1\times X)^h_Z\times_S (S-Z))\simeq 
\fib(F(\A^1\times X)\to F((\A^1\times X)\times_S (S-Z))\simeq \\
\fib(F(X)\to F(X\times_S (S-Z))\simeq
\fib(F(X^h_Z)\to F(X^h_Z\times_S (S-Z)).
\end{multline*}
Thus $i^!F$ is $\A^1$-invariant.

To prove that claim for $\tilde i_*$ consider $F\in \Sh_{\wtau}(\Sm_{S,Z})$. 
Then $\tilde i_*(F)$ is a $\tau$-sheaf on $\SmAff_S$ because for any $\tau$-covering $\widetilde X\to X$ in $\SmAff_S$ the morphism $\widetilde X^h_Z\to X^h_Z$ is a $\tau$-covering in $\SmAff_{B,Z}$.
If $F$ is $\A^1$-invariant, then for any $X\in \SmAff_S$
$\tilde i_*F(\A^1\times X)=F((\A^1\times X)^h_Z)\simeq F(X^h_Z)=\tilde i_*F(X)$.
\end{proof}

\begin{proposition}\label{prop:HtsSmBcZSmBSmU}
Given a base scheme $S$ and a topology $\tau$ on $\Aff_S$ such that $\tf\subset\tau$, 
consider a pair of adjunctions
\[\Prefr_{\wtau,\A^1}(\SmAff_{S,Z})\rightleftarrows \Prefr_{\tau,\A^1}(\SmAff_{S}) \rightleftarrows \Prefr_{\tau,\A^1}(\SmAff_{S-Z})\]
given by $\tilde i_*\dashv \tilde i^!$, $j^*\dashv j_*$.
Then for any $F\in \Prefr_{\tau,\A^1}(\SmAff_S)$,
there is a homotopy pullback square
\[\xymatrix{
\tilde i_*\tilde i^!F\ar[d]\ar[r]& F\ar[d]\\
{*}\ar[r] & j_*j^*F
}\]
\end{proposition}
\begin{proof}
The square is pullback for any $F\in \Prefr_\tf(\SmAff_S)$, since 
\[\tilde i_*\tilde i^!F \simeq 
\mathrm{fib}(F(X^h_Z)\to F(X^h_Z-X_Z))\simeq 
\mathrm{fib}(F(X)\to F(X-X_Z))\simeq 
\mathrm{fib}(F(X)\to j_*j^*F(X)).\]
The claim for $\Prefr_{\tf,\A^1}(\SmAff_S)$ follows by \Cref{lm:tildei*!j**preserveA1loctauloc}.
%
\end{proof}

\tocless\subsection{}
Consider any of the following three presheaves of various types of framed correspondences defined in different sources: 
\begin{itemize}
\item[-] the presheaves $\Fr(-,X)$; 
\item[-] denote by $h^\nfr(X)$ the presheaves of \emph{normally} framed correspondences defined independently in \cite{five-authors} and \cite{GNThomSpectra},
and 
\item[-] denote by $h^\pfr(X)$ the presheaves $\Fr^{\mathrm{st:id}}(-,X)$ defined in \cite[Definition 7]{SmModelSpectrumTP}. 
\end{itemize}
Any three variants of the framed correspondences satisfy the lifting property with respect to affine henselian pairs.
\begin{lemma}\label{lm:henshairliftFrsmoothY}
Let $h^\fr(E)$ denote anyone of the presheaves $\Fr(-,E)$, $h^\nfr(E)$, $h^\pfr(E)$
for $E\in \SmAff_{S}$, and let $X\in \Aff\cap\Sch_S$, $Y\subset X$ be a closed subscheme.
Then the morphism 
\[ h^\fr(E)(X^h_Y)\to h^\fr(E)(Y)\]
is surjective. 
\end{lemma}
\begin{proof}
The case of the presheaf $\Fr(-,E)$ is given by \cite[Lemma A.11]{DKO:SHISpecZ}.
By \cite[Theorem 5.1.5]{five-authors}, \cite[Proposition 4]{SmModelSpectrumTP} the presheaves $h^{\nfr}(E)$ and $h^\pfr(X)$ are representable by smooth affine schemes.
Then the claim for the latter two presheaves follows by \cite[Proposition 3.10]{FrRigidSmAffpairs}. 
\end{proof}
\begin{lemma}\label{lm:henshairliftFrsmooth}
For $h^\fr(E)$, $E$, and $X$ as above, let $Y,W\subset X$ be closed subschemes.
Then the morphism 
\[ h^\fr(E)(X^h_Y)\to h^\fr(E)(Y\cup W^h_{Y\cap W})\]
is surjective. 
\end{lemma}
\begin{proof}
The claim follows from \Cref{lm:henshairliftFrsmoothY}
applied to $X$ being $X^h_Y$, and $Y$ being $Y\cup W^h_{W\cap Y}$.
\end{proof}
\begin{corollary}\label{cor:LA1trivfibhenspairsmoothFr}
Let $E\in \AffSm_{S}$.
Let $U\in \AffSm_{S}$, 
$Y\subset U$ be a closed subscheme.
Then the morphism 
\[ L_{\A^1}h^\fr(V)(U^h_Y)\to L_{\A^1}h^\fr(V)({Y})\]
is a trivial fibration.
\end{corollary}
\begin{proof}
Let $K\to N$ be an injection of simplicial sets.
Then since $h^\fr$ satisfies closed gluing,
the morphism of sets
\[ L_{\A^1}h^\fr(E^h_Z)(U^h_Z)^{N} \to 
L_{\A^1}h^\fr(E^h_Z)(U^h_Z)^{K}\times_{
L_{\A^1}h^\fr(E^h_Z)(U_Z)^{K}}
L_{\A^1}h^\fr(E^h_Z)(U_Z)^{N}
\]
equals
\[
h^\fr(E)(U^h_Y\times N_{\A^1})\to 
h^\fr(E)(U^h_Y\times K)^{}\times_{
h^\fr(E)(Y\times K)}
h^\fr(E)(Y\times N)\cong
h^\fr_{S,Z}(E)(U^h_Y\times K \amalg_{Y\times K} Y\times N).
\]
The last morphism is surjective 
by \Cref{lm:henshairliftFrsmooth} applied to 
$E=V$, $X=U\times K$, $W= U\times K \amalg_{Y\times K} Y\times N$.
\end{proof}

\begin{corollary}\label{lm:LA1hfrhenspairweakeq}
For any 
$E\in \SmAff_S$, 
$U\in \SmAff_S$,
and closed subscheme $Y\subset X$,
the morphism 
$L_{\A^1}h^\tfr(E)(U^h_Y)\to L_{\A^1}h^\tfr(E)(Y)$ 
is a weak equivalence. 
\end{corollary}
\begin{proof}
By \cite[Corollary 2.3.25]{five-authors} and by \cite[Proposition 3]{SmModelSpectrumTP} there are equivalences 
$L_{\A^1}h^\tfr(E)\simeq L_{\A^1}h^\nfr(E)\simeq L_{\A^1}\Fr(E,-)\simeq L_{\A^1}h^{\pfr}(E)$ on the category $\SmAff_S$.
So summarising $L_{\A^1}h^\tfr(E)\simeq L_{\A^1}h^\fr(E)$ on the category $\SmAff_S$, for any $h^\fr\in\{\Fr,h^\nfr,h^\pfr\}$.
The then claim follows by \Cref{cor:LA1trivfibhenspairsmoothFr}.
%
\end{proof}

\tocless\subsection{}
For any category $\mathcal S_*$ that is $\Sch_S$, $\Sm_S$, or $\Aff_{S,Z}$, $\SmAff_{S,Z}$, $\Smat_{S,Z}$ 
define the $\infty$-category $\Corr^{\A\fr}(\mathcal S_*)$ as follows.
Recall form \cite{DKO:SHISpecZ} that $\Smat_{S,Z}$ denotes the subcategory of $\SmAff_{S,Z}$ spanned by schemes with trivial tangent bundle.
\begin{definition}Consider the functor of $\infty$-categories 
$v\colon \Corr^\fr(\mathcal S_*)\to \Pre^\fr_{\A^1}(\mathcal S_*)$ 
and choose a model 
for $v$ in simplicially enriched categories. 
Then define the simplicially enriched category
$\Corr^{\A\fr}(\mathcal S_*)$ with the same objects as in the simplicially enriched category $\Corr^\fr(\mathcal S_*)$ and hom-spaces given by $\Corr^{\A\fr}(X_1,X_2)=\Pre^\fr_{\A^1}(v(X_1),v(X_2))$.
Then there is the canonical functor $l\colon \Corr^\fr(\mathcal S_*)\to \Corr^{\A\fr}(\mathcal S_*)$.
Consider 
the $\infty$-category associated with $\Corr^{\A\fr}(\mathcal S_*)$ 
and the functor of $\infty$-categories $l\colon \Corr^\fr(\mathcal S_*)\to \Corr^{\A\fr}(\mathcal S_*)$.
\end{definition}
The canonical functor $\Corr^\fr(\mathcal S_*)\to \Corr^{\A\fr}(\mathcal S_*)$ that is a kind of an $\A^1$-localisation. We use formally the following lemma.
\begin{lemma}\label{lm:HA1CorrfrCorrAfr}
The adjunction
$l^*\colon \Pre_{\A^1}(\Corr^\fr(\mathcal S_{S,Z}))\rightleftarrows \Pre_{\A^1}(\Corr^{\A\fr}(\mathcal S_{S,Z}))\colon l_*$
is an equivalence,
where $\mathcal S_{S,Z}$ is $\SmAff_{S,Z}$ or $\Smat_{S,Z}$.
\end{lemma}
\begin{proof}
By the definition $l_*h^{\A\fr}(E)\simeq L_{\A^1}h^\fr(E)$, and $l^* h^\fr(E)\simeq h^{\A\fr}(E)$.
Hence $l_*l^*h^\fr(E)\simeq L_{\A^1}h^\fr(E)$. Thus $l_*l^*\simeq_{\A^1} \Id$.

Since $l^* h^\fr(E) \simeq h^{\A\fr}(E)$, $l^*$ preserves $\A^1$-homotopy equivalences. Hence
\[l^*l_*h^{\A\fr}(E)\simeq l^*L_{\A^1}h^\fr(E)\simeq_{\A^1} l^*h^\fr(E)\simeq h^{\A\fr}(E).\]
Thus $l^*l_*\simeq_{\A^1} \Id$.
\end{proof}


Consider the commutative square
\begin{equation}\label{eq:CorrfrAfrSZZ}
\xymatrix{
\Corr^{\A\fr}(\SmAff_{S,Z})\ar[r]& \Corr^{\A\fr}(\SmAff_Z)\\
\Corr^\fr(\SmAff_{S,Z})\ar[u]\ar[r]& \Corr^\fr(\SmAff_Z)\ar[u] 
,}\end{equation}
where the horizontal arrows take $X^h_Z$ to $X_Z$, and similar one for $\Smat_*$.

\begin{lemma}\label{eq:AfrSZAsimeqZ}
The functor in the first row of \eqref{eq:CorrfrAfrSZZ}
is an equivalence.
\end{lemma}
\begin{proof}
The claim follows because morphisms
$\Corr^{\A\fr}_{S,Z}(X_1,X_2)\to \Corr^{\A\fr}_{Z}(X_1,X_2)$
are weak equivalences 
by \Cref{lm:LA1hfrhenspairweakeq}.
\end{proof}

The functor in the second row of \eqref{eq:CorrfrAfrSZZ}
induces the functors of categories
\begin{equation}\label{eq:HfrASZZ}\overline i^*\colon \Pre^\fr(\SmAff_Z)\leftrightarrows \Pre^\fr(\SmAff_{S,Z})\colon \overline i_*\end{equation} 
given by $\overline i^*(h^\tfr(X^h_Z))=h^\tfr(X_Z)$, $\overline i_*(F)(X^h_Z)=F(X_Z)$, and similarly for $\Smat_*$.
Moreover, the 
the functor $\overline i_*$ preserves $\A^1$-local objects because the horizontal functors in \eqref{eq:CorrfrAfrSZZ} take $(\A^1\times X)^h_Z$ to $\A^1\times X_Z$,

\begin{proposition}\label{prop:HfrBcZsimeqHfrZ}
The functor $\overline i_*$
induces equivalence of categories
$\PrefrA(\SmAff_{S,Z})\simeq \PrefrA(\SmAff_Z)$. 
Consequently, $\Prefr_{\stau,\A^1}(\SmAff_{S,Z})\simeq \Prefr_{\tau,\A^1}(\SmAff_Z)$ for a topology $\tau$ on $\SmAff_Z$.
The case holds for the category $\Smat_*$.
\end{proposition}
\begin{proof}


The equivalence $\Prefr_{\A^1}(\SmAff_Z)\simeq \Prefr_{\A^1}(\SmAff_{B,Z})$ follows by \Cref{lm:HA1CorrfrCorrAfr} and \Cref{eq:AfrSZAsimeqZ}.
To prove the claim for a non-trivial topology $\tau$ 
we note that 
a morphism $\widetilde X\to X$ is a $\stau$-covering if and only if
$\widetilde X_Z\to X_Z$ is a $\tau$-covering.
Then 
a presheaf $F\in \Pre^\fr(\SmAff_Z)$ is a $\tau$-sheaf if and only if 
the presheaf $F(-\times_S Z)$ is $\stau$-sheaf.
Thus the proven equivalence for the trivial topology restricts to the equivalence on the subcategories of sheaves. 
\end{proof}


\tocless\subsection{} We conclude the main result of the section.
Consider the functor $i^*\colon \Prefr_{\A^1,\stau}(\SmAff_{Z})\to\Prefr_{\A^1,\tau}(\SmAff_{S})$; $i_*F(X^h_Z)=F(X_Z)$.
Then $i_*\simeq\tilde i_*\overline i_*$, and $i_*$ has the right adjoint $i^!:=\overline i^*\tilde i^!$, where $\overline i^*=(\overline i_*)^{-1}$ is the inverse functor to the equivalence $\overline i_*\colon \Prefr_{\A^1,\tau}(\SmAff_{S})\to \Prefr_{\A^1,\stau}(\SmAff_{Z})$ given by \Cref{prop:HfrBcZsimeqHfrZ}.
\begin{theorem}\label{th:tautfLoc}
Let $Z\to S$ be a closed immersion of affine schemes.
Let $\tau$ be a topology on $\Sch_S$ such that $\tau\supset \tf$, and $\wtau=\stau$ on $\SmAff_{S,Z}$.
Consider the pair of adjunctions
\[\Prefr_{\A^1,\tau}(\SmAff_{S})\rightleftarrows \Prefr_{\A^1,\tau}(\SmAff_{S})\rightleftarrows\Prefr_{\A^1,\tau}(\SmAff_{Z})\]
given by $i_*\dashv i^!$, $j^*\dashv j_*$.
Then for any $F\in \Prefr_{\A^1,\tau}(\SmAff_{S})$, there is a homotopy pullback 
\[\xymatrix{
i_* i^! F\ar[r]\ar[d] & F\ar[d]\\
{*}\ar[r] & j_* j^* F.
}\]
\end{theorem}
\begin{proof}
The claim follows by \Cref{prop:HtsSmBcZSmBSmU} and \Cref{prop:HfrBcZsimeqHfrZ}. 
\end{proof}

\begin{example}
The example of the topologies $\tau$ above are the trivial topology, Zariski topology, and the Nisnevich topology all of that were already considered in \cite{SHfrzarnis}. Moreover, since any scheme has a covering by affine schemes, as shown in \cite{SHfrzarnis} the claim for the case of $\SmAff_S$ implies the claim for $\Sm_S$ as well.
\end{example}

\section{Improved localisation theorem}\label{section:SplittingdiagramLocalisation}

In the section, we make an improvement 
of \Cref{th:tautfLoc}. 
This is the key point for results of \Cref{subsect:stabletfmotivicdecomposition,subsect:categoricalform}

The basic principle of the argument is similar to \cite{DKO:SHISpecZ}. 
Namely, we decompose the motivic equivalence from \Cref{th:tautfLoc} into separate Nisnevich local and $\A^1$-homotopy equivalences formulated in terms of 
the category $\Sm_{S,Z}$.
The main novelty is that the results completely cover the case of motivic $(s,t)$-bi-spectra.

\tocless\subsection{}
We start with recollection of some functors from \cite{DKO:SHISpecZ}. Let $S$ be a base scheme, $Z$ be a closed subscheme.
Consider $\infty$-categories with functors
\begin{equation}\label{eq:SmZScZSmZ}
\begin{array}{lclclcl}
\Sm_Z &\xleftarrow{\overarrow{i}^{-1}}& \Sm_{S,Z} &\xleftarrow{\tilde i^{-1}}& \Sm_{S} &\xrightarrow{j^{-1}}& \Sm_{S,Z};\\
X_Z &\mapsfrom& X^h_Z &\mapsfrom& X &\mapsfrom& X\times_{S}(S-Z),
\end{array}\end{equation}
and corresponding $\infty$-categories of framed correspondences with similar functors.
Then there are functors on the $\infty$-categories of presheaves
\begin{equation}\label{eq:Smaf(Z)A1Smaf(SZ)nisSmaf(S)(U)}\xymatrix{
& & \Pre^\fr(\Sm_S)\ar[dl]^{\tilde i^!}\ar[dr]^{j^*} & \\
\Pre^\fr(\Sm_Z)\ar@<1ex>[r]^{\overarrow{i}_*} &\Pre^\fr(\Sm_{S,Z})\ar@<1ex>[l]^{\overarrow{i}^*}\ar[rd]^{\tilde i_*}&  &\Pre^\fr(\Sm_{S-Z})\ar[ld]^{j_*} \\
& & \Pre^\fr(\Sm_S) &
}\end{equation}
\begin{gather}\label{eq:PreScZiSjS-Z}
\begin{array}{ll}
\tilde i^!(F)(X^h_Z) = \cofib( F(X^h_Z)\to F(X^h_Z-X_Z) ), &X^h_Z\in \Sm_{S,Z},\\
\tilde i_*F(X)=F(X^h_Z), &X\in\Sm_S,\\
j^*F(V)=F(V), &V\in \Sm_{S-Z},\\
j_*F(X)=F(X-X_Z), &X\in \Sm_S\\
\overarrow{i}_*F(X)=F(X_Z), & X\in\Sm_{S,Z}\\
\overarrow{i}^*(h^\fr(Y^h_Z)) = h^\fr(Y_Z), 
\end{array}\end{gather}
and $\overarrow{i}^*(F)$ is the left Kan extension of $F$ along $\overarrow{i}$, for $F\in \Pre^\fr(\Sm_{S,Z})$.

\tocless\subsection{}
Let 
$\Delta_{\A^1}^\bullet$ denote co-simplicial objects
$\Delta_S^\bullet$, or $\Delta_{S,Z}^\bullet$, or $\Delta_{S-Z}^\bullet$
in categories \eqref{eq:SmZScZSmZ} depending on the context.
For a scheme $V\in \Sm_{S,Z}$ consider the endofunctor $F\mapsto F^V$ on $\Pre^\fr(\Sm_{S,Z})$, where $F^V(-)=F(-\times_{S,Z} V)$. The same notation we use for simplicial schemes, then $\Lrep_{\A^1}(F)=F^{\Delta_{\A^1}^\bullet}$.

\begin{lemma}\label{lm:(Vtimes-)taulocA1inv}
For any $S$ and $Z$ as above and for any $V\in \Sm_{S,Z}$,
the functor $(-)^V$ on $\Pre^\fr(\Sm_{S,Z})$ preserves 
$\tau$-local objects, and $\A^1$-invariant objects.
\end{lemma}
\begin{proof}
The claim follows since the endofunctor $\Sm_S\to \Sm_S;$ $U\mapsto U\times V$ preserves $\tau$-coverings and morphisms of the form $X\times\A^1\to X$.
\end{proof}


\tocless\subsection{} The main results of the section are the following two theorems relating to the right and the left part of the diagram above.
\begin{theorem}\label{th:LocA1tfnisstructuresShLocsquare}

Let $\tau$ be a topology on $\Sm_S$ stronger then $\tf$-topology
\begin{itemize}
\item[(1)]
The functors $\tilde i^!$, $\tilde i_*$, $j^*$, $j_*$ preserve $\tau$-local and $\wtau$-local objects. 
There is a canonical equivalence \[i_*i^!(F)\simeq \cofib(F\to j_*j^*F)\]
for the induced functors
\[\Sh_{\wtau}^\fr(\Sm_{S,Z})\rightleftarrows \Sh_\tau^\fr(S) \rightleftarrows \Sh_\tau^\fr(S-Z).\]
\item[(2)]
For $V\in \Sm_S$, where $V_{S-Z}=V\times_S(S-Z)$
the functors $\tilde i^!$, $\tilde i_*$, $j^*$, $j_*$ on $\Sh^\fr_\wtau$ and $\Sh^\fr_\tau$ defined in part (1)
there are equivalences
\[\begin{array}{ll}
\tilde i^! F^V\simeq (\tilde i^!F)^{V^h_Z},& 
\tilde i_* F^{V^h_Z}\simeq (\tilde i_*F)^{V},\\ 
j^* F^{V_{S-Z}}\simeq (j^*F)^V,& 
j_* F^{V_{S-Z}}\simeq (j_*F)^V; 
\end{array}
\]
in particular, the above functors commute with the functors $\Omega^l_{\Gm}\Lrep_{\A^1}$, $l\geq 0$,
on $\Sm_S$ or $\Sm_{S,Z}$, or $\Sm_{S-Z}$ respectively.
\item[(3)]
For any $F\in \Pre^\fr_{\A^1}$ there is an equivalence
\[\tilde i^! \Lrep_{\nis} F\simeq \Lrep_{\nis} \tilde i^! F\]
and similarly for 
$\tilde i_*$, $j^*$, $j_*$. 
\end{itemize}

\end{theorem}
\begin{proof}
By \Cref{lm:wtaulocaobjectsShLocsquare} 
functors $\tilde i^!$ $\tilde i_*$ $j^*$ $j_*$ preserve $\tau$-local and $\wtau$-local objects.
For any 
$F\in \Pre_\tau^\fr(\Sm_S)$
and $X\in \Sm_S$,
there is a canonical fibre sequence
\[
\fib( F(X) \to F(X\times_S (S-Z)) ) \to F(X)\to F(X\times_S (S-Z))
\]
and the equivalence
\[F(X\times_S (S-Z))\simeq j_* j^* F(X).\]
For any $F\in \Sh_\tau^\fr(\Sm_S)$,
there are canonical equivalences
\[
\tilde i_* \tilde i^! F(X)\simeq
\fib( F(X^h_Z) \to F(X^h_Z\times_S (S-Z)) )\simeq 
\fib( F(X) \to F(X\times_S (S-Z)) )
.\]
Thus 
there is a natural fibre sequence
\[\tilde i_* \tilde i^! F(X)\to F(X)\to j_* j^* F(X).\]
This proves Part (1).

\Cref{lm:F(Vtimes-)ShLocsquare} gives the claim for the (endo)functors $(-)^V$.
The case of $\Omega_{\Gm}\Lrep_{\A^1}$
follows because $\Lrep_{\A^1} F(X)\cong F(\Delta_{\A^1}^\bullet\times -)$, 
where $\Delta_{\A^1}^\bullet=\Delta_S^\bullet$, or $\Delta_{S,Z}^\bullet$, or $\Delta_{S-Z}^\bullet$,
and $\Omega_{\Gm} F\cong \fib ( F(\Gm\times -)\to F(\{1\}\times -) )$.
This completes Part (2).

Part (3) holds by
\Cref{lm:tiusjustiufjdsLnis}.
\end{proof}

\begin{theorem}\label{th:LocA1tfnisstructuresDeformation}

Let $\tau$ be a topology on $\Sm_Z$.
The following holds: 
\begin{itemize}
\item[(1)]
functors $\overarrow{i}^*$ and $\overarrow{i}_*$ on $\Pre^\fr$ preserve $\A^1$-invariant objects and take $\tau$-local objects to $\stau$-local objects; 
the induced functors on the subcategories 
\begin{equation}\label{lm:overlineistauA1}\overarrow{i}_*^{\stau,\A^1}\colon \Pre^\fr_{\tau,\A^1}(Z)\simeq \Pre^\fr_{\stau,\A^1}(\Sm_{S,Z})\colon \overarrow{i}^*_{\stau,\A^1};\end{equation}
are equivalences;
\item[(2)] $\overarrow{i}^*$ is exact and conservative with respect to $\A^1$-homotopy equivalences;
\item[(3)]
$\overarrow{i}_*$ 
is exact and conservative with respect to Nisnevich local equivalences on $\Pre^\fr$;
$\overarrow{i}^*$ 
is Nisnevich exact on $\Pre^\fr$,
and is conservative with respect to Nisnevich local equivalences on the subcategory $\Pre^\fr_{\A^1}$;
\item[(4)]
both functors $\overline{i}_*$ and $\overarrow{i}^*$ 
commute with the functor $\Omega_{\Gm}\Lrep_{\A^1}$ on $\Pre^\fr$.
\end{itemize}

\end{theorem}
\begin{proof}
The claim in part (1) on $\A^1$-invariant objects is provided by \Cref{prop:iAGm}(2,3).
We proceed with the case of $\stau$-local objects.
The functor $\overarrow{i}_*$ takes $\tau$-local objects to $\stau$-local because the functor $\SmAff_{S,Z}\to \SmAff_Z;$ $X^h_Z\mapsto X_Z$ takes $\stau$-coverings to $\tau$-coverings.
In follows by the definition of $\stau$-topology that
for any $\tau$-covering in $\SmAff_Z$ there is a refinement $\widetilde X\to X$ and 
an $\stau$-covering $\widetilde Y\to Y$ such that $\widetilde X=\widetilde Y_Z$, $X=Y_Z$. 
Hence
the functor $\overarrow{i}^*$ takes $\stau$-local objects to $\tau$-local.

The part (4) is the part (1) of \Cref{prop:iAGm}. 

Since $i_*$ commutes with $\Lrep_{\A^1}$, it is exact with respect to $\A^1$-homotopy equivalences, and since
it induces equivalences on subcategories of $\A^1$-invariant objects in addition, it is conservative with respect to $\A^1$-homotopy equivalences. This completes Part (2).

Part (3) is by \Cref{prop:ovidsus(simeqNis)}.
\end{proof}








\subsection{Deformation $\Sm_Z\leftarrow \Sm_{S,Z}$}\label{sect:DeformationZSZ}

\begin{lemma}\label{lm:iusA1invariant}
The functor $\overline i^*\colon \Pre(\SmAff_{S,Z})\to \Pre(\SmAff_Z)$
preserves $\A^1$-invariant objects.
\end{lemma}
\begin{proof}
Let $X\in \Sm_S$.
Then morphisms $X_Z\cong X_Z\times 0\to X_Z\times \A^1$
induce morphisms 
\[\begin{array}{lrll}
i_X^*\colon& i^* F((X\times\A^2)_Z)&\to& i^* F(X_Z),\\ p^*_X\colon& i^* F(X_Z)&\to& i^* F((X\times\A^2)_Z)
.\end{array}\]
Immediately, $i_X^*p^*_X\simeq \id_{i^* F(X_Z)}$.
We are going to show that $p^*_X i^*_X\simeq \id_{i^* F(X_Z\times\A^2)}$.

Consider the comma category ${X_Z}\downarrow\Sch_S$ of $\Sch_S$ under $X_Z$, and the subcategory 
spanned by the objects in ${X_Z}\downarrow\Sch_S$ given by morphisms of schemes $X_Z\to Y^h_Z$ for $Y\in \SmAff_S$. The latter category we call the \emph{indexing category} and denote by ${X_Z}\downarrow\Sm_{S,Z}$.
Similarly consider the indexing category ${(X\times\A^2)_Z}\downarrow\Sm_{S,Z}$ 
spanned by the objects given by $(X\times\A^2)_Z\to (Y^\prime)^h_Z$
in the comma category ${(X\times\A^2)_Z}\downarrow\Sch_S$.
Then 
\[\begin{array}{lll}
i^*F(X_Z)&\simeq& \varinjlim_{X_Z\to Y^h_Z, Y\in \Sm_S} F(Y),\\
i^*F((X\times\A^1)_Z)&\simeq& \varinjlim_{(X\times\A^1)_Z\to Y^h_Z, Y\in \Sm_S} F(Y),
\end{array}\]
where the injective limits are over the indexing categories ${X_Z}\downarrow\Sm_{S,Z}$ and ${(X\times\A^2)_Z}\downarrow\Sm_{S,Z}$.
Morphisms $p^*_X$ and $i^*_X$ are induced by
\[\begin{array}{lll}\varinjlim_{(X\times \A^2)_Z\to Y^h_Z} F(Y) &\to& \varinjlim_{X_Z\to Y^h_Z} F(Y),\\
\varinjlim_{X_Z\to Y^h_Z} F(Y) &\to& \varinjlim_{(X\times\A^2)_Z\to (Y^\prime)^h_Z} F(Y^\prime)\end{array}\]
where the first morphism is induced by the functor of the indexing categories
\[{(X\times\A^2)_Z}\downarrow\Sm_{S,Z}\to {X_Z}\downarrow\Sm_{S,Z}; [(X\times \A^2)_Z\to Y^h_Z] \mapsto [X_Z\cong (X\times 0)_Z\to Y^h_Z]\]
and the second one is induced by the functor 
\[{X_Z}\downarrow\Sm_{S,Z}\to {(X\times\A^2)_Z}\downarrow\Sm_{S,Z}; [X_Z\to Y^h_Z] \mapsto [(X\times\A^2)_Z\to (Y\times\A^2)^h_Z]\]
and natural equivalences $F(Y^h_Z)\to F((Y\times\A^2)^h_Z)$ provided by that $F$ is $\A^1$-invariant on $\Sm_{S,Z}$.

Then the composite $p^*_X i^*_X$
is induced by the endofunctor on $(X\times\A^2)_Z\downarrow \Sm_{S,Z}$ that takes
$(X\times \A^2)_Z\to Y^h_Z$ 
to object given by the composite morphism of schemes
\begin{equation}\label{eq:X0A1->YA1}
(X\times \A^2)_Z\cong ((X\times 0)\times \A^2)_Z\to ((X\times \A^2)\times \A^2)_Z\to ((Y)\times\A^2)^h_Z 
.\end{equation}
and the natural equivalences
$F(Y^h_Z)\to F((Y\times\A^2)^h_Z)$.
The composite morphisms of the sequences of the form \eqref{eq:X0A1->YA1}
are $\A^1$-homotopy equivalent to
the composite morphisms 
\[(X\times \A^2)_Z\cong ((X\times \A^2)\times 0)_Z \to ((X\times \A^2)\times \A^2)_Z\to (Y\times\A^2)^h_Z.\]
The latter morphism equals to the composite of
\begin{equation}\label{eq:XA1->Y0->YA1}
(X\times \A^2)_Z\to Y^h_Z\cong (Y\times 0)^h_Z\to (Y\times\A^2)^h_Z.
\end{equation}
Since $F$ is $\A^1$-invariant the latter morphism in composition \eqref{eq:XA1->Y0->YA1} 
induces the canonical equivalence 
$F((Y\times 0)^h_Z)\simeq F((Y\times\A^2)^h_Z)$.
Hence the endofunctor of the indexing category
that takes 
$(X\times \A^2)_Z\to Y^h_Z$ to the composite of \eqref{eq:XA1->Y0->YA1}
induces the autoequivalence of $i^* F((X\times\A^2)^h_Z)$.

Thus for each $X\in \Sm_S$ the canonical projection induces the equivalence 
$i^* F((X\times\A^2)^h_Z)\simeq i^* F(X^h_Z)$. Hence $i^* F$ is $\A^1$-invariant.
\end{proof}


\begin{proposition}\label{prop:iAGm}
The following holds for the functors $\overarrow{i}^*\colon \Pre^\fr(\mathcal S_{S,Z})\to \Pre^\fr(\mathcal S_{Z})\colon \overarrow{i}_*$,
where $\mathcal S_{*}=\Sm_{*},\SmAff_{*},\Smat_{*}$:
\begin{itemize}
\item[(1)]
$\overarrow{i}^*$ and $\overarrow{i}_*$ commute with $\Lrep_{\A^1}$ and $\Omega^i_{\Gm}\Lrep_{\A^1}$, for $ i\geq 0 $,
\item[(2)]
$\overarrow{i}^*$ and $\overarrow{i}_*$ preserve and detect $\A^1$-local objects, 
\item[(3)]
$\overarrow{i}^*$ and $\overarrow{i}_*$ induce the equivalence on the subcategories of $\A^1$-local objects
$\Prefr_{\A^1}(\Sm_{S,Z})\simeq \Prefr_{\A^1}(Z)$.
\end{itemize}
\end{proposition}
\begin{proof}

The functor $\overarrow{i}_*$ commutes with the functors 
$\Omega^i_{\Gm}\Lrep_{\A^1}$, for $ i\geq 0 $,
because 
$(\Gm\times\Delta_B^n)^h_Z\times_B Z=(\Gm\times\Delta_B^n)_Z$.
Since $\overarrow{i}_*$ commutes with $\Lrep_{\A^1}$, it preserves $\A^1$-invariant objects.

Consider the commutative diagram
\begin{equation}\label{eq:CorrAfrfrSZZ}
\xymatrix{
\Corr^{\A\fr}(\Sm_{B,Z})\ar[r]^{\overarrow{i}^{\A}} & \Corr^{\A\fr}(\Sm_{Z}) \\
\Corr^{\fr}(\Sm_{B,Z})\ar[u]^{l_{B,Z}}\ar[r]^{\overarrow{i}} & \Corr^{\fr}(\Sm_{Z})\ar[u]^{l_Z} 
.}\end{equation}
By \Cref{eq:AfrSZAsimeqZ} $\Corr^{\A\fr}(\Sm_{B,Z})\simeq \Corr^{\A\fr}(\Sm_{Z})$,
hence the first row in the induced diagram 
\begin{equation}\label{eq:HHldsids}
\xymatrix{
\Pre(\Corr^{\A\fr}(\Sm_{B,Z}))\ar[d]^{l_*} & \Pre(\Corr^{\A\fr}(\Sm_{Z}))\ar[d]^{l_*}\ar[l]^{\overarrow{i}_*^{\A}} \\
\Pre(\Corr^{\fr}(\Sm_{B,Z})) & \Pre(\Corr^{\fr}(\Sm_{Z}))\ar[l]^{\overarrow{i}_*} 
}\end{equation}
is an equivalence.
By \Cref{lm:lAsimeq} the first row 
is equivalent to 
\begin{equation}\label{eq:HHAid*simeq}\Pre_{\A^1}(\Sm_{B,Z}) \xleftarrow{\simeq} \Pre_{\A^1}(\Sm_{Z}),\end{equation}
and the vertical morphisms in \eqref{eq:HHldsids} are equivalent to the canonical embeddings.
Moreover, since $\overarrow{i}_*$ induces equivalence on the subcategories of $\A^1$-invariant objects, it detects $\A^1$-invariant objects.
Thus the claim of points (1-3) is proven for $\overarrow{i}_*$.

Consider the commutative diagram
\begin{equation*}
\xymatrix{
\Pre(\Corr^{\A\fr}(\Sm_{B,Z}))\ar[r]^{\overarrow{i}^*_{\A}} & \Pre(\Corr^{\A\fr}(\Sm_{Z})) \\
\Pre(\Corr^{\fr}(\Sm_{B,Z}))\ar[u]^{l^*_{B,Z}}\ar[r]^{\overarrow{i}^*} & \Pre(\Corr^{\fr}(\Sm_{Z}))\ar[u]^{l^*_Z} 
}\end{equation*}
provided by \eqref{eq:CorrAfrfrSZZ}.
Since the first row above is equivalent to the inverse for the founcotr \eqref{eq:HHAid*simeq},
$\overarrow{i}^*$ preserves $\A^1$-homotopy equivalences. By \Cref{lm:iusA1invariant} $\overarrow{i}^*$ preserves $\A^1$-invariant objects.
Hence $\overarrow{i}^*$ commutes with $\Lrep_{\A^1}$,
and the induced morphism on the subcategories of $\A^1$-invariant objects is inverse to \eqref{eq:HHAid*simeq}. This proves points (2-3).
To prove (1) we note that since $\overarrow{i}^{\A}_*$ commutes with $\Omega_{\Gm}$, the same holds for $\overarrow{i}_{\A}^*$.
%

\end{proof}
\begin{lemma}\label{lm:lAsimeq}
The functor 
$l_*\colon \Pre^{\A\fr}(\Sm_{B,Z})\to \Pre^\fr(\Sm_{B,Z})$ lands in $\Pre^\fr_\A(\Sm_{B,Z})$, 
and induces the equivalence \[\Pre^{\A\fr}(\Sm_{B,Z})\simeq \Pre^\fr_\A(\Sm_{B,Z}).\]
In addition, there are equivalences $l_*l^*\simeq \Lrep_{\A}$, $l^*l_*\simeq \Id$.
\end{lemma}
\begin{proof}
In view of the definition of the functor
$\Corr^{\fr}(\Sm_{B,Z})\xrightarrow{l} \Corr^{\A\fr}(\Sm_{B,Z})$
there is an equivalence
$l_*(h^{\A\fr}(E))\simeq L_{\A^1}^\fr(E)$
for any $E\in \Sm_{B,Z}$.
So $l_*$ takes representable presheaves to $\A^1$-invariant ones.
Since any presheaf is a colimit of prepresentable ones, the functor $l_*$ preserves colimits, and the subcategory of $\A^1$-invariant presheaves $\Pre^\fr_\A(\Sm_{B,Z})$ is closed under colimits, the first claim follows.

On other side by the definition it follows that 
$l^*(h^\fr(E))\simeq h^{\A\fr}(E)$.
So $l_*l^*(h^\fr(E))\simeq L_{\A^1}h^{\fr}(E)$
Thus $l_*l^*\simeq L_{\A^1}$ is equivalent to identity on $\Pre^\fr_{\A^1}$.  

Since $l$ is essentially surjective on objects, $l_*$ is conservative. 
More accurately let us write that if $l_*F\simeq l_*G$ then 
\[\begin{array}{lll}
F(E)\simeq& \Pre^{\A1\fr}(h^{\A\fr}(E),F)\simeq& \Pre^\fr(h^\fr(E),l_*F),\\
G(E)\simeq&& \Pre^\fr(h^\fr(E),l_*(G))
\end{array}\] 

To prove the last claim consider the morphism $F\to l^*l_*F$. By the above 
\[l_*(F)\to l_*(l^*l_*F)\simeq l_*l^*(l_*F)\simeq l_*F,\] hence $F\simeq l^*l_*F$ since $l_*$ is conservative.
%
%
\end{proof}
\begin{proposition}\label{prop:ovidsus(simeqNis)}
The functor $\overarrow{i}_*\colon \Pre^\fr(\mathcal S_Z)\to \Pre^\fr(\mathcal S_{S,Z})$
preserves and is conservative with respect to Nisnevich local equivalences. 
The functor $\overarrow{i}^*\colon \Pre^\fr(\Sm_{S,Z})\to \Pre^\fr(\Sm_{Z})$
preserves Nisnevich local equivalences and is conservative with respect to 
Nisnevich local equivalences on the subcategory of $\A^1$-invariant objects. 
Here $\mathcal S_{*}=\Sm_{*},\SmAff_{*},\Smat_{*}$.
\end{proposition}
\begin{proof}

Let $F\simeq_{\nis} G$ be a Nisnevich local equivalence in $\Pre^\fr(\mathcal S_Z)$. 
Let $X\in \Sm_S$, then the morphism $\overarrow{i}_*F(X^h_x)\to \overarrow{i}_*G(X^h_x)$ is equivalent to the morphism $F((X^h_x)_Z)\to G((X^h_x)_Z)$ that is an equivalence since $(X^h_x)_Z\cong (X_Z)^h_x$ is local henselian.

Given a morphism $F\to G$ that maps to a Nisnevich local equivalence $\overarrow{i}_*F\simeq_{\nis} \overarrow{i}_*G$.
Consider a local henselian essentially smooth scheme $X^h_x$ then there is a scheme $\widetilde X\in \Sm_S$ with a morphism $x\to \widetilde X$ such that 
$(\widetilde X_x)_Z\simeq X_x$.
Hence $(\widetilde X^h_x)_Z\simeq X^h_x$.
Then the morphism $F(X^h_x)\to G(X^h_x)$ is equivalent to the morphism $\overarrow{i}_*F(\widetilde X^h_x)\simeq \overarrow{i}_*G(\widetilde X^h_x)$.
So the claim is proven for $\overarrow{i}_*$.

%
Since $\overarrow{i}^*$ preserves Nisnevich squares, according to \Cref{def:PrefrNislocsimeqNis} and \Cref{th:restrfrLocniscommute} $\overarrow{i}^*$ preserves Nisnevich local equivalences.
The last claim follows because of \Cref{prop:iAGm}.(3), and because $\overarrow{i}_*$ preserves Nisnevihc local equivalences.
\end{proof}



\begin{lemma}
There are a canonical equivalences
\[\Pre_{\A^1}(\SmAff_{S,Z})\simeq \Pre_{\A^1}(\Smat_{S,Z}), \; \Pre^\fr_{\A^1}(\SmAff_{S,Z})\simeq \Pre^\fr_{\A^1}(\Smat_{S,Z}).\]
\end{lemma}
\begin{proof}
For any scheme $X^h_Z\in \SmAff_{S,Z}$ there is an $\A^1$-equivalent scheme $(X^\prime)^h_Z\in \Smat_{S,Z}$.
Namely, if $X^\prime$ is the total space of the vector bundle $N$ over $X$ such that $N\oplus T_X$ is trivial; since in this case $T_{X^\prime}$ is the inverse image of $N\oplus T_X$ along $X^\prime\to X$.
Hence the functors
\begin{equation}\label{eq:HSmAfftoHSmat}\Pre_{\A^1}(\SmAff_{S,Z})\to \Pre_{\A^1}(\Smat_{S,Z}), \Pre^\fr_{\A^1}(\SmAff_{S,Z})\to \Pre^\fr_{\A^1}(\Smat_{S,Z})\end{equation}
are conservative.
Since $\Smat_{S,Z}\to \SmAff_{S,Z}$ is fully faithful,
the composite functors
\[\begin{array}{lclcl}
\Pre_{\A^1}(\Smat_{S,Z})&\to& \Pre_{\A^1}(\SmAff_{S,Z})&\to& \Pre_{\A^1}(\Smat_{S,Z}),\\ 
\Pre^\fr_{\A^1}(\Smat_{S,Z})&\to& \Pre^\fr_{\A^1}(\SmAff_{S,Z})&\to& \Pre^\fr_{\A^1}(\Smat_{S,Z})
\end{array}\]
are equivalent to the identity.
Then since \eqref{eq:HSmAfftoHSmat} is conservative, the composite functors
$\Pre_{\A^1}(\SmAff_{S,Z})\to \Pre_{\A^1}(\Smat_{S,Z})\to \Pre_{\A^1}(\SmAff_{S,Z})$,
$\Pre^\fr_{\A^1}(\SmAff_{S,Z})\to \Pre^\fr_{\A^1}(\Smat_{S,Z})\to \Pre^\fr_{\A^1}(\SmAff_{S,Z})$
are equivalent to the identity too.
\end{proof}


\subsection{Localisation theorem for $\tf$-sheaves}

Let $\tau$ be a topology on $\Sm_S$ stronger then $\tf$-topology.
Recall that $\wtau$ is the weakest topology on $\Sm_{S,Z}$ such that the functor $\Sm_S\to \Sm_{S,Z}$; $X\mapsto X^h_Z$ 
is continuous with respect to the topology $\tau$ on $\Sm_S$.
To uniform notations we use the symbol $\wtau$ for the topology $\tau$ on categories $\Sm_S$ and $\Sm_{S-Z}$ too.

\begin{lemma}\label{lm:wtaulocaobjectsShLocsquare}
The functors $\tilde i^!$, $\tilde i_*$, $j^*$, $j_*$
on $\Pre^\fr(\Sm_S)$, $\Pre^\fr(\Sm_{S-Z})$, and $\Pre^\fr(\Sm_{S,Z})$ from \eqref{eq:PreScZiSjS-Z}
preserve $\wtau$-local objects.
\end{lemma}
\begin{proof}
We repeat the argument of \cite[]{DKO:SHISpecZ} in the setting of $\infty$-categories $\Corr^\fr(-)$.
The claim for $\tilde i_*$, $j^*$ and $j_*$ follows
because 
the functors 
\begin{equation}\label{eq:XhZX(S-Z)X}\begin{array}{ll}
\Sm_S\to \Sm_{S,Z}; &X\mapsto X^h_Z;\\
\Sm_S\to \Sm_{S-Z}; &X\mapsto X_{S-Z};\\
\Sm_{S-Z}\to \Sm_S; &X\mapsto X;\\
\end{array}\end{equation}
preserve $\wtau$-coverings.

Let $F\in \Sh^\fr_\tau(\Sm_S)$.
By the definition any $\wtau$-covering admits a refinement $\widetilde X^h_Z\to X^h_Z$
induced by $\tau$-covering $\widetilde X\to X$.
Then $\widetilde X_{S-Z}\to X_{S-Z}$ is $\tau$-covering too.
So since $\tau$ is stronger then $\tf$-topology
there are equivalences
\begin{align*}
\check{C}_{\widetilde X^h_Z}(X^h_Z,\tilde i^! F)&\simeq
\fib( \check{C}_{\widetilde X^h_Z}(X^h_Z,F)\to \check{C}_{\widetilde X^h_Z-X_Z}(X^h_Z-X_Z,F) )\\& \simeq 
\fib( \check{C}_{\widetilde X}(X,F)\to \check{C}_{\widetilde X-X_Z}(X-X_Z,F) )\\& \simeq 
\fib( F(X)\to F(X-X_Z) )\\& \simeq 
\fib( F(X^h_Z)\to F(X^h_Z-X_Z) )\\& \simeq 
\tilde i^! F(X^h_Z).
\end{align*}
\end{proof}

\begin{lemma}\label{lm:tiusjustiufjdsLnis}
The functors  $\tilde i_*$, $j^*$ on 
$\Pre^\fr$
commute with $\Lrep_{\nis}$. 
For any $F\in \Pre^\fr_{\A^1}$ 
\[\tilde i^! \Lrep_{\nis} F\simeq \Lrep_{\nis} \tilde i^! F,
j_* \Lrep_{\nis} F\simeq \Lrep_{\nis} j_* F.\]
\end{lemma}
\begin{proof}
By \Cref{lm:wtaulocaobjectsShLocsquare}
$\tilde i^!$ and $j_*$ preserve Nisnevich local objects
hence $\tilde i^! \Lrep_\nis F$ and $j_* \Lrep_\nis F$ are Nisnevich local.
So it is enough to prove that 
$\tilde i^! \Lrep_\nis F\to \Lrep_\nis \tilde i^! F$ and
$j_* \Lrep_\nis F\to \Lrep_\nis j_* F$
induces equivalences on the Nisnevich topology points.
A point is given by the scheme $X^h_x$ for $X\in \Sm_S$, $x\in X_Z$, then
\begin{gather*}
\tilde i^! F(X^h_x)\simeq
\fib ( F(X^h_x)\to F((X^h_x)_{S-Z}) ),\\
\tilde i^! \Lrep_\nis F(X^h_x)\simeq
\fib ( \Lrep_\nis F(X^h_x)\to \Lrep_\nis F((X^h_x)_{S-Z}) )\stackrel{}{\simeq}
\fib ( F(X^h_x)\to F((X^h_x)_{S-Z}) ),
\end{gather*}
and
\begin{gather*}
\Lrep_\nis j_* F(X^h_x)\simeq
j_* F(X^h_x)\simeq
F((X^h_x)_{S-Z}),\\
j_* \Lrep_\nis F(X^h_x)\simeq
\Lrep_\nis F((X^h_x)_{S-Z}) \stackrel{}{\simeq}
F((X^h_x)_{S-Z}).
\end{gather*}
If $X\in \Sm_S$, $x\in X_{S-Z}$, then $j_* F(X^h_x)\simeq F(X^h_x)$, and $j_* \Lrep_\nis F(X^h_x)\simeq \Lrep_\nis F(X^h_x)\simeq F(X^h_x)$.

Since by \Cref{lm:wtaulocaobjectsShLocsquare}
functors
$\tilde i_*$, $j^*$
preserve Nisnevich local objects,
it is enough to prove that
they preserve Nisnevich local equivalences.
The claim follows because
the first and the third functors in \eqref{eq:XhZX(S-Z)X}
preserve essentially smooth local henselian schemes.
\end{proof}

\begin{lemma}\label{lm:F(Vtimes-)ShLocsquare}
Functors $\tilde i^!$, $\tilde i_*$, $j^*$, $j_*$ on $\Sh^\fr_\wtau$ and $\Sh^\fr_\tau$ given by \Cref{lm:wtaulocaobjectsShLocsquare}
commute with the functors $(-)^V$, for $V\in \Sm_S$ or $\Sm_{S,Z}$, $\Sm_{S-Z}$.
\end{lemma}
\begin{proof}
The case of $\tilde i_*$, $j^*$, $j_*$ follows because
the functors \eqref{eq:XhZX(S-Z)X} commute with the functor $X\mapsto X\times V$.
The case of $\tilde i^!$ follows because
\begin{align*}
\tilde i^!F( (X\times V)^h_Z )&\simeq
\fib( F( (X\times V)^h_Z )\to F( (X\times V)^h_Z-(X\times V)_Z ) )\\ &\simeq
\fib( F( X\times V )\to F( (X\times V)-(X\times V)_Z ) )\\& \simeq
\tilde i^!F^V( X^h_Z ).
\end{align*}
\end{proof}

\section{Cone theorem}\label{section:ConeThoreom}

Consider the category of smooth open pairs 
$\Sm^\mathrm{pair}_S$
with objects given by 
$\langle X,V
U\rangle$, where $X\in \Sm_S$, and $U$ is an open subscheme;
we repeat the definition of framed correspondences according to the above notation.
\begin{definition}
A \emph{Nisnevich framed correspondence} of level $n$ between pairs $\langle X,U\rangle,\langle Y,V\rangle\in \Sm^\mathrm{pair}_S$ 
is given by the data $(Z,\phi,g)$:\begin{itemize}
\item[(i)] a closed subscheme $Z\subset \A^n\times X$ finite over $X$, and such that $Z\times_X U=\emptyset$.
\item[(ii)] a set of regular functions $\varphi_i\in \calO( (\A^n\times X)^h_Z )$, $i=1,\dots n$, 
and a regular map $g\colon (\A^n\times X)^h_Z\to Y$
such that
$(\A^n\times X)^h_Z  \times_{\varphi_1,\dots \varphi_n,g} (0\times (Y\setminus V)) = Z$.
\end{itemize}
We define $W=Z(\varphi)$ and
denote the Nisnevich framed correspondence as above by 
$(Z,\phi,g)=(Z,W,\phi,g)$.
Note that the data $W$ is extra.
\end{definition}
Denote by 
$\Fr(\langle X,U\rangle ,\langle Y,V\rangle )=\Fr(\langle Y,V\rangle)(\langle X,U\rangle)$ the set of Nisnevich framed correspondences from $\langle X,U\rangle$ to $\langle Y,V\rangle$.
We consider here only the case of $\langle X,\emptyset\rangle$ and denote by 
$\Fr(\langle Y,V\rangle )=\Fr(-,\langle Y,V\rangle )$ the presheaf on $\Sm_S$.

The following result follows from a set of equivalences proven in the next two subsections.
\begin{theorem}\label{th:Fr(Y/V)}
For $\langle Y,V\rangle\in \Sm^\mathrm{pair}_S$ 
the canonical morphism
of pointed presheaves
$\Fr(X)/\Fr(U) \to \Fr(\langle X/U\rangle)$
is a motivic equivalence.
\end{theorem}
\begin{proof}

The claim follows because of the sequence of motivic equivalences
\[\begin{array}{lll}
\Fr(X)/\Fr(U)& 
\stackrel{}{\simeq}& \cofib ( \Fr(U) \to \Fr(X) )\\
&\stackrel{Lm. \ref{lm:FraffFr}}{\simeq}& \cofib ( h^\fraff(U) \to h^\fraff(X) )\\
&\stackrel{Cor. \ref{cor:frafftghfr}}{\simeq}& \cofib ( h^\tgfraff(U) \to h^\tgfraff(X) )\\
&\stackrel{Lm. \ref{lm:tghtg(Y/V)}}{\simeq}& h^\tgfraff(\langle X/U\rangle )\\
&\stackrel{Prop. \ref{prop:Y/VtgfrCMon}}{\simeq}& h^\tghfraff(\langle X/U\rangle )\\
&\stackrel{Prop. \ref{prop:freq:Nr=Tg}}{\simeq}& h^\lownrfraff(\langle X/U\rangle )\\
&\stackrel{Prop. \ref{prop:freq:Nis=Nr=Id}}{\simeq}& h^\fraff(\langle X/U\rangle )\\
&\stackrel{Lm. \ref{lm:FraffFr}}{\simeq}& \Fr(\langle X/U\rangle ),
\end{array}\]
where the presheaves 
$h^\fraff$, $h^\tgfraff$, $h^\tghfraff$, $h^\lownrfraff$
are defined in
\Cref{def:Fraff,def:tgfr,def:tghfr,def:lownr}.
\end{proof}

\subsection{Nisnevich local equivalences}

\begin{definition}\label{def:Fraff}
A \emph{henselian framed affine correspondence} of level $n$ between pairs
 $\langle X,U\rangle,\langle Y,V\rangle\in \Sm^\mathrm{pair}_S$, 
is given by the data $(Z,\phi,g)$:\begin{itemize}
\item[(i)] a closed subscheme $Z\subset \A^n\times X$ finite over $X$, and such that $Z\times_X U=\emptyset$.
\item[(ii)] a set of regular functions $\varphi_i\in \calO( (\A^n\times X)^h_Z )$, $i=1,\dots n$, 
and a morphism of schemes $g\colon (\A^n\times X)^h_Z\to Y$ 
that decomposes as \[(\A^n\times X)^h_Z\to Y^\prime\to Y\] for $Y^\prime\in \Sm_S\cap \Aff$ and 
such that
\[(\A^n\times X)^h_Z  \times_{\varphi_1,\dots \varphi_n,g} (0\times (Y\setminus V)) = Z.\]
\end{itemize}
We define $W=Z(\varphi)$ and
denote the above framed correspondence by 
$(Z,\phi,g)=(Z,W,\phi,g)$.
Note that the data $W$ is extra.

\end{definition}
Let $\Fr^{\mathrm{aff}}_n(\langle X,U\rangle ,\langle Y,V\rangle )$ denote the set of henselian framed affine correspondences of level $n$.
Similarly to \Cref{def:FrCor} define the set 
\[\Fr^{\mathrm{aff}}(\langle X,U\rangle ,\langle Y,V\rangle )=\varinjlim_n \Fr^{\mathrm{aff}}_n(\langle X,U\rangle ,\langle Y,V\rangle ),\] 
and the presheaf 
$h^\fraff(\langle Y,V\rangle )=\Fr^{\mathrm{aff}}(-,\langle Y,V\rangle )$.

\begin{lemma}\label{lm:FraffFr}
For any $\langle Y,V\rangle\in \Sm^\mathrm{pair}_S$,
the canonical morphism
\[h^\fraff(\langle Y/V\rangle)\to \Fr(\langle Y/V\rangle)\]
of presheaves on $\Sm_S$ is a Nisnevich local equivalence.
\end{lemma}
\begin{proof}
By the definition for any $X$ 
the set $\Fr^\mathrm{aff}_n(X,\langle Y/V\rangle)$
is a subset in $\Fr_n(X,\langle Y/V\rangle)$.
To prove the Nisnevich equivalence 
we need to show the isomorphism
$\Fr_n(X,\langle Y/V\rangle)\cong \Fr^\mathrm{aff}_n(X,\langle Y/V\rangle)$
for a local henslian scheme $X$.
Given $c=(Z,\varphi,g)\in \Fr_n(X,\langle Y/V\rangle)$.
Since $Z$ is finite over $X$, 
it is corpoduct of local henselian schemes. 
Then $(\A^n_X)^h_Z= \coprod_{z\in Z_{0}}(\A^n_Z)^h_z$, where $z$ rans over the set of closed points of $Z$.
For each closed point $z$ in $Z$, choose an affine open neighbourhood $Y^\prime_{g(z)}$ of $g(z)$ in $Y$, and let 
$Y^\prime=\coprod_{z\in Z} Y^\prime_{g(z)}$.
Then the morphism $g\colon (\A^n_X)^h_Z\to Y$ passes through $Y^\prime$.
Thus $c\in \Fr^\mathrm{aff}_n(X,\langle Y/V\rangle)$.
\end{proof}

\begin{definition}
\label{def:tghfr}
A \emph{finite henselian tangentially framed correspondence} 
between pairs $\langle X,U\rangle,\langle Y,V\rangle\in \Sm^\mathrm{pair}_S$
is defined by the data
\begin{itemize}
\item[(i)] 
A finite syntomic $f\colon W\to X$,
a closed subscheme $Z\subset W$ that intersects non-emptily each connected component of $W$,
and a morphism $g\colon W^h_Z\to Y$ 
that decomposes as $W^h_Z\to Y^\prime\to Y$ for $Y^\prime\in \Sm_S\cap \Aff$, and 
such that $Z\cong W^h_Z\times_Y (Y\setminus V)$,
and $U\times_X Z=\emptyset$;
\item[(ii)] 
a trivialisation $0\simeq \mathcal L_{f}$ of the cotangent complex $\mathcal L_{f}$ in $K$-theory $\infty$-gruppoid on $W$.
\end{itemize}
\end{definition}
Let $h^\tghfr(\langle Y/V\rangle )$ denote the presheaf with sections given by finite henselian tangentially framed correspondences.

\begin{definition}\label{def:tgfr}
A \emph{finite tangentially framed correspondence} 
between pairs $\langle X,U\rangle,\langle Y,V\rangle\in \Sm^\mathrm{pair}_S$ 
is defined by the data
\begin{itemize}
\item[(i)] 
A span $X\xleftarrow{f} W \xrightarrow{g} Y$ in $\SchPair$ with $f$ finite syntomic, 
and such that $g$
decomposes as $W\to Y^\prime\to Y$ for $Y^\prime\in \Sm_S\cap \Aff$, and
the subscheme $Z\defeq W\times_Y (Y\setminus V)$ is finite over $X$ and $U\times_X Z=\emptyset$.
\item[(ii)] 
a trivialisation $0\simeq \mathcal L_{f}$ of the cotangent complex $\mathcal L_{f}$ in $K$-theory $\infty$-gruppoid on $W$.
\end{itemize}
Define a presheaf $h^\tgfr(\langle Y/V\rangle )$
with sections given by finite tangentially framed correspondences.
\end{definition}

\begin{remark}
There is the equality of presheaves on $\Sm_S$
$h^\tgfr(Y)= h^\tghfr(\langle Y/V\rangle)$, where $V=\emptyset$.
\end{remark}

\begin{lemma}\label{lm:tghtg(Y/V)}
The morphism 
\[h^\tgfr(\langle Y\rangle )\to h^\tghfr(\langle Y\rangle )\]
is a Nisnevich equivalence in $\Pre(\Sm_S)$.
\end{lemma}
\begin{proof}
Let $U$ be an essentially smooth local henselian scheme over $S$. 
Consider a scheme $W$ finite over $U$
and a closed subscheme $Z$ that intersects non-emptily each connected component of $W$.
Since $U$ is local henselian, the scheme $W$ is a coproduct of local henselian schemes; 
so the henselisation $W^h_Z$ equals the minimal clopen subscheme in $W$ that contains $Z$; 
hence it is $W$ since $Z$ intersects each connected component.
Thus $W^h_Z\cong W$.

Then it follows that the morphism simplicial sets
$h^\tgfr(Y)(U)\to h^\tghfr(Y)(U)$
is an equivalence for any essentially smooth local henselian $U$. Hence 
$h^\tgfr(Y)\simeq_\nis h^\tghfr(Y)$
is a Nisnevich local equivalence.
\end{proof}

\begin{proposition}\label{prop:Y/VtgfrCMon}
Let $\langle Y,V\rangle\in \Sm^\mathrm{pair}_S$,
then the morphism 
\[\cofib ( h^\tgfr(V) \to h^\tgfr(Y) )\longrightarrow h^\tgfr(\langle Y/V\rangle)\] 
induces a Nisnevich equivalence in $\Pre(\Sm_S)$. 
\end{proposition}
\begin{proof}
The claim follows by \Cref{prop:CMonTgfr} and \Cref{lm:VveeY/V}.
\end{proof}

Recall that the $\infty$-category of $\infty$-commutative monoids $\CMon$ is 
canonically equivalent to the $\infty$-category 
$\Pre_\Sigma(\Gamma^{\mathrm{op}})$,
where $\Gamma=\mathcal{F}in_+^\mathrm{op}$, and $\mathcal{F}in_+$ is the category of finite pointed sets.
Let 
\[\Pre^\tgfr_\Sigma(\Sm_S), \Pre^{\CMon,\tgfr}_\Sigma(\Sm_S)\] be 
the $\infty$-categories of radditive (simplicial) preseahves on $\Corr^\tgfr(\Sm_S)$ 
and radditive presheaves of $\infty$-commutative monoids on $\Corr^\fr(\Sm_S)$
respectively.
Consider the adjunction 
\begin{equation}\label{eq:CMonfrfrv}
v^*\colon \Pre^\tgfr_\Sigma(\Sm_S)\simeq \Pre^{\CMon,\tgfr}_\Sigma(\Sm_S)\colon v_*,
\end{equation}
where $v_*$ denotes the forgetful functor, and $v^*$ is the left adjoint.
\begin{proposition}\label{prop:CMonTgfr}
Adjunction \eqref{eq:CMonfrfrv}
is an equivalence of $\infty$-categories.
\end{proposition}
\begin{proof}
The category $\Pre^{\CMon,\tgfr}_\Sigma(\Sm_S)$ is canonically equivalent to
the category of bi-radditive presheaves
$\Pre_\Sigma(\Corr^\fr(\Sm_S)\times \Gamma^\mathrm{op})$,
i.e. $\infty$-functors $F\colon \Corr^\fr(\Sm_S)\times \Gamma^\mathrm{op}\to \mathrm{sSet}$ such that 
$F((X\amalg Y,K))\simeq F((X,K))\times F((Y,K))$, and
$F((X,K\amalg L))\simeq F((X,K))\times F((X,L))$, for any $X,Y\in \Sm_S$, $K,L\in \Gamma^\mathrm{op}$%
.
In view of the mentioned equivalence 
\[\Pre^{\CMon,\tgfr}_\Sigma(\Sm_S)\simeq \Pre_\Sigma(\Corr^\fr(\Sm_S)\times \Gamma^\mathrm{op})\]
the functors $v_*$ and $v^*$ are induced by restriction along the functors
\[\begin{array}{lclllcr}
t&\colon &\Corr^\fr(\Sm_S)\times \Gamma^\mathrm{op}&\to& \Corr^\fr(\Sm_S)&\colon& e\\
&&(X,K)&\longmapsto&X\times K&\\
&&(X,[1])&\longmapsfrom&X&
\end{array}
\]
So we can write \[v_*= e_*, v^*=t_*.\]
Then by the definition
\[t e (X) = X, (e\circ t)((X,K))\cong (X\times K,[1]).\]
So immediately we have
\[e_*t_*(F)\simeq F\]
for any $F\in \Pre(\Corr^\fr(\Sm_S))$ and 
\[t_*e_*(F)\simeq F,\]
for any 
$F\in \Pre_\Sigma(\Corr^\fr(\Sm_S)\times \Gamma^\mathrm{op})$,
because of equivalences 
\[F((X\times K,[1]))=F((\coprod_{\alpha\in K}X,[1])\simeq \prod_{_{\alpha\in K}}F((X,[1]))\simeq F((X,\coprod_{\alpha\in K}[1]))=F((X,[1]\times K)).\]
\end{proof}

\begin{lemma}\label{lm:VveeY/V}
The canonical morphism gives 
the Nisnevich equivalence 
\begin{equation}\label{eq:cofib(tgfrV->tgfrY)=tgfrY/V}\cofib( v^*h^\tgfr(V)\to v^*h^\tgfr(Y) )\simeq v^*h^\tgfr(\langle Y/V\rangle)\end{equation}
in $\Pre^{\CMon,\tgfr}(\Sm_S)$.
\end{lemma}
\begin{proof}
By \Cref{lm:V+Y-V:v^*(tgfr)(U)} there is an equivalence
$v^*h^\tgfr(Y)(U) \simeq v^*h^\tgfr(V)(U) \oplus v^*h^\tgfr(\langle Y/V\rangle)(U)$,
for $U=X^h_x$ for any $X\in \Sm_{S}$, and $x\in X$.
So morphism \eqref{eq:cofib(tgfrV->tgfrY)=tgfrY/V}
induces the equivalence on the stalks at essentially smooth henselian local schemes.
\end{proof}

\begin{lemma}\label{lm:V+Y-V:v^*(tgfr)(U)}
For any $\langle Y,V\rangle\in \Sm^\mathrm{pair}_S$, and a local henselian $U=X^h_x$, $X\in \Sm_{S}$, $x\in X$,
the sequence of $\infty$-commutative monoids
\[v^*h^\tgfr(Y)(U)\to v^*h^\tgfr(V)(U)\to v^*h^\tgfr(\langle Y/V\rangle)(U)\]
splits, and there is an equivalence
\begin{equation}\label{eq:V(U)oplusY/V(U)}
p\colon v^*h^\tgfr(Y)(U) \simeq v^*h^\tgfr(V)(U) \oplus v^*h^\tgfr(\langle Y/V\rangle)(U),
\end{equation}
where $\oplus$ denotes (co)product in the category $\CMon$ of $\infty$-commutative monoids.
\end{lemma}
\begin{proof}
By the definition there are canonical morphisms
\[
i_V\colon 
h^\tgfr(V)\to h^\tgfr(Y), 
i_{Y\setminus V}\colon 
h^\tgfr(\langle Y/V\rangle)\to h^\tgfr(Y)
\]
Define morphisms
\[
(-)^V\colon h^\tgfr(Y)\to h^\tgfr(V), 
(-)^{Y\setminus V}\colon h^\tgfr(Y)\to h^\tgfr(\langle Y/V\rangle)\]
as follows.
Any scheme $W$ finite over the local henselian $U$ 
is a coproduct of local henselian schemes, and
for any closed subscheme $Z\subset W$
there is a canonical splitting 
\[W=W^Z\amalg W^{-Z}\]
such that $W^Z$ is the minimal clopen subscheme such that $W^Z\supset Z$.
So for any finite syntomic span 
$(W,f,g)$ we have
\begin{equation}\label{eq:WsupV+(Y-V)}
W = W^V \amalg W^{Y\setminus V}, \K(W)\simeq \K(W^V)\times \K(W^{Y\setminus V}),
\end{equation}
where 
$W^V=W^Z$ for $Z=W\times_Y (Y\setminus V)$,
and 
the equivalence at the right side above 
is provided by the morphisms 
$(-)\big|_{W^V}\colon \K(W)\simeq \K(W^V)$ and $(-)\big|_{W^{Y\setminus V}}\colon \K(W)\simeq \K(W^{Y\setminus V})$.
Then there is a unique lift $g^\prime\colon W^V\to V$ of $g\big|_{W^V}\colon W^V\to Y$.
So there are morphisms 
\[\begin{array}{llll}
(-)^V\colon& h^\tgfr(Y)(U)&\to& h^\tgfr(V)(U);\\ 
&(W,f,\phi,g) &\mapsto& (W^V,f\big|_{W^V},\phi\big|_{W^V},g^\prime) \\
(-)^{Y\setminus V}&\to& h^\tgfr(Y)(U)\to h^\tgfr(\langle Y/V\rangle);\\ 
&(W,f,\phi,g)&\mapsto& (W^{Y\setminus V},f\big|_{W^{Y\setminus V}}, \phi\big|_{W^{Y\setminus V}},g\big|_{W^{Y\setminus V}})
\end{array}\]

We claim that the morphisms
\begin{equation}\label{eq:q}\begin{array}{llll}
q\colon &
h^\tgfr(Y)(U)&\to& V^\tgfr(U)\oplus h^\tgfr(\langle Y/V\rangle)(U) ;\\&  
c &\mapsto& (c^V, c^{Y\setminus V})\\
q^{-1}\colon &
h^\tgfr(V)(U)\oplus h^\tgfr(\langle Y/V\rangle)(U)&\to&  h^\tgfr(Y)(U) ;\\&  
(c_1,c_2) &\mapsto& c_1 + c_2
\end{array}\end{equation}
induce the pair of the inverse equivalences $v^*(q)$ and $v^*(q^{-1})$.
The claim follows because of the equalities
\[
q\circ i_{Y\setminus V}\simeq id_{h^\tgfr(\langle Y/V\rangle)(U)}, 
q \circ i_V\simeq id_{h^\tgfr(V)(U)} , 
\]
provided immediately by the definitions,
and the equivalence
\begin{equation}
\label{eq:tgfrV+(Y-V)}
id_{v^*h^\tgfr(Y)(U)}\simeq (-)^V+(-)^{Y\setminus V}\colon 
c\simeq j(c^V)+i(c^{Y\setminus V})
\end{equation} 
provided by \eqref{eq:WsupV+(Y-V)}. 
To explain equivalence \eqref{eq:tgfrV+(Y-V)} in detail, let us 
consider firstly the presheaf $v^*(h^\syncor(\langle Y/V\rangle))$. 
Since $h^\syncor(Y)$ is the presheaves of the classical pointed sets, 
$v^*h^\syncor(Y)$ is the presheaf of classical commutative monoids. 
So \eqref{eq:WsupV+(Y-V)} implies the classical equality 
\begin{equation}\label{eq:syncorV+(Y-V)}
id_{v^*h^\syncor(Y)(U)}= (-)^V+(-)^{Y\setminus V}.\end{equation} 
Then because of the equivalence $K(W)\simeq K(W^V)\times K(W^{Y\setminus V})$ 
equality \eqref{eq:syncorV+(Y-V)} provides equivalence \eqref{eq:tgfrV+(Y-V)}.
\end{proof}

\subsection{$\A^1$-homotopy equivalence}

\begin{definition}\label{def:lownr}
A level $n$ \emph{lower normally framed correspondence} between pairs 
$\langle X,U\rangle,\langle Y,V\rangle\in \Sm^\mathrm{pair}_S$ is given by the following data:\begin{itemize}
\item[(i)] closed subschemes $Z\subset W\subset \A^n\times X$
such that $U\times_X Z=\emptyset$, 
$Z$ is finite over $X$,
and $Z$ intersects each connected component of $W$.
\item[(ii)] a set of regular functions $\varphi_i\in \calO( (\A^n\times X)^{1\mathrm{th}}_Z )$, $i=1,\dots n$,
such that $Z(\varphi_1,\dots \varphi_n)=W$,
\item[(iii)] a morphism $g\colon W\to Y$ 
that decomposes as $W\to Y^\prime\to Y$ for $Y^\prime\in \Sm_S\cap \Aff$, and
such that $Z\cong W\times_Y (Y\setminus V)$,
and $U\times_X Z=\emptyset$.
\end{itemize}

Denote by $h^\lownrfr(\langle Y,V\rangle)(\langle X,U\rangle)$ the set of lower normally framed correspondences from 
$\langle X,U\rangle$ to $\langle Y,V\rangle$

A level $n$ \emph{upper normally framed correspondence} between pairs 
$\langle X,U\rangle,\langle Y,V\rangle\in \Sm^\mathrm{pair}_S$ 
is given by the following data:\begin{itemize}
\item[(i)] closed subschemes $Z\subset W\subset \A^n\times X$ finite over $X$, 
such that $U\times_X Z=\emptyset$, and $Z$ intersects each connected component of $W$.
\item[(ii)] a set of regular functions $\phi_i\in \calO( (\A^n\times X)^{1\mathrm{th}}_W )$, $i=1,\dots n$,
such that $Z(\phi_1,\dots \phi_n)=W$,
\item[(iii)] and a morphism $g\colon W^h_Z\to Y$
that decomposes as $W^h_Z\to Y^\prime\to Y$ for $Y^\prime\in \Sm_S\cap \Aff$, and
such that $Z\cong W^h_Z\times_Y (Y\setminus V)$,
and $U\times_X Z=\emptyset$.
\end{itemize}
The scheme $Z$ in the above definitions is called \emph{the support} of the correspondence, and $W$ (if it is) is called \emph{the big-support}.

\end{definition}

\begin{proposition}\label{prop:freq:Nr=Tg}
For any $\langle Y,V\rangle\in \Sm^\mathrm{pair}_S$, 
the morphism $h^\uppnrfr(\langle Y/V\rangle) \to h^\tghfr(\langle Y/V\rangle)$ 
is a motivic equivalence in $\Pre(\Sm_S)$.
\end{proposition}
\begin{proof}
The claim follows from \cite[Lemmas 2.3.12, 2.3.23]{five-authors}
similarly to \cite[Proposition 2.3.25]{five-authors}.
\end{proof}

A morphism $f$ satisfies \emph{the lifting property with respect to the closed immersions on affines} if
for any closed immersion $Y\to X$ in $\SmAff_S$, the morphism $\pi_0 F(X)\to \pi_0( G(X) \times_{G(Y)} F(Y) )$ is surjective.

\begin{lemma}\label{lm:AAffmot-criterion}
Let $S$ be a scheme and $f\colon F\to G$ be a morphism in $\Pre(\Sm_S)$.
Assume that $f$ satisfies the lifting property with respect to the closed immersions on affines, and closed glueing.
Then the restriction of $f$ in $\Pre(\AffSm_S)$ is $\A^1$-equivalence, and $f$ is a motivic equivalence in $\Pre(\Sm_S)$.
\end{lemma}
\begin{proof} See the arguments of \cite[Corollary 2.2.20]{five-authors}.
\end{proof}

\begin{proposition}\label{prop:freq:Nis=Nr=Id}
For any $\langle Y,V\rangle\in \Sm^\mathrm{pair}_S$,
the restrictions of the morphisms
\begin{equation}\label{eq:Frflownfruppernrfr}
\Lrep_{\A^1}h^\fraff(\langle Y/V\rangle)\to 
\Lrep_{\A^1}h^\lownrfr(\langle Y/V\rangle)\leftarrow 
\Lrep_{\A^1}h^\uppnrfr(\langle Y/V\rangle)
\end{equation}
to the category $\Sm_S\cap \Aff$ are trivial fibrations.
Consequently, morphisms of presheaves \eqref{eq:Frflownfruppernrfr} on $\Sm_S$
are motivic equivalences.
\end{proposition}
\begin{proof}

We are going to show that the morphism
$\Fr(\langle Y/V\rangle)\to h^\lownrfr(\langle Y/V\rangle)$
satisfies assumption of \Cref{lm:AAffmot-criterion}.
Let 
\[(Z,W,\phi,g)\in h^\lownrfr(Y/V)(X), (Z_0, \tilde \phi_0,\tilde g_0)\in \Fr(Y/V)(X_0)=\Fr(X_0,Y/V)\]
define an element in 
$h^\lownrfr(Y/V)(X)\times_{h^\lownrfr(Y/V)(X_0)} \Fr(Y/V)(X_0)$.
We lift it to an element in $\Fr(Y/V)(X)$.
Since $Y\times \A^n$ is smooth over $S$, 
it follows by \cite[Lemma B.10]{DKO:SHISpecZ} 
that
there is a regular map 
$(\tilde \phi,\tilde g)\colon (\A^n\times X)^h_Z \to \A^n\times Y$ 
such that 
\[\begin{array}{ll}
\tilde \phi\big|_{(\A^n\times X_0)^h_{Z_0} }=\tilde \phi_0\big|_{(\A^n\times X_0)^h_{Z_0} }, &
\tilde \phi\big|_{(\A^n\times X)^\fth_{Z}}=\phi\big|_{(\A^n\times X)^\fth_Z}, \\
\tilde g\big|_{(\A^n\times X_0)^h_{Z_0} }=\tilde g_0\big|_{(\A^n\times X_0)^h_{Z_0} }, &
\tilde g\big|_{W}=g.
\end{array}\]
Define now the required lift as the correspondence
\[(Z, \tilde \phi, \tilde g)\in \Fr(Y/V)(X).\]

We are going to show that the morphism
$h^\uppnrfr(\langle Y/V\rangle)\to h^\lownrfr(\langle Y/V\rangle)$
satisfies assumption of \Cref{lm:AAffmot-criterion}.
Let a pair of elements
\[(Z,W,\phi,g)\in \Corr^\lownrfr(X,Y/V), (Z_0, \tilde \phi_0,\tilde g_0)\in \Corr^\uppnrfr(X_0,Y/V)\]
defines the element in 
$\Corr^\lownrfr(X,Y/V)\times_{\Corr^\lownrfr(X_0,Y/V)} \Corr^\uppnrfr(X_0,Y/V)$.
By Serre's theorem on ample bundles \cite[Chapter~III, Theorem~8.8]{Hartshorne}
for large enough $d$
there are sections $s_i\in \Gamma(\PP^{n}\times X \mathcal O(d))$, $i=1,\dots n$,
\[
s_i\big|_{\PP^{n-1}\times X}=t_i^d, 
s_i\big|_{R}=t_\infty^d v_i
\]
where 
\[R= (\A^n\times X_0)^\fth_{W_0}\amalg_{(\A^n\times X_0)^\fth_{Z_0}} (\A^n\times X)^\fth_Z\]
and 
\[v_i=\tilde \phi^0_i\amalg_{\phi^0_i} \phi_i\in \Gamma(R,\mathcal O_R).\]
\end{proof}

\begin{corollary}\label{cor:frafftghfr}
For each $\langle Y,V\rangle\in \Sm^\mathrm{pair}_S$, there is a canonical motivic equivalence 
$h^\fraff(\langle Y/V\rangle) \to h^\tghfr(\langle Y/V\rangle)$
in $\Pre(\Sm_S)$.
\end{corollary}

\section{Counterexample to the Nisnevich cancellation over a base.}\label{section:ContrexampleCancellationDMeff}

Let $\Pre^{\tr,\Ab}(B)$ denote the category of additive functors $\Cor(B)\to \Ab$, and
$\Shv^{\tr,\Ab}(B)$ denote the subcategory of Nisnevich sheaves. 
Let $w_{\A^1}$ denote the class of $\A^1$-homotopy equivalences in $D(\Pre^{\tr,\Ab}(B))$. 
Then \[\DM(B)=D(\Shv^{\tr,\Ab}(B))[w_{\A^1}^{-1}][\Gm^{\wedge -1}],\]
where $\Gm^{\wedge 1}=\Cone(\{1\}\to \Gm)\in D(\Shv^{\tr,\Ab}(B))$,
is the category of Voevodsky's motives.
In the section, we show that 
the functor
\begin{equation}\label{eq:DAShvtoDM}D(\Shv^{\tr,\Ab}(B))[w_{\A^1}^{-1}]\to \DM(B)\end{equation}
is not fully faithful
for any separable noetherian scheme $B$ of positive finite Krull dimension.

\begin{lemma}\label{lm:wedgeGinjweldgeGl}
Let $\calC$ be a monoidal stable $\infty$-category and $\mathbb G$ be an object,
and $\calC[\mathbb G^{\wedge -1}]$ denote the universal $\infty$-category such that the functor $-\wedge \mathbb G^{\wedge -1}$ is an autoequivalence.
The following conditions are equivalent:
\begin{itemize}
\item[(1)] 
the functor \begin{equation}\label{eq:wedgeGcalC}\calC\to \calC\colon X\mapsto X\wedge \mathbb G\end{equation} is a fully faithful embedding;
\item[(2)] 
the functor \[\calC\to\calC[\mathbb G^{\wedge -1}]\] is a fully faithful embedding;
\item[(3)] there is an equivalence \[\calC\to \calC(1)\colon X\mapsto X\wedge \mathbb G,\]
where $\calC(1)$ is the subcategory $\calC$ spanned by the objects of the form $X\wedge \mathbb G$, $X\in \calC$.
\end{itemize}
\end{lemma}
\begin{proof}
Points (1) and (3) are equivalent by the definition.
Since for any $X,Y\in \calC$,
\[\Map_{\calC[\mathbb G^{\wedge -1}]}(X,Y)\simeq\varinjlim_{l}\Map_{\calC}(X\wedge \mathbb G^{\wedge l},Y\wedge \mathbb G^{\wedge l}),\]
point (1) implies point (2).
Conversely, 
if for all $X,Y\in \mathcal C$,
\begin{equation}\label{eq:piicalCXYcalCG-1XY}\pi_i\Map_{\calC}(X,Y) \simeq \pi_i\varinjlim_{l}\Map_{\calC}(X\wedge \mathbb G^{\wedge l},Y\wedge \mathbb G^{\wedge l}), i\in \mathbb Z,\end{equation}
then
\begin{equation}\label{eq:piicalCXYcalCXGYG}\pi_i\Map_{\calC}(X,Y) \hookrightarrow \Map_{\calC}(X\wedge \mathbb G^{\wedge 1},Y\wedge \mathbb G^{\wedge 1}).\end{equation}
Then all the transition maps in the diagram of the limit $\pi_i\varinjlim_{l}\Map_{\calC}(X\wedge \mathbb G^{\wedge l},Y\wedge \mathbb G^{\wedge l})$ are injective, and consequently, surjective, because of the isomorphism \eqref{eq:piicalCXYcalCG-1XY}.
Thus \eqref{eq:piicalCXYcalCXGYG} is an isomorphism.
\end{proof}

\begin{proposition}\label{prop:nisCancOmegaGmnis}
Let $B$ a noetherian separated scheme of finite Krull dimension.
Assume that the functor \eqref{eq:DAShvtoDM}
is fully faithful.
Then the functor
\begin{equation}\label{OmegaGmDPretrA1}
\Omega_{\Gm}\colon D(\Pre^{\tr,\Ab}(B))[w_{\A^1}^{-1}]\to D(\Pre^{\tr,\Ab}(B))[w_{\A^1}^{-1}]
\end{equation}
preserves Nisnevich equivalences.
\end{proposition}
\begin{proof}
Consider the category $D(\Pre^\tf(B))$
and the functor
\begin{equation}\label{eq:SigmaGm1DPretfB}
\Sigma^1_{\Gm}\colon D(\Pre^\tf(B))[w_{\A^1}^{-1}]\to D(\Pre^\tf(B))[w_{\A^1}^{-1}]
\end{equation}
By \cite{Voe-cancel} the functor 
\eqref{eq:SigmaGm1DPretfB} is an embedding.
So the pair
\begin{equation}\label{eq:Sigma1Omega1GmPretrB}\Sigma^1_{\Gm}\colon D(\Pre^\tf(B))[w_{\A^1}^{-1}]\to D(\Pre^\tf(B))[w_{\A^1}^{-1}](1)\colon\Omega^1_{\Gm}\end{equation}
is an equivalence,
where
the right side is by the definition the image of \eqref{eq:SigmaGm1DPretfB}.
If $D(\Shv^{\tr,\Ab}(B))[w_{\A^1}^{-1}]\to \DM(B)$ is fully faithful, then by \Cref{lm:wedgeGinjweldgeGl}
the functor
\begin{equation}\label{SigmaGm1ShvtrB}\Sigma_{\Gm}^1\colon D(\Shv^{\tr,\Ab}(B))[w_{\A^1}^{-1}]\to D(\Shv^{\tr,\Ab}(B))[w_{\A^1}^{-1}]\end{equation} is fully faithful.
Then
\begin{equation}\label{eq:Sigma1Omega1GmShvtrB}\Sigma^1_{\Gm}\colon D(\Shv^\tf(B))[w_{\A^1}^{-1}]\to D(\Shv^\tf(B))[w_{\A^1}^{-1}](1)\colon\Omega^1_{\Gm}\end{equation}
is an equivalence,
where 
\[D(\Shv^\tf(B))[w_{\A^1}^{-1}](1)=D(\Pre^\tf(B))[w_{\A^1}^{-1}](1)[w_{\nis}^{-1}]\]
is equivalent to the image of the functor \eqref{SigmaGm1ShvtrB}.
Thus since both pairs \eqref{eq:Sigma1Omega1GmPretrB}, \eqref{eq:Sigma1Omega1GmShvtrB} are equivalences, and the left adjoint functor $\Sigma^1_{\Gm}$ in \eqref{eq:Sigma1Omega1GmPretrB} is Nisnevich exact, it follows that $\Omega^1_{\Gm}$ in \eqref{eq:Sigma1Omega1GmPretrB} is Nisnevich exact.
\end{proof}

\begin{lemma}\label{ex:A1CornoncontractabletfSquareGmB}
Let $B$ be a local scheme, $z\in B$ be the closed point.
Then $H^1_\nis(B\times \Gm,\calF)\neq 0$ for some $\A^1$-homotopy invariant presheaf with transfers $\calF$.
\end{lemma}
\begin{proof}
The argument is like in \cite[\S 13]{DKO:SHISpecZ}zero regular function, $f(z)=0$.
Let $g=(1-ft)\in \calO_{B\times\Gm}(B\times\Gm)=\calO_B(B)[t,t^{-1}]$.
Then the square Nisnevich (Zariski)
\[\xymatrix{
B\times\Gm-Z(fg)& B\times\Gm-Z(f)\ar[l]\\
B\times\Gm-Z(g)\ar[u]& B\times\Gm\ar[u]\ar[l]\\
}\]
defines non-trivial element in $H^1_\nis(B\times\Gm, \calF)$, where 
\[\calF=\Coker( \Cor_B(-\times\A^1,B\times\Gm-Z(fg))\to \Cor_B(-,B\times\Gm-Z(fg)) ).\]
\end{proof}

\begin{proposition}
Let $B$ be a noetherian separated scheme of positive Krull dimension.
Then 
the functor \eqref{OmegaGmDPretrA1}
is not Nisnevich exact, 
and the functor \eqref{eq:DAShvtoDM}
is not fully faithful.
\end{proposition}
 \begin{proof}
Without loss of generality we can assume that $B$ is local.
Any $\A^1$-invariant presheaf with transfers defines canonically an element in $D(\Pre^{\tr,\Ab}(B))[w_{\A^1}^{-1}]$,
and 
\[\begin{array}{rll}
\calF(B\times\Gm)&=&\mathrm{Hom}_{D(\Pre^{\tr,\Ab}(B))[w_{\A^1}^{-1}]}(B\times\Gm, \calF),\\
0&=&\mathrm{Hom}_{D(\Pre^{\tr,\Ab}(B))[w_{\A^1}^{-1}]}(B\times\Gm, \calF[i]), i>0,\\
H^i_\nis(B\times\Gm, \calF)&=&\mathrm{Hom}_{D(\Pre^{\tr,\Ab}(B))[w_{\A^1}^{-1}]}(B\times\Gm, \Lrep_\nis\calF[i]),
\end{array}\]
where $\Lrep_\nis$ denotes the composition
\[D(\Pre^{\tr,\Ab}(B))[w_{\A^1}^{-1}]\to D(\Shv^{\tr,\Ab}(B))[w_{\A^1}^{-1}]\to D(\Pre^{\tr,\Ab}(B))[w_{\A^1}^{-1}],\]
where the first functor is the localisation, and second functor is the canonical embedding.
If \eqref{OmegaGmDPretrA1} would be Nisnevich exact, then 
\[\mathrm{Hom}_{D(\Pre^{\tr,\Ab}(B))[w_{\A^1}^{-1}]}(B\times\Gm, \calF)\cong \mathrm{Hom}_{D(\Pre^{\tr,\Ab}(B))[w_{\A^1}^{-1}]}(B\times\Gm, \Lrep_\nis\calF);\]
hence
\[H^1_\nis(B\times\Gm, \calF)\cong \mathrm{Hom}_{D(\Pre^{\tr,\Ab}(B))[w_{\A^1}^{-1}]}(B\times\Gm, \calF[1])\cong 0.\]
So by \Cref{ex:A1CornoncontractabletfSquareGmB} the first claim follows.
By \Cref{prop:nisCancOmegaGmnis} the first claim implies the second.
\end{proof}

\section{Computing of the stable motivic localisation.} 
\label{sect:stablemotiviclocalisation}

Throughout the section 
$k$ denotes a field that satisfy the 
strict homotopy invariance 
assumption (1) 
from \Cref{sect:AssumptionsBase},
and $S$ denotes a noetherian separated scheme of finite Krull dimension
that satisfy assumptions (1) and (2) from \Cref{sect:AssumptionsBase}.

\subsection{Reformulation lemma}\label{sect:reformlemma}

By \Cref{lm:LreptLprepcpreserveexact}, that follows in the text,
the strict homotopy invariance theorem over $k$ in sense of point (1) from \Cref{sect:AssumptionsBase} 
can be formulated 
in two equivalent ways: (1) to say that $\Lrep_\nis$ on $\Pre^{S^1,\fr}_\Sigma(k)$ preserves $\A^1$-invariant spectra, and (2) to say that $\Lrep_{\A^1}$ on $\Pre^{S^1,\fr}_\Sigma(k)$ preserves Nisnevich local equivalences.
The way (1) is the standard formulation. The reformulation (2) was used already by the author for the argument for the case of $\mathrm{GW}$-transfers in \cite{DruDMGWeff}.
The reformulation allows to make the Nisnevich localisation procedure being external one with respect to $\A^1$-invariantisation and $\Gm$-stabilisation;
let us copy the diagram for $\mathbf{DM}^\mathrm{GW}(k)$ from \cite{DruDMGWeff} applying it to $\SH^\fr(k)$
\begin{equation}\label{eq:Pre(k)Shvnis(k)A1Gm}
\xymatrix{
\Spt^{S^1,\fr}(k)\ar[r]\ar[d]& \Spt^{S^1,\fr}_{\A^1}(k) \ar[r]\ar[d]& \Spt^{S^1,\fr}_{\A^1}(k)[\Gm^{\wedge -1}]\ar[d] \\
\Spt^{S^1,\fr}_\nis(k)\ar[r]& \SH^{S^1,\fr}_{\A^1,\nis}(k) \ar[r]& \rmSH^{S^1,\Gm,\fr}_{\A^1,\nis}(k)
.}\end{equation}
The first row is a 'naive' framed motivic category that deal with 'naive' $\A^1$-homotopies and $\Gm$-loops.
%
This is useful for the base scheme case here, because it allows to make some part of arguments on $\A^1$-localisation and $\Gm$-stabilisation being unique, and combine them. 
Firstly, we write the diagram \eqref{eq:Pre(k)Shvnis(k)A1Gm} for base scheme case with $\Spt^{S^1,\fr}_\tf(S)$ at the first row instead of $\Spt^{S^1,\fr}(k)$. Then we continue the way of thoughts suggested bu \cite{DruDMGWeff} further and short the diagram using that 
the localisation in presentable categories is a particular case of the stabilisation, the $\A^1$-localisation and $\Gm$-stabilisation can be united into one stabilisation procedure with respect to the functor 
\begin{equation}\label{eq:Delta_SwedgeGm}\Omega_{\Gm}\Lrep_{\A^1}\colon F(-)\mapsto F((-)_+\wedge(\Delta_S^\bullet)_+\wedge \Gm).\end{equation}
So we get the diagram 
\begin{equation}\label{eq:Pre(S)Shvnis(S)A1Gm}\xymatrix{
\Spt^{S^1,\fr}_\tf(S)\ar[r]\ar[d]& \Spt^{S^1,\fr}_{\A^1,\tf}(S)[\Gm^{\wedge -1}]\ar[d] \\
\Spt^{S^1,\fr}_\nis(S)\ar[r]& \rmSH^{S^1,\Gm,\fr}_{\A^1,\nis}(S)
,}\end{equation}
with rows given by the stabilisation with respect to the endofunctor \eqref{eq:Delta_SwedgeGm}.

The used lemmas in the abstract form sounds as follows.
\begin{lemma}\label{lm:LreptLprepcpreserveexact}
Let $\mathcal C$ be an $\infty$-category, and
$\mathcal C_t$, and $\mathcal C_c$, $\mathcal C_{tc}=\mathcal C_t\cap \mathcal C_c$
be reflective subcategories. 
Let
$\Lrep_t$, $\Lrep_c$, $\Lrep_{tc}$ denote the localisation endofunctors on $\mathcal C$ that lands in $\mathcal C_t$, $\mathcal C_c$, $\mathcal C_{tc}$ respectively.

The following conditions are equivalent:
\begin{itemize}
\item[(1)] the canonical morphism $\Lrep_{tc}\leftarrow \Lrep_t \Lrep_c$ is an equivalence;
\item[(2)] $\Lrep_t$ preserves $\mathcal C_c$;
\item[(3)] $\Lrep_c$ preserves $t$-equivalences.
\end{itemize}
\end{lemma}
\begin{proof}
Let us show that (2) implies (1).
Since $\Lrep_t$ preserves $\mathcal C_c$, and $\Lrep_c$ lands in $\mathcal C_c$, 
it follows that $\Lrep_t \Lrep_c$ lands in $\mathcal C_c$.
Then since $\Lrep_c$ is equivalent to the identity after the restriction to $\mathcal C_c$, 
the canonical morphism $\Lrep_c (\Lrep_t \Lrep_c)\leftarrow \Lrep_t \Lrep_c$ is an isomorphism.
Hence $\Lrep_c (\Lrep_t \Lrep_c)$ lands in $\mathcal C_t$, and then 
the canonical morphisms $\Lrep_t(\Lrep_c (\Lrep_t \Lrep_c)) \leftarrow \Lrep_c (\Lrep_t \Lrep_c)\leftarrow \Lrep_t\Lrep_c$ are isomorphisms.
The composition of the latter composite morphism with the isomorphism 
$(\Lrep_t\Lrep_c) (\Lrep_t \Lrep_c))\simeq \Lrep_t(\Lrep_c (\Lrep_t \Lrep_c))$
is equivalent to the canonical morphism 
$(\Lrep_t\Lrep_c) (\Lrep_t \Lrep_c))\leftarrow \Id(\Lrep_t \Lrep_c)$.
Hence the latter morphism is an isomorphism.
Thus since $\Lrep_{tc} \simeq  \varinjlim_l \Id(\Lrep_t\Lrep_c)^l$
$\Lrep_{tc}\simeq \Lrep_t\Lrep_c$.

Let us show that (1) implies (2).
Since $\Lrep_{tc}$ lands in $\mathcal C_c$, it preserves $\mathcal C_c$. Then by (1) $\Lrep_t \Lrep_c$ preserves $\mathcal C_c$.

Let us show that (1) implies (3).
Let $F\to G$ be a t-equivalence. Then $\Lrep_{t}F\simeq \Lrep_{t}G$, and $\Lrep_{tc}F\simeq \Lrep_{tc}G$.
Hence by (1) $\Lrep_t (\Lrep_c F)\simeq \Lrep_s (\Lrep_c G)$. Thus $\Lrep F\to \Lrep G$ is a t-equivalence.

Let us show that (3) implies (1).
The canonical morphism $\Lrep_t \Lrep_c \leftarrow \Lrep_c$ is a t-equivalence.
Since $\Lrep$ preserves t-equivalences, $\Lrep_c (\Lrep_t \Lrep_c)\leftarrow \Lrep_c\Lrep_c\colon v_1$ is a t-equivalence.
Since the canonical morphism $\Lrep_c\Lrep_c\leftarrow \Lrep_c\colon v_2$ is an isomorphism,
and $v_1v_2$ is equivalent to the canonical morphism 
$\Lrep_c (\Lrep_t \Lrep_c)\leftarrow \Lrep_c \Id$, the latter morphism is a t-equivalence.
The composite
\[(\Lrep_t\Lrep_c) (\Lrep_t \Lrep_c)\simeq \Lrep_t(\Lrep_c (\Lrep_t \Lrep_c))\leftarrow \Lrep_t (\Lrep_c \Id)\simeq (\Lrep_t \Lrep_c) \Id\] is 
equivalent to the canonical morphism
$(\Lrep_t\Lrep_c) (\Lrep_t \Lrep_c)\leftarrow (\Lrep_t \Lrep_c) \Id$.
Thus it follows that $\Lrep_{tc} \simeq  \varinjlim_l (\Lrep_t\Lrep_c)^l\Id\simeq 
\Lrep_t\Lrep_c$.
\end{proof}

\begin{lemma}\label{lm:LocsubsetStab}
Given an $\infty$-category $C$ and a reflective subcategory, with the localisation functor $L_t\colon C\to C_t$, that is left adjoint to the canonical embedding, then $C_t\simeq C[\Lrep_t^{-1}]$, where $\Lrep_t$ is the composite $E_tL_t\colon C\to C_t\to C$ of $L_t$ and the canonical embedding $E_t\colon C_t\to C$.
\end{lemma}
\begin{proof}
The category $C[\Lrep_t^{-1}]$ is the universal category equipped with the functor $s\colon C\to C[\Lrep_t^{-1}]$ and the endofunctors $\Lrep_t$ and $\Lrep_t^{-1}$ such that $\Lrep_t^{-1}\Lrep_t\simeq \mathrm{Id}_{C[\Lrep_t^{-1}}$ and $s\Lrep_t\simeq \Lrep_ts$.
Since the category $C_t$, the functor $L_t\colon C\to C_t$ and $\Lrep_t=\Lrep_t^{-1}=\mathrm{Id}_{C_t}$ have the above list of properties, there is the canonical functor $C[\Lrep_t^{-1}\to C_t$.
The functor $L_t$ is idempotent. Hence $L_t\simeq L_t^{-1}$ on $C[L_t^{-1}]$ are autoequivalences. Hence $s$ takes $t$-equivalences to isomorphisms, and there is the canonical functor $C_t\simeq C[t^{-1}]\to C[\Lrep_t^{-1}]$.
Thus we get the commutative diagram
\[\xymatrix{
C\ar[r]\ar[dr]& C_t\ar@<1ex>[d]\\
& C[\Lrep_t^{-1}]\ar@<1ex>[u]
,}\]
and the composite endomorphisms are equivalences by the universal properties.
\end{proof}

\begin{lemma}
Let $C$ be an $\infty$-category,
$C_t$ be a reflective subcategory, $L_t\colon C\to C_t$ be the localisation functor left adjoint to the canonical embedding $E_t\colon C_t\to C$. 
Let $\Sigma\dashv \Omega$ be an adjunction 
$\Sigma C \leftrightarrows C \Omega$.
Denote $\Omega_t=L_t \Omega E_t\colon C_t\to C_t$.

Suppose that $\Sigma$ and $\Omega$ preserve the class of $t$-equivalences,
then $L_t$ induces the functor \[L_t^\Omega=L_t[\Omega^{-1}]\colon C[\Omega^{-1}]\to C_t[\Omega_t^{-1}],\]
such that \[\Omega_{t}^\infty\Sigma^\infty_{t}L_t\simeq L_t\Omega_{}^\infty\Sigma^\infty_{}.\]
\end{lemma}
\begin{proof}
Since the functors $\Sigma$ and $\Omega$ are exact with repsepct to $t$-equivalences, it follows that 
the adjunction $\Sigma\dashv \Omega$ induces the adjunction $\Sigma_t\dashv \Omega_t$ of endofunctors on $C_t$.
The claim follows.
\end{proof}

\begin{lemma}\label{lm:ReformulationOmega}
Let $C$ be an $\infty$-category.

Let $C_t$ be a reflective subcategory, $L_t\colon C\to C_t$ be the localisation functor left adjoint to the canonical embedding $E_t\colon C_t\to C$, and define the class of $t$-equivalences in $C$ as the class of morphisms that maps to isomorphisms via $L_t$. 

Let $\Sigma\dashv \Omega$ be the adjunction of endofunctor  $C\to C$, that we denote by $G$. 
Denote $\Sigma_t=L_t \Sigma E_t\colon C_t\to C_t$.

Suppose that $\Sigma$ preserves $t$-equivalences.
The following conditions are equivalent

\begin{itemize}
\item[(1)]$L_t \Omega^\infty \simeq \Omega^\infty_t L^\Omega_t$, 
\item[(2)] $\Lrep_t \Lrep_{\Omega} \simeq \Lrep_{\Omega,t}\colon \PSpt_\Omega(C)\to \PSpt_\Omega(C)$ 
\item[(3)]$\Omega$ preserves $t$-equivalences 
\end{itemize}

\end{lemma}
\begin{proof}

Point (3) equivalently means the equivalences
of functors
\begin{equation}\label{eq:LtnOmega(t)}\Omega L_t \Omega \cong \Omega_t L_t.\end{equation}
The latter equivalences implies the equivalence
$L_t \Omega^\infty \simeq \Omega^\infty_t L^\Omega_t$.
So \Cref{lm:ReformulationOmega}(3) implies \Cref{lm:ReformulationOmega}(1) 
This part of the \Cref{lm:ReformulationOmega} is only what we essentially use in the further applications in the section.

We continue the argument to complete the proof of \Cref{lm:ReformulationOmega} in the full form.
\eqref{eq:LtnOmega(t)} equivalently says that 
level-wise application of $L_t$ takes $\Omega$-objects in the category $\PSpt_\Sigma(C)$ to $\Omega_t$-objects in $\PSpt_\Sigma(C)$.
The latter is equivalent to \Cref{lm:ReformulationOmega}(2) by \Cref{lm:LreptLprepcpreserveexact}.

Thus we have the equivalences of \Cref{lm:ReformulationOmega}(3), \Cref{lm:ReformulationOmega}(2), and \eqref{eq:LtnOmega(t)}, and the implication to \Cref{lm:ReformulationOmega}(1)
To conclude the claim of the lemma we show the implication from \Cref{lm:ReformulationOmega}(1) to \Cref{lm:ReformulationOmega}(2). 
We mention that
for $F\in \PSpt_{\Omega}(C)$ the terms of $\Lrep_{\Omega,t} F$ are given by
\[E_t\Omega_t^\infty L^\Omega_t\Shift^{-l} F,\]
where $\Shift^{-l}$ denote the shift on $\PSpt_\Omega(C)$.
So \Cref{lm:ReformulationOmega}(2) follows because by \Cref{lm:ReformulationOmega}(1) 
$E_t\Omega_t^\infty L^\Omega_t\Shift^{-l} F\simeq E_t L_t\Omega^\infty \Shift^{-l} F=\Lrep_t \Omega^\infty \Shift^{-l} F$,
that equals the terms of 
$\Lrep_{t} L_{\Omega} F$.
\end{proof}

\subsection{Properties of the endofunctor $\Omega_{\Gm}\Lrep_{\A^1}$.}

\newcommand{\Xlocx}{X^\mathrm{loc}_{x}} 

In this subsection, we consider the endo-functor $\Omega_{\Gm}^l\Lrep_{\A^1}$ on $\Pre^{S^1,\fr}(\Sm_{S,z})$
for a scheme $S$ as above and $z\in S$.

Firstly, we show by the standard argument used already in \cite{GW-cancel,DruDMGWeff}
that for a base field $k$ that satisfy assumption (1) from \Cref{sect:AssumptionsBase} 
there is the Nisnevich local stable weak equivalence
\[\Omega_{\Gm}F\to \Omega_{\Gm}\Lrep_\nis(F)\] for $F\in \Pre^{S^1,\fr}_{\A^1}(\Sm_k)$,
or equivalently, $\Omega_{\Gm}$ preserves Nisnevich equivalences on $F\in \Pre^{S^1,\fr}_{\A^1}(\Sm_k)$.
We prove the claim in the general form of \Cref{lm:OmegaGmsimeqOmegaGmLrepnis} and
next using the above reformulation lemmas we deduce the property of the functor $\Omega_{\Gm}^l\Lrep_{\A^1}$ on $\Pre^{S^1,\fr}(k)$ in \Cref{lm:PregpfrkOmegaLrepNisexSHI}.

\begin{lemma}\label{lm:OmegaGmsimeqOmegaGmLrepnis}
For any $F\in \Pre^{S^1,\fr}_{\A^1}(\Sm_k)$ 
and any open subscheme $V$ in $\A^1_k$
there is the Nisnevich local stable equivalence of $S^1$-spectra presheaves \begin{equation}\label{eq:OmegaGm(Lnis)}
F^V\simeq_{S^1,\nis} (\Lrep_\nis F)^V.
\end{equation}
\end{lemma}
\begin{proof}
The presheaf $F$ is $\A^1$-invariant, and by the strictly homotopy invariance theorem from point (1) of \Cref{sect:AssumptionsBase} $\Lrep_\nis(F)$ is $\A^1$-invariant too. 
Denote 
$\overline F=
\cofib(F\to \Lrep_\nis(F)
.$
Then the presheaves $\overline F$ and $\overline{F}^V$ 
are $\A^1$-invariant. 
Let $X\in \Sm_k$, $x\in X$, and let $\eta$ be the generic point of the local scheme $\Xlocx$ of $X$ at $x$, and let $\iota$ denote the generic point of $V\times\eta$. 
Then for all $l\in \mathbb Z$ the canonical morphisms
\begin{equation}\label{eq:F(VXlocxotVetaotiota)}\pi_l \overline{F}^V(X^\mathrm{loc}_{x}) \to \pi_l \overline{F}^V(\eta)\to \pi_l \overline{F}(\iota)\end{equation}
are injections of abelian groups by \cite[Theorem 3.15(1,3)]{hty-inv}. 
Since $F(\iota)\simeq L_\nis F(\iota)$, the right side group in \eqref{eq:F(VXlocxotVetaotiota)} is trivial, and consequently the left side one too.
Hence $\overline{F}^V\simeq 0$.
The claim follows because 
\[\overline{F}^V\simeq \cofib (F^V\to (\Lrep_\nis F)^V ).\]
\end{proof}

\begin{lemma}
\label{lm:PregpfrkOmegaLrepNisexSHI}
The functors $\Omega_{\Gm}^l\Lrep_{\A^1}$ 
on the categories 
$\Pre^{S^1,\fr}(k)$ and 
on the category $\Pre^{\gp\fr}(k)$
over a field $k$ as above
preserve Nisnevich equivalences.
\end{lemma}
\begin{proof}
Let $F_0\to F_1$ be a Nisnevich local equivalence in $\Pre^{S^1,\fr}(k)$. Then $F\defeq \fib (F_0\to F_1)\simeq_{\nis} 0$.
We are going to show that $\Omega_{\Gm}^l\Lrep_{\A^1}\simeq_\nis 0$. Then $\Omega_{\Gm}^l\Lrep_{\A^1} F_0\simeq_\nis \Omega_{\Gm}^l\Lrep_{\A^1}F_1$.

By strict homotopy invariance theorem over $k$ $\Lrep_\nis$ on $\Pre^{s,\fr}(\Sm_k)$ preserves $\A^1$-invariant objects.
Then by \Cref{lm:LreptLprepcpreserveexact} $\Lrep_{\A^1}$ preserves Nisnevich local equivalences;
so \begin{equation}\label{eq:LA1nisexact}\Lrep_{\A^1}F\simeq_{\nis,S^1} 0.\end{equation}
Then by \Cref{lm:OmegaGmsimeqOmegaGmLrepnis} \[\Omega_{\Gm}^l\Lrep_{\A^1} F\simeq \Omega_{\Gm}^l\Lrep_\nis\Lrep_{\A^1} F\stackrel{\eqref{eq:LA1nisexact}}{\simeq} 0.\]

The claim for $\Pre^{\gp\fr}(k)$ follows from the case of $\Pre^{S^1,\fr}(k)$ because
the delooping construction for group-like simplicial spaces gives a functor $E\colon \Pre^{\gp\fr}(\Sm_k)\to \Pre^{s,\fr}(\Sm_k)$ such that $\Omega^\infty_{S^1} E(F)\simeq F$ for each $F\in \Pre^{\gp\fr}(\Sm_k)$, and because the functor $E$ preserves Nisnevich local equivalences.
\end{proof}

Further we apply results of \Cref{section:SplittingdiagramLocalisation} to deduce the claim for $\Pre^{\gp\fr}_\tf(\Sm_{S,z})$ and $\Pre^{\gp\fr}_\tf(\Sm_S)$.

\begin{lemma}\label{lm:OmegaGmLrepA1(simeq_nis)ShtfScz}
Given a scheme $S$ as above, and a point $z\in S$,
the functor $\Omega_{\Gm}^l\Lrep_{\A^1}$ preserves Nisnevich equivalences on the category $\Pre^{\gp\fr}_\tf(\Sm_{S,z})$.
\end{lemma}
\begin{proof}
Since $\tf$-topology is trivial on $\Sm_{S,z}$, 
$\Sh^{\gp\fr}_\tf(\Sm_{S,z})\simeq\Pre^{\gp\fr}(\Sm_{S,z})$.
For each $F$ in $\Pre^{\gp\fr}(\Sm_{S,z})$
there are natural equivalences 
\begin{equation}\label{eq:OmegaGmLA1overarrowi}
\Omega_{\Gm}^l\Lrep_{\A^1} F
\stackrel{\Cref{th:LocA1tfnisstructuresDeformation}(1,2)}{\simeq}
\overarrow{i}_*\overarrow{i}^*( \Omega_{\Gm}^l\Lrep_{\A^1} F ) 
\stackrel{\Cref{th:LocA1tfnisstructuresDeformation}(4)}{\simeq}
\overarrow{i}_*(\Omega_{\Gm})^l\Lrep_{\A^1} (\overarrow{i}^*F)).
\end{equation}

Given a Nisnevich equivalence $F\to G$ in $\Pre^{\gp\fr}(\Sm_{S,z})$,
by \Cref{th:LocA1tfnisstructuresDeformation}(3) $\overarrow{i}^*F\simeq_{\nis}\overarrow{i}^*G$; 
by \Cref{lm:PregpfrkOmegaLrepNisexSHI} 
$\Omega_{\Gm}^l\Lrep_{\A^1} (\overarrow{i}^*F)\simeq_{\nis} \Omega_{\Gm}^l\Lrep_{\A^1} (\overarrow{i}^*G)$;
by \Cref{th:LocA1tfnisstructuresDeformation}(3) 
\[
\overarrow{i}_*(\Omega_{\Gm}^l\Lrep_{\A^1} (\overarrow{i}^*F))
\simeq_{\nis}
\overarrow{i}_*(\Omega_{\Gm}^l\Lrep_{\A^1} (\overarrow{i}^*G)).\]
Then by \eqref{eq:OmegaGmLA1overarrowi} the claim follows. 
\end{proof}

\subsection{Formulas in the category $\Shv_\tf^\fr(S)$.}
In this subsection, we compute the 
unit of the adjunction 
\[\Sigma^\infty_{S^1,\Gm, (\A^1,\nis)}\colon \bfShv_\tf^\fr(S) \leftrightarrows \SH_{\nis,\A^1}^{S^1,\Gm,\fr}(S)\colon \Omega^\infty_{S^1,\Gm, (\A^1,\nis)}\]
in terms of the unstable motivic localisation $\Lrep^{\fr,\tf}_{\A^1}$ on $\Pre_\tf^\fr(S)$.
Denote by unit of the adjunction of $\infty$-categories
\[\Sigma^\infty_{S^1,\Gm, (\A^1,\nis)}\colon \Pre_\tf^\fr(S) \leftrightarrows \rmSH_{\nis,\A^1}^{S^1,\Gm,\fr}(S)\colon \Omega^\infty_{S^1,\Gm, (\A^1,\nis)}\]
by 
$(\Omega\Sigma)^{\tf,\fr}_{S^1,\Gm,(\A^1,\nis)}$.

Denote by 
$(\Omega\Sigma)^{\tf,\fr}_{S^1,\Gm,(\A^1,\tf)}$
the unit of the adjunction 
\[\Sigma^\infty_{S^1,\Gm, (\A^1,\tf)}\colon \Pre_\tf^\fr(S) \leftrightarrows \rmSH_{\tf,\A^1}^{S^1,\Gm,\fr}(S)\colon \Omega^\infty_{S^1,\Gm,(\A^1,\tf)}.\]

\begin{proposition}\label{prop:Pretffr(SmSz)SGA1nis}
There is the natural equivalence
\[(\Omega\Sigma)^{\tf,\fr}_{S^1,\Gm,(\A^1,\nis)} \simeq \Lrep_{\nis}(\Omega\Sigma)^{\tf,\fr}_{S^1,\Gm,(\A^1,\tf)} \]
on $\Pre_\tf^\fr(\Sm_{S,z})$.
\end{proposition}
\begin{proof}
By the definition and \Cref{th:restrfrLocniscommute} 
$\Pre_{\nis,\A^1}^\fr(\Sm_{S,z})$, and
$\rmSH_{\nis,\A^1}^{S^1,\Gm,\fr}(\Sm_{S,z})$
are
the
Nisnevich localisations of 
$\Pre_{\tf,\A^1}^\fr(\Sm_{S,z})$, and
$\rmSH_{\tf,\A^1}^{S^1,\Gm,\fr}(\Sm_{S,z})$.
By \Cref{lm:LocsubsetStab} 
\[\rmSH_{\tf,\A^1}^{S^1,\Gm,\fr}(\Sm_{S,z})\simeq \Pre_{\tf,\A^1}^\fr(\Sm_{S,z})[(\Omega_{\Gm}\Lrep_{\A^1})^{-1}].\]
Then since
by \Cref{lm:OmegaGmLrepA1(simeq_nis)ShtfScz} $\Omega_{\Gm}\Lrep_{\A^1}$ on $\Pre_\tf^\fr(\Sm_{S,z})$ preserves Nisnevich local equivalences
the claim follows.
\end{proof}

\begin{lemma}\label{lm:CommuteSSz}
For any scheme $S$ and a point $z\in S$
the functor $\Pre^\fr(\Sm_S)\to \Pre^\fr(\Sm_{S_z})$ commutes with $\Lrep_{\A^1}$, $\Lrep_{\tf}$, $\Lrep_{\nis}$, $(\Omega\Sigma)^{\tf,\fr}_{S^1,\Gm,(\A^1,\tf)}$, $(\Omega\Sigma)^{\tf,\fr}_{S^1,\Gm,(\A^1,\nis)}$.
\end{lemma}
\begin{proof}
The claim holds because the base change commutes with the endo-functors $F\mapsto F^V$ for any $V\in \Sm_S$, and with the localisation with respect to any cd-topology localisation.
\end{proof}

\begin{proposition}\label{prop:Pretffr(S)SGA1nis}
There is the natural equivalence
\[(\Omega\Sigma)^{\tf,\fr}_{S^1,\Gm,(\A^1,\nis)} \simeq \Lrep_{\nis}(\Omega\Sigma)^{\tf,\fr}_{S^1,\Gm,(\A^1,\tf)} \]
on $\Pre_\tf^\fr(S)$.
\end{proposition}
\begin{proof}
We prove the claim by induction on $\dim S$. Assume that claim is proven for all schemes of dimension less then $\dim S$.

By \Cref{lm:CommuteSSz} without loss of generality we can assume that $S$ is local.
Let $z\in S$ be the closed point.
Then by \Cref{prop:Pretffr(SmSz)SGA1nis}
the equivalence holds on $\Pre_\tf^\fr(\Sm_{S,z})$.
By the inductive assumption 
the equivalence holds on $\Pre_\tf^\fr(\Sm_{S-z})$.
By \Cref{th:LocA1tfnisstructuresShLocsquare}(2,3) 
it follows that 
\[\tilde i_*\tilde i^!(\Omega\Sigma)^{\tf,\fr}_{S^1,\Gm,(\A^1,\nis)}\simeq\tilde i_*\tilde i^!\Lrep_{\nis}(\Omega\Sigma)^{\tf,\fr}_{S^1,\Gm,(\A^1,\tf)},
j_*j^*(\Omega\Sigma)^{\tf,\fr}_{S^1,\Gm,(\A^1,\nis)}\simeq j_*j^*\Lrep_{\nis}(\Omega\Sigma)^{\tf,\fr}_{S^1,\Gm,(\A^1,\tf)}\]
on the category $\Spt^{\Gm,\fr}_\tf(\Sm_{S})$.
The claim follows by \Cref{th:LocA1tfnisstructuresShLocsquare}(1).
\end{proof}

We complete the section with the formulas for $(\Omega\Sigma)^{\tf,\fr}_{S^1,\Gm,(\A^1,\tf)}(F)$.

\begin{proposition}\label{prop:OmegaSigmafrtf_SGA1}
For any 
$F\in \Pre_\tf^\fr(S)$
there is the natural equivalence
\[(\Omega\Sigma)^{\tf,\fr}_{S^1,\Gm,(\A^1,\tf)}(F)\simeq \varinjlim_{l}\Omega^l_{\Gm}\Lrep^{\fr,\tf}_{\A^1}(\Sigma^l_{\Gm}F^\gp).\]
\end{proposition}
\begin{proof}
The claim follows because the functors $\Omega^l_{\Gm}$ preserve the subcategory $\Pre^\fr_{\A^1,\tf}(S)$ in $\Pre_\tf(S)$.
\end{proof}

Consider the canonical functor $\Pre^\fr_{\tf}(\mathcal S)\to \Pre^\fr(\mathcal S)$.
\begin{lemma}\label{lm:LreptffrA1LrepA1tf}
There is the canonical equivalence
\[\gamma_*\Lrep^{\fr,\tf}_{\A^1}\simeq \Lrep_{\A^1,\tf} \gamma_*.\]
\end{lemma}
\begin{proof}
Since $\Lrep^\fr_{\A^1}(F)=F(-\times\Delta^\bullet_S)$, $\Lrep_{\A^1}(F)=F(-\times\Delta^\bullet_S)$, there is the equivalence $\gamma_*\Lrep^{\fr}_{\A^1}\simeq \Lrep_{\A^1} \gamma_*\colon \Pre^\fr(S)\to \Pre_{\A^1}(S)$.
The claim follows because of \Cref{th:restrfrLocniscommute}.
\end{proof}

\subsection{Formulas in $\Pre^\fr_\Sigma(S)$.}

\begin{theorem}\label{th:oopspacesGmS1andT}
For $X\in \Sm_S$, and an open subscheme $U$ the $\infty$-loop space in the $\infty$-category of Nisnevich sheaves $\Pre_\nis(\Sm_S)$ of the motivic suspension spectrum of $X$ in $\Spt^{\PP^1}_{\A^1,\nis}(S)$
is Nisnevich locally equivalent to 
\begin{equation}\label{eq:OmegastOmegaPinftyLmot(Y)}\begin{array}{lllll}
\varinjlim_{l}&\Omega^l_{\Gm^{\wedge 1}}&(\Lrep_{\A^1,\tf}&\Fr((X/U)\wedge \Gm^{\wedge l}))^\gp&\simeq\\
\varinjlim_{l}&\Omega^l_{\Gm^{\wedge 1}\wedge S^1}&\Lrep_{\A^1,\tf}&\Fr((X/U)\wedge T^{\wedge l}).
\end{array}\end{equation}
\end{theorem}
\begin{proof}
The general case follows from the case of $U=\emptyset$ in view of \Cref{th:Fr(Y/V)}.
The first equivalence follows by \Cref{prop:Pretffr(S)SGA1nis,prop:OmegaSigmafrtf_SGA1,lm:LreptffrA1LrepA1tf}.
The second follows from the first in view of \Cref{th:Fr(Y/V)}
applied to pairs $(X\times\A^n)/(X\times\A^n-Y\times 0)$, where $Y=X\setminus U$.
\end{proof}

\begin{theorem}\label{th:loopspaceGmSAST}
Let $S$ be as in \Cref{sect:AssumptionsBase}.
For $X\in \Sm_S$, and an open subscheme $U$ the $\infty$-loop space in the $\infty$-category of Nisnevich sheaves $\Pre_\nis(\Sm_S)$ of the motivic suspension spectrum of $X$ in $\Spt^{\PP^1}_{\A^1,\nis}(S)$
is Nisnevich locally equivalent to 
\begin{equation}\label{eq:OmegaGSAOmegaPinftyLmot(Y)}\begin{array}{lllll}
\varinjlim_{l}&\Omega^l_{\Gm^{\wedge 1}\wedge \Delta^1_S/\delta\Delta^1_S}&\Lrep_{\tf}&\Fr((X/U)\wedge T^{\wedge l})&\in \sShv_\bullet(\Sm_S).
\end{array}\end{equation} 
\end{theorem}
\begin{proof}

We use the argument form \cite{Morel-connectivity}.

Consider the category $\Pre^{S^1,\fr}_{\tf}(B)$. Let $L^{{S^1,\tf}(1)}_{\A^1}F=\cofib(F(\A^1\times-)\rightrightarrows F)$, where the morphisms are given by the zero and the unit sections.
The canonical morphism $L^{{S^1,\tf}(1)}_{\A^1}F\leftarrow F$ is an equivalence, if and only if $F$ is $\A^1$-invariant, and 
for \begin{equation}\label{eq:LA1cubeinftylimit}(L^{{S^1,\tf}(1)}_{\A^1})^\infty = \varinjlim_l (L^{{S^1,\tf}(1)}_{\A^1})^l \end{equation}
$L^{{S^1,\tf}(1)}_{\A^1}(L^{{S^1,\tf}(1)}_{\A^1})^\infty F\simeq (L^{{S^1,\tf}(1)}_{\A^1})^\infty F$. Hence \eqref{eq:LA1cubeinftylimit} is equivalent to the $\A^1$-localisation functor $L^{S^1,\tf}_{\A^1}$ on $\Pre^{S^1,\fr}_{\tf}(B)$.

Denote by $S^1_{\A^1}=\cofib(\{0,1\}\to \A^1)$, and then $\Omega_{S^1_{\A^1}}=\fib(F(\A^1\times-)\rightrightarrows F)$.
Because of the equivalence 
\begin{equation}\label{eq:OmegaS1A}\Omega_{S^1_{\A^1}}F[1]=\Omega_{S^1_{\A^1}}(\Sigma_{S^1}F)\simeq L^{{S^1,\tf}(1)}_{\A^1}F,\end{equation}
where $F[1]=F\wedge S^1$,
we can say that the limit of powers of \eqref{eq:OmegaS1A} 
\[\varinjlim_l \Omega_{S^1_{\A^1}}^lF[l]\]
is equivalent to $L^{{S^1,\tf}}_{\A^1}$. 
Consider the category $\Pre^{S^1,\fr}(B)$, and the $\A^1$-$\tf$-motivic localisation $L^{S^1}_{\A^1,\tf}$.
Then we can write that $L^{S^1}_{\A^1,\tf}=L^{S^1,\tf}_{\A^1}L^{S^1}_{\tf}$, and since the functor $\Omega_{S^1_{\A^1}}$ on $\Pre^{S^1,\fr}(B)$ preserves $\tf$-local objects, it follows that \[L^{S^1}_{\A^1,\tf}\simeq \varinjlim_l \Omega_{S^1_{\A^1}}^l L^{S^1}_{\tf}F[l]\] on $\Pre^{S^1,\fr}(B)$. Moreover, since the functor
\[F\mapsto \varinjlim_l \Omega_{S^1_{\A^1}}^l L^{S^1}_{\tf}F[l] \]
on the category $\Pre^{\fr}(B)$ commutes with $\Omega_{S^1}$ it preserves $\Omega_{S^1}$-spectra,
and for any $\Omega_{S^1}$-spectrum $F$ the $S^1$-spectrum of framed presheaves $\varinjlim_l \Omega_{S^1_{\A^1}}^l L^{S^1}_{\tf}F[l]$ is $S^1$-stable $\A^1$-$\tf$-motivically local.

Then we have the sequence of equivalences
\[\begin{array}{lcl}
\varinjlim_{l}\Omega^{l}_{\Gm}\Omega_{S^1}^\infty \Lrep_{\A^1,\tf} \Fr(\Sigma^\infty_{\Gm}(X/U))^\gp & \simeq&
\varinjlim_{l_1}\varinjlim_{l_2}\Omega^{l_1}_{\Gm}\Omega^{l_2}_{S^1_{\A^1}} \Lrep_{\tf} \Fr(\Sigma^{l_1}_{\Gm}\Sigma^{l_2}_{S^1}(X/U))^\gp\\ &\simeq&
\varinjlim_{l_1}\varinjlim_{l_2}\Omega^{l_1}_{\Gm}\Omega^{l_2}_{S^1_{\A^1}} \Lrep_{\tf} \Fr(\Sigma^{l_1}_{\Gm}\Sigma^{l_2}_{S^1}(X/U))\\ &\simeq& 
\varinjlim_{l}\Omega^{l}_{\Gm}\Omega^{l}_{S^1_{\A^1}} \Lrep_{\tf} \Fr(\Sigma^{l}_{\Gm}\Sigma^{l}_{S^1}(X/U))\\ &\simeq& \varinjlim_{l}\Omega^{l}_{\Gm\wedge S^1_{\A^1}} \Lrep_{\tf} \Fr(\Sigma^{l}_{T}(X/U))
\end{array}\]
\end{proof}

\begin{corollary}\label{cor:generatorsand_relations}
Let $S$ be 
local  and $\dim S=1$.
Let $z$ denote the closed point
Then for any $Y\in \Sm_S$ there is the isomorphism of Nisnevich sheaves
\begin{equation}\label{eq:pi-1(Y)conginjlimH1tf(GmSB^l,ZF(YTl))}\pi_{-1}(Y)\cong \varinjlim_{l} H_\tf^1( (\Gm^{\wedge 1}\wedge S_{S}^1)^{\wedge l}\times -, \ZF(Y\wedge T^{\wedge l})_\tf ). 
\end{equation}
\end{corollary}
\begin{proof}
The formula 
from \Cref{th:loopspaceGmSAST} 
provides that $\pi_{-l}(Y)=0$ for $l>0$, 
and consequently by the topological Hurewicz theorem $\pi_{-1}(Y)$ equals $\pi_{-1}(\mathrm{H}\mathbb{Z}\wedge Y)$.
Then the formula from \Cref{th:loopspaceGmSAST} yields the isomorphism \eqref{eq:pi-1(Y)conginjlimH1tf(GmSB^l,ZF(YTl))}%
.
\end{proof}
\begin{remark}Note that for any presheaf $F$ there is the isomorphism
\[ H_\tf^1(U,F_\tf) \cong \coker \left( F( U\times_S(S-z) )\oplus F( U^h_z ) \to F( U^h_{z}\times_S(S-z) ) \right) .
\]
The latter result allows to prove that the sheaf $\pi_{-1}$ is non-trivial for smooth $Y$ in general over a one-dimensional scheme $S$.
We left this for future work.
\end{remark}

\begin{theorem}\label{th:LsmotLnisLtfsmotFPretfr}
Let
\[\calF\in \PSpt^{\Gm,S^1}_{\Sigma}(\Corr^\tfr(S)),\text{ or }\calF\in \PSpt^{T}_{\Sigma}(\Corr^\tfr(S)),\] 
see \Cref{section:StableLocalisationDefinitionsNotations} for the definition of the latter categories,
there are equivalences
\[\begin{array}{lcl}
\Lrep^\fr_\smot(F)&\simeq&\Lrep^\fr_\nis \Lrep^\fr_{\Gm}\Lrep^\fr_{\A^1,\tf}(F)^\gp,\\
\Lrep^\fr_\smot(F)&\simeq&\Lrep^\fr_\nis \Lrep^\fr_{T}\Lrep^\fr_{\A^1,\tf}(F).
\end{array}\]
\end{theorem}
\begin{proof}
Note that by \Cref{lm:LreptffrA1LrepA1tf} and \Cref{th:restrfrLocniscommute} and since the functor $F\mapsto F((-)_+\wedge \Gm)$ commutes with $\gamma_*$ the superscript $\fr$ in the above formulas could be skipped with replacing of $F$ by $\gamma_*F$.
By \Cref{th:oopspacesGmS1andT} 
it follows that 
$\Lrep_\smot(F)\simeq\Lrep^\fr_\nis \Lrep^\fr_{\A^1,\tf,\Gm,S^1}(F)$
because the terms of the spectrum 
are equivalent to $({\Gm,S^1})$-$\infty$-loop space of $S^{l_1}\wedge \Gm^{\wedge l_2}\wedge Y$.
The first equivalence in the theorem follows because of \eqref{prop:OmegaSigmafrtf_SGA1} and because of the equivalence
\begin{equation}\label{eq:gphfrS1}
\Lrep_{S^1}(h^\tfr(\Sigma^\infty_{S^1}(X/U)))\simeq h^\tfr(\Sigma^\infty_{S^1}(X/U))^\gp.
\end{equation}
The case of $T$-spectra follows because
the $\A^1$-homotopy equivalence $T=\A^1_S/(\A^1_S-0)\simeq S^1\wedge \Gm$,
and consequently the equivalences of categories
$\Pre^{T}_{\A^1}(\Sm_S)\simeq \Pre^{S^1\wedge \Gm}_{\A^1}(\Sm_S)$,
$\Pre^{T}_{\A^1}(\Sm_S)\simeq \Pre^{S^1\wedge \Gm}_{\A^1}(\Sm_S)$.
\end{proof}

\begin{corollary}\label{th:SHnisA1S1GmSHfrtfA1S1Gm}
Let $S$ be a base scheme satisfying assumptions from \Cref{sect:AssumptionsBase}, for example, $\Spec \mathbb Z$.

There are reflective embeddings
\[\SH_{\A^1,\nis}(S)\to \SH_{\A^1,\tf}(\Corr^\tfr_S), \quad\SH^\mathrm{veff}_{\A^1,\nis}(S) \to \SH^\mathrm{veff}_{\A^1,\tf}(\Corr^\tfr_S)\]
with images spanned by Nisnevich sheaves in the right-side categories; the left adjoint functors preserve the motives of smooth schemes,
and for a local henselian essentially smooth $U$, and any $Y\in \Sm_S$ \[\pi_{*,*}^{\SH_{\A^1,\nis}(S)}(Y)(U)\simeq\pi_{*,*}^{\SH^\fr_{\A^1,\tf}(S)}(Y)(U).\]
where $\SH^{\fr}_{\A^1,\tf}(S)=\SH_{\A^1,\tf}(\Corr^\tfr_S)$
\end{corollary}

\begin{theorem}\label{LsmotLnisLtfsmotF}
Let $\Lrep_{\nissmot}$ be the the stable motivic localisation endofunctor on the categories of bi-spectra of $T$-spectra of simplicial presehaves.
For any
radditive bi-spectrum or $T$-spectrum of quasi-stable framed presheaves $\calF$ in sense of \cite{Framed}, 
there is the canonical scheme-wise level-wise equivalence
\begin{equation*}\label{eq:Lnismot(F)simeqLnisLtfsmot(F)}\Lrep_{\nissmot}(\calF)\simeq \Lrep_\nis\Lrep_{\tfsmot}(\calF),\end{equation*}
where $\Lrep_{\tfsmot}$ is the stable $\tf$-motivic localisation with respect to $\tf$-topology instead of the Nisnevich topology,
and $\Lrep_{\nis}$ is equivalent to the level-wise Nisnevich localisation, and consequently, it preserves Nisnevich locally connective objects.
In particular, for $X\in \Sm_S$, and an open subscheme $U$
\[\Lrep_\mot(\Sigma_{\Gm,S^1}^\infty )=\Lrep_\nis \Lrep_{\A^1,\tf,\Gm}(\Fr(\Sigma^\infty_{\Gm,S^1}X/U))^\gp \simeq \Lrep_\nis \Lrep_{\A^1,\tf,\Gm}(h^\fr(\Sigma^\infty_{\Gm,S^1}X/U)^\gp),\]
and
\[\Lrep_\mot(\Sigma_{T}^\infty )=\Lrep_\nis \Lrep_{\A^1,T}(\Fr(\Sigma^\infty_{T}(X/U))^\gp \simeq \Lrep_\nis \Lrep_{\A^1,\tf,\Gm}(h^\fr(\Sigma^\infty_{\Gm,S^1}X/U)^\gp),\]
where $\gp$ denotes the group-completion of $\infty$-monoid or $S^1$-spectrum.
\end{theorem}
\begin{proof}
By \cite[Corollary 2.3.25]{five-authors} there is the scheme-wise equivalence of presheaves 
$\Lrep_{\A^1}h^\tfr(Y)\big|_{\SmAff_S}\simeq \Lrep_{\A^1}\Fr(-,Y)\big|_{\SmAff_S}$. 
Hence the subcategory of 
$\A^1$-invariant radditive stable framed simplicial presheaves 
in the category of simplicial presheaves on $\Sm_S$
is equivalent to the image of $\Pre_{\A^1}(\Corr^\tfr(S))$.
Hence the first claim follows from \Cref{th:LsmotLnisLtfsmotFPretfr}.
The second claim 
follows from the first because of \Cref{th:Fr(Y/V)}
and 
because of \eqref{eq:gphfrS1} and
\[
\Lrep_{S^1}(\Fr(-\Delta^\bullet_S, \Sigma^\infty_{S^1}(X/U)))\simeq \Fr(-\Delta^\bullet_S, \Sigma^\infty_{S^1}(X/U))^\gp,
\]
by \cite[Theorem 6.5]{Framed}, and see also \cite[Appendix B]{SmAffOpPairs}.
\end{proof}

\begin{theorem}
Let $S$ be a scheme that satisfy assumptions (1) and (2) from \Cref{sect:AssumptionsBase}.
The functor 
\[\SH_{\tf,\A^1}^\fr(S)\to \SH_{\nis,\A^1}^\fr(S)\]
preserves homotopy t-structure.
\end{theorem}
\begin{proof}
By \Cref{th:LsmotLnisLtfsmotFPretfr} the localisation functor in induced by 
the functor $\Lrep^\fr_\nis\colon \Pre^\fr(S)\to \Pre^\fr(S)$ 
in view of the embeddings of the $\infty$-categories
$\rmSH_{\tf,\A^1}^{S^1,\Gm,\fr}(S)$ and $\rmSH_{\nis,\A^1}^{S^1,\Gm,\fr}(S)$
into the category $\Pre^{S^1,\Gm,\fr}(S)$.
So the claim follows since $\Lrep_\nis$ takes connective objects to Nisnevich locally connective objects.
\end{proof}

\appendix
\section{Stable localisation}\label{section:StableLocalisationDefinitionsNotations}

Let $P$ be an $\infty$-category, and \[\Sigma \colon P\leftrightarrows P\colon\Omega\] is an adjunction. 

Let $\mathbb Z_{\geq 0}^\mathrm{op}$ denote the category with objects non-negative integers and unique morphism $l_0\leftarrow l_1$ for each pair $l_0\leq l_1$.
Consider the functor $\mathbb Z_{\geq 0}^\mathrm{op}\to \mathrm{Cat}_\infty$ given by the sequence of functors
\begin{equation}\label{eq:seqPre(C)Omega}
P\xleftarrow{\Omega}P\xleftarrow{\Omega}P\xleftarrow{\Omega}\dots.\end{equation}
Denote by $P_\Omega$
the $\infty$-category fibred over $\mathbb Z_{\geq 0}^\mathrm{op}$ defined by \eqref{eq:seqPre(C)Omega} in view of the Grothedieck correspondence. 
Let $P_\Omega\colon F\mapsto F[-1]$ denote the functor that 
shifts the sequence \eqref{eq:seqPre(C)Omega} by one to the left.

\begin{definition}
Put $\PSpt^\Omega(P)= \Func_{\mathbb Z_{\geq 0}^\mathrm{op}}(\mathbb Z_{\geq 0}^\mathrm{op}, P_\Omega)$, that is the category of sections of the fibred category $P_\Omega$ over $\mathbb Z_{\geq 0}^\mathrm{op}$.
\end{definition}
For $F\in \PSpt^\Omega(P)$ denote by $F[-1]\in \PSpt^\Omega(P)$ the object given by
the composite of the sequence
\[\mathbb Z_{\geq 0}^\mathrm{op}\simeq \mathbb Z_{\geq 1}^\mathrm{op}\hookrightarrow \mathbb Z_{\geq 0}^\mathrm{op}\xrightarrow{F}\Pre_\Omega
.\]
Denote by $\Omega\colon P_\Omega\to P_\Omega$ the functor induced by $\Omega$ on $P$.

\begin{definition}
Denote by 
$\nu\colon \Id_{\PSpt^\Omega(P)}\to \Omega(-)[-1]$
the canonical natural transformation.
Denote by 
$\Spt^\Omega(P)$ the subcategory of $\PSpt^\Omega(P)$
spanned by objects $F$ such that 
$F\simeq \Omega F[-1]$.
\end{definition}
Denote $\Omega^l(-)[l]\colon \PSpt^\Omega(P)\to \PSpt^\Omega(P)$ the $l$-th power of the functor $\Omega(-)[1]$.

\begin{lemma}\label{lm:SptS(C)}
(1) The subcategory $\Spt^\Omega(P)$ is a reflective subcategory of $\PSpt^\Omega(P)$,
and 
the localisation functor
\[\PSpt^\Omega(P)\to \Spt^\Omega(P)\]
takes $F$ to $\varinjlim_{l} \Omega^l F[l]$;
(2) 
there is a natural equivalence 
$P[\Omega^{-1}]\simeq \Spt^\Omega(P)$.
\end{lemma}
\begin{proof}
(1) 
Consider the canonical natural transformation $F\to \varinjlim_{l} \Omega^l F[l]$.
The first claim follows, because 
the functor 
$F\mapsto \varinjlim_{l} \Omega^l F[l]$
is idempotent, since \[\Omega(\varinjlim_{l} \Omega^l F[l])[1]\cong \varinjlim_{l} \Omega^{l+1} F[l+1]\cong \varinjlim_{l} \Omega^l F[l],\]
and takes the above natural transformation to the identity.

(2) The second claim follows because of equivalences
\begin{equation*}\label{eq:POmegastablimseqlimarrowSpt}P[\Omega^{-1}]\simeq 
\operatorname{lim}\left(
P\xleftarrow{\Omega}P\xleftarrow{\Omega}P\xleftarrow{\Omega}\dots
\right)\simeq 
\operatorname{eq}\left(
\xymatrix{\prod\limits_{l\geq 0}P& \ar@<0.5ex>[l]^{\mathrm{shift}\circ \Omega}\ar@<-0.5ex>[l]_{\mathrm{shift}}\prod\limits_{l>0}P}
\right)\simeq
\Spt^\Omega(P)
,\end{equation*}
where $\mathrm{shift}\colon \prod_{l>0}P\simeq \prod_{l\geq 0}P$.
\end{proof}

Denote by  
\[L_S\colon P\to P[\Omega^{-1}]\] the localisation functor provided by \Cref{lm:SptS(C)},
and by $\mathcal L_S\colon \PSpt^\Omega(P)\to \PSpt^\Omega(P)$ the composite of $L_S$ and the canonical embedding.


\begin{definition}\label{den:PreGmT(SmS)}
Suppose $S\to C$ be a morphism of $\infty$-categories equipped with the canonical functor $S\times C\to C$.
Let $P=\Pre(C)$, $\alpha\colon V\to U$ be a morphism in $S$, and $\Omega\colon P\to P$ be the functor 
\[\Omega(F)(X) = \fib(F(X\times U)\to F(X\times V))\] for $F\in P$, $X\in \Sm_S$.
Denote 
\[\Pre^{V/U}(C)=\PSpt^\Omega(\Pre(C)), \Pre^{V/U}_{V/U}(C)=\Spt^\Omega(\Pre(C)).\]
\end{definition}

\section{Trivial fibre and Nisnevich localisation of framed presheaves.}

Consider the adjunction
\[\gamma^*\colon \Pre(\mathcal S)\leftrightarrows \Pre^\fr(\mathcal S)\colon \gamma_*\]
for $\mathcal S$ being $\Sm_S$, or $\SmAff_S$, $\Smat_S$.
The Nisnevich and trivial fibre localisations on categories $\Pre(S)$ and $\Pre^\fr(S)$ agree in view of the forgetful functor $\gamma^*$.
The case of Nisnevich topology is proven in \cite{five-authors}, and the case of trivial fibre topology follows by the same argument.

\begin{definition}\label{def:PrefrNislocObj}
An object $F\in \Pre^\fr(\mathcal S)$ is $\tf$-local or Nisnevich local if $\gamma_*F$ is Nisnevich local.
\end{definition}
\begin{definition}\label{def:PrefrNislocsimeqNis}
Define the class of \emph{Nisnevich local equivalences} in $\Pre^\fr(\mathcal S)$ as the image of the Nisnevich local equivalences along the functor $\gamma^*$ given by the left Kan extension.
\end{definition}
\begin{lemma}\label{lm:gamma-d-*Nisloc}
An object $F\in \Pre^\fr(\mathcal S)$ is $\tf$-local or Nisnevich local if and only if 
for any $\tf$-local or Nisnevich local equivalence $\mathcal X\to \mathcal Y\in \Pre^\fr(\mathcal S)$,
the induced morphism on Hom-spaces $[\mathcal Y,F]\to[\mathcal X,F]$ is an equivalence.
\end{lemma}
\begin{proof}
The claim follows because of the adjunction $\gamma^*\dashv\gamma_*$ and \Cref{def:PrefrNislocsimeqNis}.
\end{proof}
The subcategories of $\tf$-local or Nisnevich local objects $\Pre^\fr_\tf(\mathcal S)$ and $\Pre^\fr_\nis(\mathcal S)$ in $\Pre^\fr(\mathcal S)$ are reflective subcategories, 
denote by $L^\fr_\tf(\mathcal S)$, $\Lrep^\fr_\tf(\mathcal S)$, $L^\fr_\nis(\mathcal S)$, $\Lrep^\fr_\nis(\mathcal S)$ the localisation functors and endo-functors. Then by \Cref{lm:gamma-d-*Nisloc} $\Pre^\fr_\tf(\mathcal S)$ and $\Pre^\fr_\nis(\mathcal S)$ are the localisations of $\Pre^\fr(\mathcal S)$ with respect to the classes of $\tf$-local or Nisnevich local equivalence.

For a morphism $\widetilde X\to X$, denote by 
$h_{\widetilde X}(X)$ the subpresheaf of the representbale presheaf $h(X)$ in $\Pre(\mathcal S)$ that is the image of the canonical morphism $h(\widetilde X)\to h(X)$, and by $h_{\widetilde X}^\fr(X)$ the image of $h^\fr(\widetilde X)$ in $h^\fr(X)$.

\begin{lemma}\label{lm:gammadstfloc}
For a $\tf$-covering $\widetilde X\to X$,
the canonical morphism $h^\fr_{\widetilde X}(X)\to h^\fr(X)$ 
goes to a $\tf$-local equivalence in $\Pre(\mathcal S)$ along the functor $\gamma_*$.
\end{lemma}
\begin{proof}
We have to check the equivalence $h^\fr_{\widetilde X}(X)(U^h_z)\to h^\fr(X)(U^h_z)$ for $U\in \Sm_S$, $z\in S$.
The claim is provided by \cite[Lemma A.11]{DKO:SHISpecZ}.
\end{proof}
\begin{theorem}\label{th:restrfrLocniscommute}
There are canonical equivalences of functors from $\Pre^\fr(\mathcal S)\to \Pre(\mathcal S)$: $\gamma_* \Lrep^\fr_\tf\simeq \Lrep_\tf \gamma_*$, $\gamma_* \Lrep^\fr_\nis\simeq \Lrep_\nis \gamma_*$.
\end{theorem}
\begin{proof}
The case of $\Lrep_\nis$ is proven in \cite{five-authors}.
The case $\tf$ is similar, let us explain.

The $\infty$-category $\Pre^\fr_\tf(\mathcal S)$ by the definition is the subcategory of $\Pre^\fr(\mathcal S)$ spanned by the objects that image along $\gamma_*$ is in $\Pre_\tf(\mathcal S)$.
Hence the functor $\gamma_*$ preserves trivial fibre local objects,
and the class of the trivial fibre local equivalences $w_\tf^\fr$ in $\Pre^\fr(\mathcal S)$ is generated by the image of the one in $\Pre(S)$ along  $\gamma^*$.
So $w_\tf^\fr$ is generated by the morphisms of the form $h^\fr_{\widetilde X}(X)\to h^\fr(X))$ for $\tf$-coverings $\widetilde X\to X$.
By \eqref{lm:gammadstfloc} $\gamma_*(w^\fr_\tf)\subset w_{\tf}$, where $w_\tf$ denotes the class of trivial fibre local equivalences in $\Pre(\mathcal S)$.
Thus $\gamma_*$ preserves trivial fibre local objects and equivalences; hence it commutes with the trivial fibre localisation.
\end{proof}

\begin{corollary}\label{cor:gamma-d-sNiscons}
The functor $\gamma_*\colon \Pre^\fr(\Sm_S)\to \Pre(\Sm_S)$ 
is exact and conservative with respect to $\tf$-local and Nisnevich local equivalences.
\end{corollary}
\begin{proof}
Let $F\to G$ be a $\tf$-local or Nisnevich local equivalence in $\Pre^\fr(\Sm_S)$.
So we write $F\simeq_\tau G$, where $\tau = \tf,\nis$.
Then $L_\tau(F)\simeq L_\tau(G)$, and by \Cref{th:restrfrLocniscommute} $L_\tau\gamma_*(F)\simeq L_\tau\gamma_*(G)$, hence $\gamma_*(F)\simeq_\tau \gamma_*(G)$.

Let $F\to G$ be a morphism in $\Pre^\fr(\Sm_S)$ such that $\gamma_*(F)\simeq_\tf\gamma_*(G)$.
Then $\Lrep_\tf(\gamma_*F)\simeq \Lrep_\tf(\gamma_*G)$. By \Cref{th:restrfrLocniscommute}
$\gamma_*(\Lrep^\fr_\tf F)\simeq \gamma_*(\Lrep^\fr_\nis G)$.
Since $\gamma$ is conservative, $\Lrep^\fr_\tf F\simeq \Lrep^\fr_\tf G$.
Thus $F\simeq_\tf G$.

The case of $\Lrep_\nis$ is similar.
\end{proof}

\end{document}